\documentclass[10pt]{article}
\usepackage{graphicx}
\usepackage{latexsym}
\newtheorem{lemma}{Lemma}
\newtheorem{openquestion}{Open question}
\def\qed{$\Box$}

\newtheorem{theorem}{Theorem}

\title{Faultfree Tromino Tilings of Rectangles}
\author{Mridul Aanjaneya\footnote{Department of Computer Science and 
Engineering, IIT Kharagpur, 721302, India [Email: mridul.aanjaneya@gmail.com]} \and 
Sudebkumar Prasant Pal\footnote{Department of Computer Science and 
Engineering, and Centre for Theoretical Studies, 
IIT Kharagpur, 721302, India [Email: spp@cse.iitkgp.ernet.in]}} 
\date{}

\begin{document}
\maketitle
\begin{abstract}
In this paper we consider faultfree tromino tilings of rectangles and characterize
rectangles that admit such tilings. We introduce the notion of {\it crossing numbers}
for tilings and derive bounds on the 
crossing numbers of faultfree tilings. We develop an iterative scheme for generating
faultfree tromino tilings for rectangles and derive the closed form expression for the 
exact number of faultfree tromino tilings for $4\times3t$ rectangles and the exact
generating function for 
$5\times 3t$ rectangles, $t\geq 1$. Our iterative scheme generalizes to arbitrary rectangles; 
for $6\times 6t$ and $7\times 6t$ rectangles, $t\geq 1$, we derive
generating functions for estimating lower bounds on the number of faultfree tilings. 
We also derive an upper bound on the number of tromino tilings of an $m\times n$ rectangle, 
where $3|mn$ and $m,n>0$.  
\end{abstract}

\section{Introduction}
\label{intro}

Tilings of the plane are of interest both to statistical physicists and recreational 
mathematicians. Solomon Golomb, in a 1953 talk at the Harvard Math Club, defined a 
class of geometric figures called {\it polyominoes}, namely, connected figures formed 
of congruent squares placed so each square shares one side with at least one other square. 
Dominoes, which use two squares, and tetrominoes (the {\it Tetris} pieces), which use 
four squares, are well known to game players. Golomb first published a paper about 
polyominoes in {\it The American Mathematical Monthly} \cite{gol}. Later, Martin Gardner 
popularized polyominoes in his {\it Scientific American} columns called ``Mathematical 
Games" (see, for example, \cite{mar1}, \cite{mar2}). 

Tiling regular shapes such as rectangles using polyomino tiles is an interesting field 
of research. Stanley and Ardila \cite{ardila}, and Do \cite{math} mention several results 
about tiling regular 
shapes with polyominoes. Many of the initial questions asked about polyominoes concern the 
number of {\it n-ominoes} (those formed from $n$ squares), and what shapes can be tiled 
using just one of the polyominoes, possibly leaving one or two squares uncovered. A rectangle 
from which one square has been removed is called a {\it deficient rectangle}. Golomb \cite{gol} 
proved that deficient squares whose side length is a power of two can be tiled 
using L-shaped tromino tiles. 
Chu and 
Johnsonbaugh first extended Golomb's work to the general cases of deficient squares \cite{john}. 
They later went on to rectangles and proved a slightly weaker version \cite{chu} of the 
{\it Deficient Rectangle Theorem}, which was proved formally by Ash and Golomb in \cite{ash}. 
A rectangle from which a domino has been removed is called a {\it domino deficient rectangle}. 
Aanjaneya \cite{mridul} characterized domino deficient rectangles that admit tilings with 
L-shaped trominoes. The first significant result on tiling enumeration was obtained 
independently in 1961 by Fisher and Temperley \cite{fish} and by Kasteleyn \cite{kas}. They 
found that the number of tilings of a $2m\times 2n$ rectangle with $2mn$ dominoes is equal to

\begin{eqnarray*}
4^{mn}\prod_{j=1}^{m}\prod_{k=1}^{n}\{\cos^2\frac{j\pi}{2m+1} + \cos^2\frac{k\pi}{2n+1}\}
\end{eqnarray*}

Given a tiling, a line which cuts the rectangle into two pieces and yet does not pass through 
any of the tiles is called a {\it fault line}. A tiling of a rectangle which has no fault lines 
is called a {\it faultfree tiling}. In {\it The Mathematical Gardner} \cite{graham}, 
a collection of essays on recreational mathematics in honour of Martin Gardner, Ron Graham considers 
the problem of tiling an $m\times n$ rectangle with $a\times b$ rectangles without any fault lines.
Specifically, he proved that:

\begin{theorem}
Let $a$ and $b$ be relatively prime positive integers. A faultfree tiling of an $m\times n$ rectangle 
with $a\times b$ rectangles exists if and only if:

1. Either $m$ or $n$ is divisible by $a$, and either $m$ or $n$ is divisible by $b$. 

2. Each of $m$ and $n$ can be expressed as $xa+yb$ in atleast two ways, where $x$ and $y$ are positive 
integers; and 

3. For the case where the tiles are dominoes, the recatngle is not $6\times 6$. 
\end{theorem} 

In this paper, we extend the concept of faultfree tilings of rectangles to L-shaped trominoes and characterize 
rectangles that admit a faultfree tromino tiling. We propose a method for generating tilings of larger 
rectangles, using rules for extending tilings of smaller rectangles inductively, by adding only a constant 
number of trominoes, columns and rows. The scheme uses rules akin to (regular) grammars. 
In section 2, we show that $3\times n$ rectangles, where $n\geq 2$, do not admit any faultfree
tilings (Lemma 1). We also introduce the notion of a {\it critical tromino}, in the 
leftmost three columns of faultfree tilings, to characterize patterns 
in faultfree tilings for $m\times n$ rectangles (Lemma 2). Based on this 
characterization, we establish the feasibility of faultfree tilings for all
$m\times n$ rectangles, where $m,n\geq 4$, and $3|mn$, in Theorem 1. 

The notion of {\it crossings} is inherent in the very definition of a faultfree tiling. In Section 3, we introduce the 
notion of {\it horizontal} and {\it vertical crossing numbers}, thereby developing a measure of faultfree-ness 
for a tiling. The {\it horizontal (vertical) crossing number} is defined as the minimum of the total number 
of crossings on any horizontal (vertical) grid line. We address the problem of finding an upper bound on the 
minimum crossing numbers across rows and columns. We derive an upper bound of 2 for the horizontal and vertical 
crossing numbers for $m\times n$ rectangles, where $m,n\geq 10$; such a tiling is achieved by our method of extending faultfree tilings of smaller rectangles inductively. 
We also determine the extremal values achievable for crossing numbers. We show that no rectangle admits a 
faultfree tiling with crossing number exceeding $min\{2t-1, \frac{2n}{3}\}$, for $3t\times n$ rectangles, where $t\geq 1$. 

We then move on to address the problem of enumeration of faultfree tromino tilings of rectangles. 
In Section  4, we develop an iterative scheme for generating faultfree tromino tilings of rectangles and derive 
the known closed form expression of $8.6^{t-3}$, where $t\geq 3$, for the number of faultfree tromino tilings of 
$4\times 3t$ rectangles. We also derive the known exact generating function for the number of faultfree tromino tilings of $5\times 3t$ rectangles, where $t\geq 1$. This generating function is:

\begin{eqnarray*}
G(z)& = & \frac{24z^{2}(3+10z+15z^{2})}{1-2z-31z^{2}-40z^{3}-20z^{4}}
\end{eqnarray*}   

These counting problems were addressed earlier by Moore in \cite{moore}; we derive results
identical to his results using our own iterative scheme for extending tilings of smaller
rectangles into those of larger rectangles. In addition, we further move on 
to generalize our iterative scheme to arbitrary rectangles which have not been considered 
so far in the literature. The complexity of the task of enumeration 
emerges naturally from our scheme. For $6\times 6t$ rectangles, where $t\geq2$,  
we prove that the number of faultfree tromino tilings is at least $128.(t+1).144^{t-2}$. For $7\times 6t$ rectangles we derive a generating function where the coefficient of $z^{t}$ gives 
a lower bound on the number of faultfree tromino tilings. 

A natural question to ask is whether the number of domino tilings of a given 
$m\times n$ rectangle is more or less than the number of tromino tilings. 
Intuitively,  one might say that we have more freedom in tiling a rectangle with 
dominoes than with trominoes, by this we mean that there are possibly more ways of 
retiling a given area with dominoes than with trominoes.  However, no such comparisons 
have been reported so far in the literature. In Section 6, we present such a comparison, 
thereby producing an upper bound of $2^{\frac{4mn}{3}}\{min[N_D(m,2n),N_D(2m,n)]\}$ on the 
number of tromino tilings of a given $m\times n$ rectangle, where $3|mn$, $m,n>0$, and $N_D(m,n)$ denotes 
the number of domino tilings of the $m\times n$ rectangle. 

\subsection{Definitions and notation}

\begin{figure}[htbp]
\centerline{\scalebox{.2}{\includegraphics{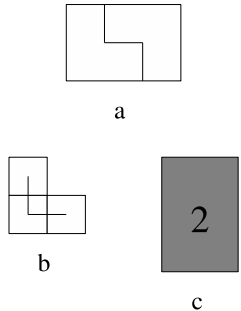}}}
\caption{Notations for a tromino tiling of $R(i,j)$.}
\end{figure}

Firstly, the reader should note that a $2\times3$ rectangle can be tiled with trominoes
(as shown in Figure 1(a)). We denote a rectangle with $i$ rows and $j$ columns by $R(i,j)$.
We will indicate decompositions into non-overlapping subrectangles by means of an additive 
notation. For example, a $3i\times 2j$ rectangle can be decomposed into $ij3\times 2$ subrectangles 
and we write this fact as $R(3i,2j) = \sum_{a=1}^{i}\sum_{b=1}^{j}R(3,2) = ijR(3,2)$. It follows 
from this and Figure 1(a) that any $2i\times 3j$ or $3i\times 2j$ rectangle can be tiled with trominoes. 
From now on, any rectangle decomposed into a combination of $3i\times 2j$ subrectangles, $2i\times 3j$ 
subrectangles and trominoes will be considered as successfully tiled by trominoes. Denote the $1\times 1$
square lying in row $i$ and column $j$ as $(i,j)$. To make the notation simple, trominoes are depicted 
in the rest of the paper as a composition of two lines forming an L-shape across an actual tromino 
(as shown in Figure 1(b)). We refer to a given arrangement of trominoes as a {\it pattern}. 
We define a {\it flip} of a pattern to be its reflection in the horizontal axis. 

Given a tiling, a line which cuts the rectangle into two pieces and yet does not pass through 
any of the tiles is called a {\it fault line}. A tiling of a rectangle which has no fault lines 
is called a {\it faultfree tiling}. We number each of the horizontal and 
vertical grid lines starting from $1$. So the left edge becomes the $1st$ 
vertical grid line and the upper edge becomes the $1st$ horizontal grid line.
If a grid line passes through a tromino, we say that the given tromino {\it crosses} over the 
grid line under consideration. We define the {\it critical horizontal (vertical) grid line} as 
the horizontal (vertical) grid line having the minimum number of crossings. 
We define the {\it horizontal (vertical) crossing number} as the total 
number of crossings on the critical horizontal (vertical) grid line. We call a tromino which crosses 
the $3rd$ horizontal (vertical) grid line a {\it critical tromino}. $R(2,3)$s are shown 
as gray-shaded rectangles labeled by the number $2$ since they can be tiled in two ways 
(as shown in Figure 1(c)). Parts having arbitrary tiling with tiles crossing over to the 
remaining portion of the rectangle are shown as dark shaded areas.  

\section{Incremental generation of faultfree tromino tilings}
\label{secffttr}

In this section, we show in Lemma 1 that the rectangle $R(3,n)$, where $n>2$, does not 
admit any faultfree tilings. We also introduce the notion of a {\it critical tromino} 
to characterize patterns of faultfree tilings in the leftmost three columns of a rectangle 
$R(m,n)$ in Lemma 2. Based on this notion, we establish the feasibility of faultfree 
tilings for all rectangles $R(m,n)$, where $m,n\geq 4$ and $3|mn$, in Theorem 2.  

\begin{lemma}
No rectangle $R(3,n)$ where $n>2$ has a faultfree tiling. 
\end{lemma}
\begin{proof}
\label{lemma3comman}
The only permissible orientations of a tromino 
covering $(1,1)$ and allowing a tiling are shown in Figure 2(ii)(a) and (b). 
In either case, the $3rd$ vertical grid line
becomes a fault line since there is only one way 
to cover $(3,1)$ in the first case and $(2,1)$ in the second case. \hfill \qed
\end{proof}

\begin{figure}[htbp]
\label{patterns}
\centerline{
\scalebox{.25}{\includegraphics{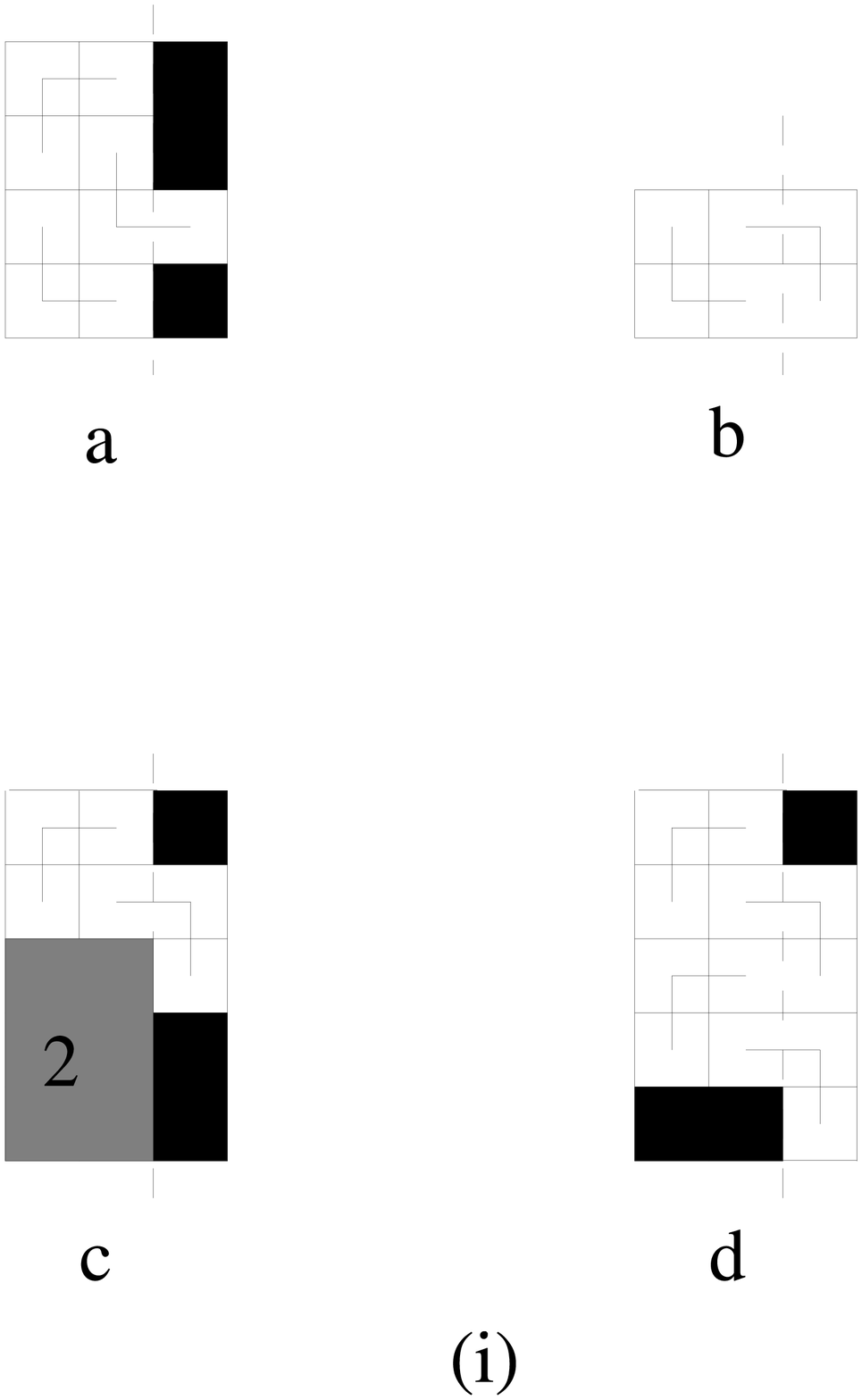}}
\scalebox{.25}{\includegraphics{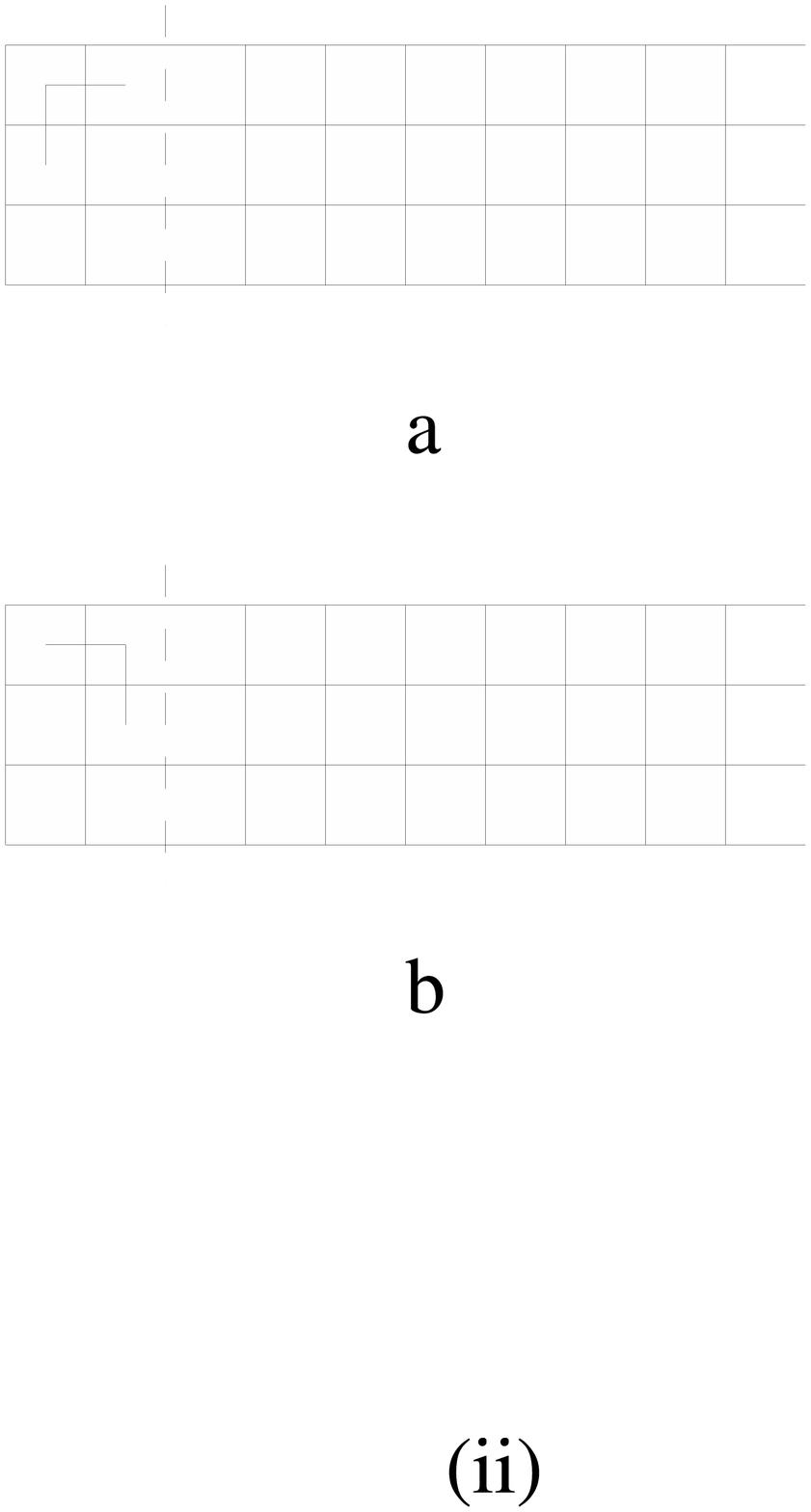}}
}
\caption{(i) Possible arrangements of trominoes around a {\it critical} tromino. (ii) Possible 
orientations of a tromino covering $(1,1)$ in $R(3,n)$.}
\end{figure}

Now consider any arbitrary faultfree tiling $R(m,n)$ where $m,n\geq 4$.
Since no tromino will fit in one column, we focus on that portion of this 
rectangle where a tromino crosses the $3rd$ vertical grid line. We call such a 
tromino a {\it critical tromino}. The only possible patterns arising for 
arrangements of trominoes to the left of a critical tromino are those shown 
in Figure 2(i). The $3rd$ vertical grid line is shown as a fine-dashed line.  

\begin{lemma}
\label{lemmagenerators}
The only possible patterns arising for arrangements of trominoes
around a critical tromino are those shown in Figure 2(i)
and their flipped counterparts.
\end{lemma}
\begin{proof}
Depending upon the number of squares (1 or 2) of the critical tromino 
to the left of the $3rd$ vertical grid line, two cases arise. 
In one case two of 
its squares lie to the left of this line. We set our coordinate system 
such that these two squares become $(2,2)$ and $(3,2)$. The only way to cover 
the squares $(2,1)$ and $(3,1)$ is shown in 
Figure 2(i)(a). 

In the other case only one square lies 
to the left of the $3rd$ vertical grid line. Without loss of generality, we assume 
that the vertical edge of the tromino faces downwards. Setting our coordinate system
such that this square is $(2,2)$, different cases 
arise depending on which 
tromino covers the square $(2,1)$.   
If the tromino covering $(2,1)$ also covers $(3,1)$, then we get the case shown in 
Figure 2(i)(b). If the tromino covering $(2,1)$ also covers 
$(1,1)$, the following three cases 
arise depending on how $(3,1)$ 
is covered. The case when the tromino covering $(3,1)$ also 
covers $(3,2)$ and $(4,2)$ is shown in Figure 2(i)(c). If the tromino 
covering $(3,1)$ covers $(3,2)$ and $(4,1)$, then 
two cases arise depending on how $(4,2)$ is covered. One of these cases is illustrated 
in Figure 2(i)(d) in which the tromino covering $(4,2)$ also covers $(4,3)$ and $(5,3)$. 
If the tromino covering $(4,2)$ also covers $(5,2)$ and $(5,1)$ we get back Figure 2(i)(c), 
if it also covers $(5,2)$ and $(5,3)$ or $(4,2)$, $(4,3)$ and $(5,2)$ then from row $3$ to row $6$ we 
get the pattern as in Figure 2(i)(a). Since only four 
orientations of this tromino are possible, we are done. \hfill \qed
\end{proof}

\begin{figure}[htbp]
\label{generators}
\centerline{
{\scalebox{.25}{\includegraphics{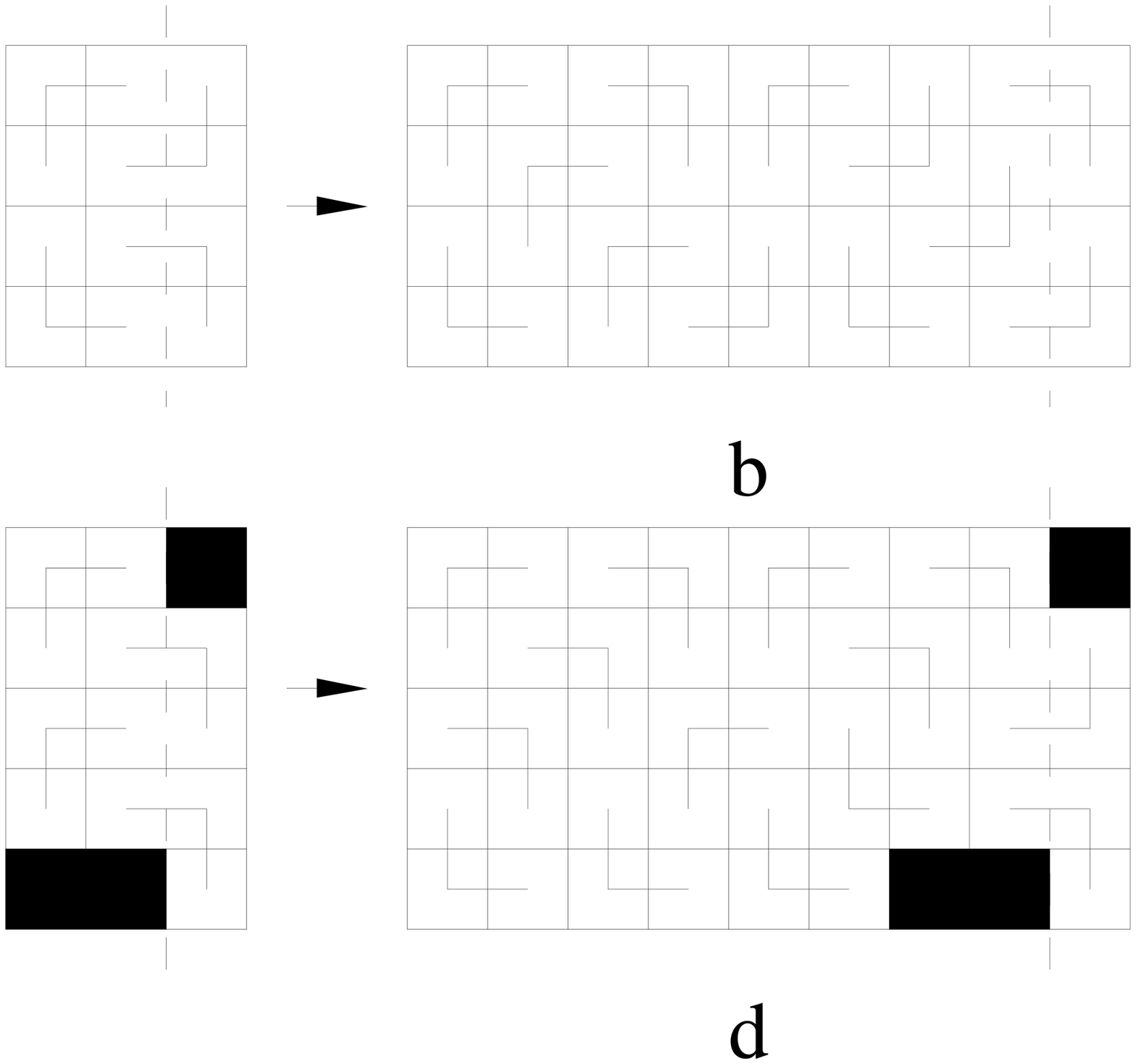}}}
{\scalebox{.25}{\includegraphics{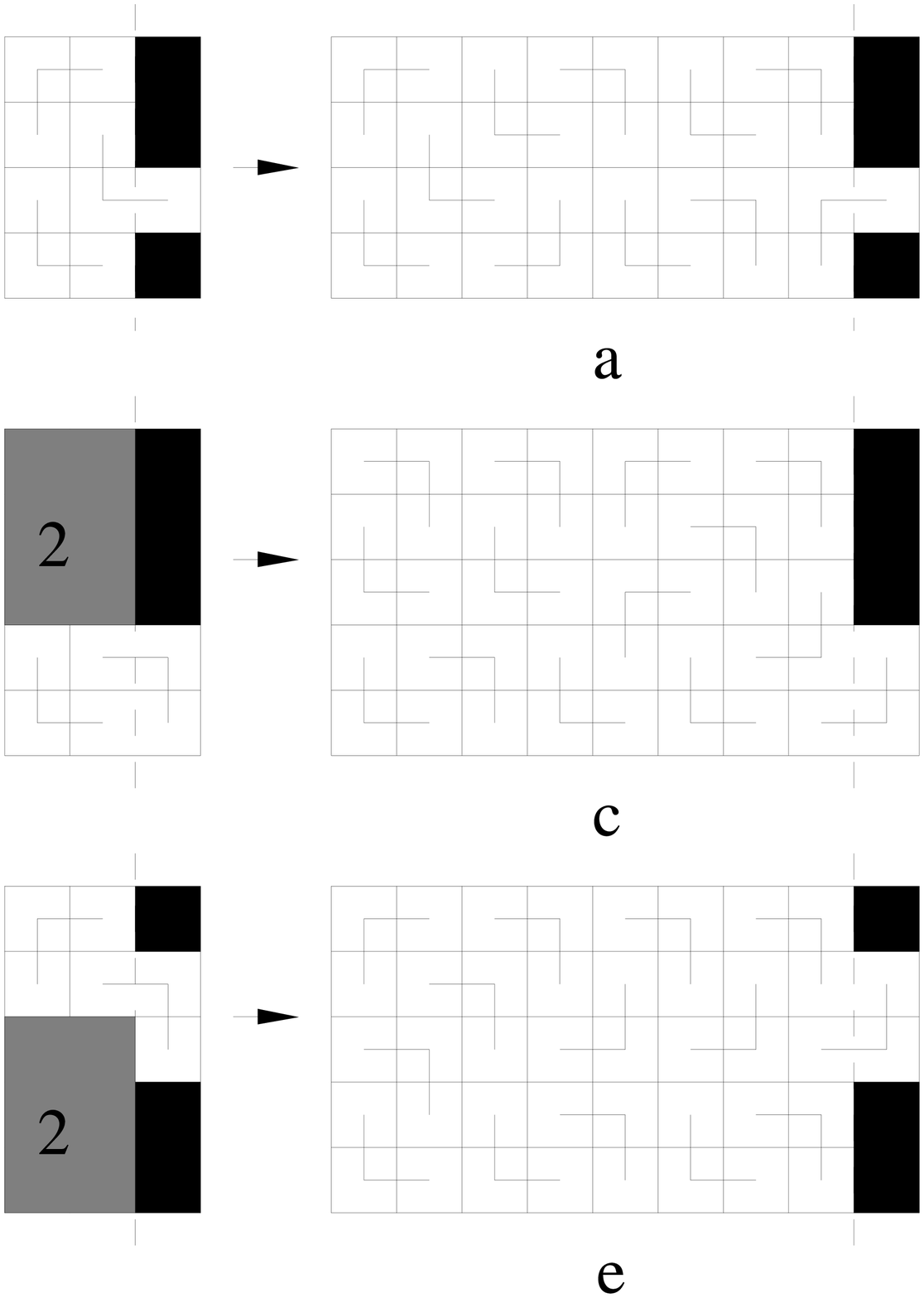}}}
}
\caption{Extension schemes for the patterns of Figure 2(i)}
\end{figure}



We follow an incremental approach where we produce a bigger faultfree 
tiling from a smaller one. 
We focus our attention 
on critical trominoes and consider the four patterns characterized in 
Lemma \ref{lemmagenerators}. Corresponding to a critical tromino, we identify 
the arrangement of trominoes around it; this arrangement must be one of 
the four patterns characterized in Lemma \ref{lemmagenerators} and depicted 
in Figure 2(i). Consider the case in Figure 2(i)(a). 
We extend the tiling of the four rows by six columns as shown in Figure 3(a). 
Note that the three dark squares are left untouched and the rest of the 
squares in the three columns as well as the six new columns are tiled in 
such a way that all the vertical grid lines have some crossing(s). 
We call such an extention pattern, a {\it generator}.
All the remaining rows of 
the smaller faultfree tiled rectangle 
can now be extended by tiling with $R(3,2)$ rectangles or any 
other suitable tiling scheme. Since the newly introduced columns already have 
some crossings, any tiling scheme for extending remaining rows will produce a 
faultfree tiling. Thus we get a faultfree tiling of a 
rectangle with six more columns. Corresponding to the 
case in Figure 2(i)(b), we may identify 
the exhaustive cases (to the left of the arrow) as in Figure 3(b) and (c). If the arrangement 
of trominoes above the pattern shown in Figure 2(i)(b) is not 
one amongst the two cases shown in Figure 3(b,c) then it must be one 
amongst the other three cases shown in Figure 2(i)(a,c,d) 
by Lemma \ref{lemmagenerators}. All the five cases 
enumerated in Figure 3 are called {\it generators}. In each 
of these generators the three column tiling on the left of the arrow is 
extended to the tiling shown to the right of the arrow. We  
extend each of the patterns by six columns as shown in 
Figure 3, ensuring that the
larger tiling produced is also faultfree. 
We will refer to this scheme of extension by six columns as {\it incremental 
generative scheme}. 
Although we may extend a tromino tilable rectangle 
by a minimum of three columns, we consider extensions only 
by six columns due to the following reason. 
Extending the entire rectangle by three columns  
may require us to extend a block of 
an odd number of consecutive rows $i$ through $j\geq i$ (separately), by three 
columns, where the rows numbered below $i$ and above $j$ (if any), are 
extended by one of the other cases in Figure 3.
So we may be required to tile $R(3,j-i+1)$ where $j-i+1$ is odd. Since 
such a rectangle cannot be tiled (see Chu and Johnsonbaugh \cite{chu}), 
we increase six columns at a time as in Figure 3.
Thus we can extend any rectangle $R(i,j)$ using our incremental generative
scheme to produce any one 
of $R(i,j+6)$, $R(i+6,j)$ or $R(i+6,j+6)$. 
Now we are ready to argue and show that all rectangles $R(m,n)$, where $m,n\geq4$ 
and $3|mn$, can have a faultfree tiling. 
All multiples of $3$ greater 
than or equal to $6$ can be written in the form $6t+6$ or $6t+9$, 
where $t\geq0$. So, all multiples of $3$ can 
be obtained using our incremental generative approach if we start 
with a rectangle having $6$ or $9$ as one of its sides. 
From the Chu-Johnsonnbaugh Theorem \cite{chu}, we know that any tileable
rectangle must have area divisible by $3$, and so at least one side 
must be a multiple of $3$. All integers 
greater than or equal to $4$ can be written in the form $6t+i$, 
where $t\geq0$ and $4\leq i \leq9$, since all 
these $i's$ are mutually incongruent modulo $6$. 
We identify $11$ {\it basis} cases $R(i,6)$ and $R(i,9)$ 
where $4\leq i\leq9$ as shown in Figure 4. 
Starting from these rectangles and using our 
incremental generative approach we can obtain 
any rectangle $R(6t_1+i,3t)$, where $6t_1+i\geq4$ and 
$t\geq2$. We summarize our result in the following theorem.

%

\begin{theorem}{\bf [Faultfree tromino tiling theorem]}
\label{theoremfaultfreetrominotiling}
All rectangles $R(m,n)$, where $m,n \geq 4$ and $3|mn$, admit a faultfree tiling. 
\end{theorem}
 
\begin{figure}[htbp]
\label{basis}
\centerline{
{\scalebox{.3}{\includegraphics{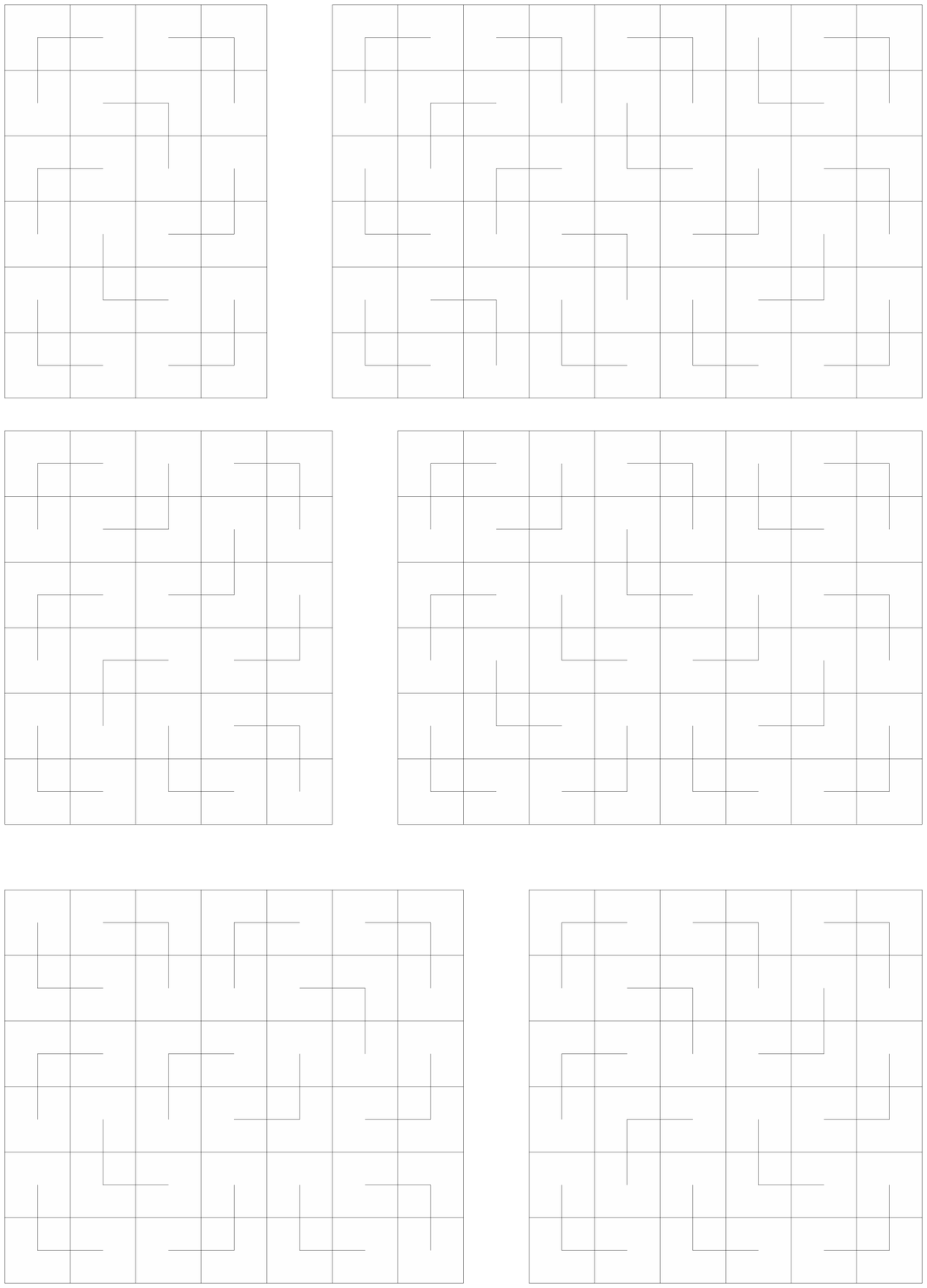}}}
{\scalebox{.3}{\includegraphics{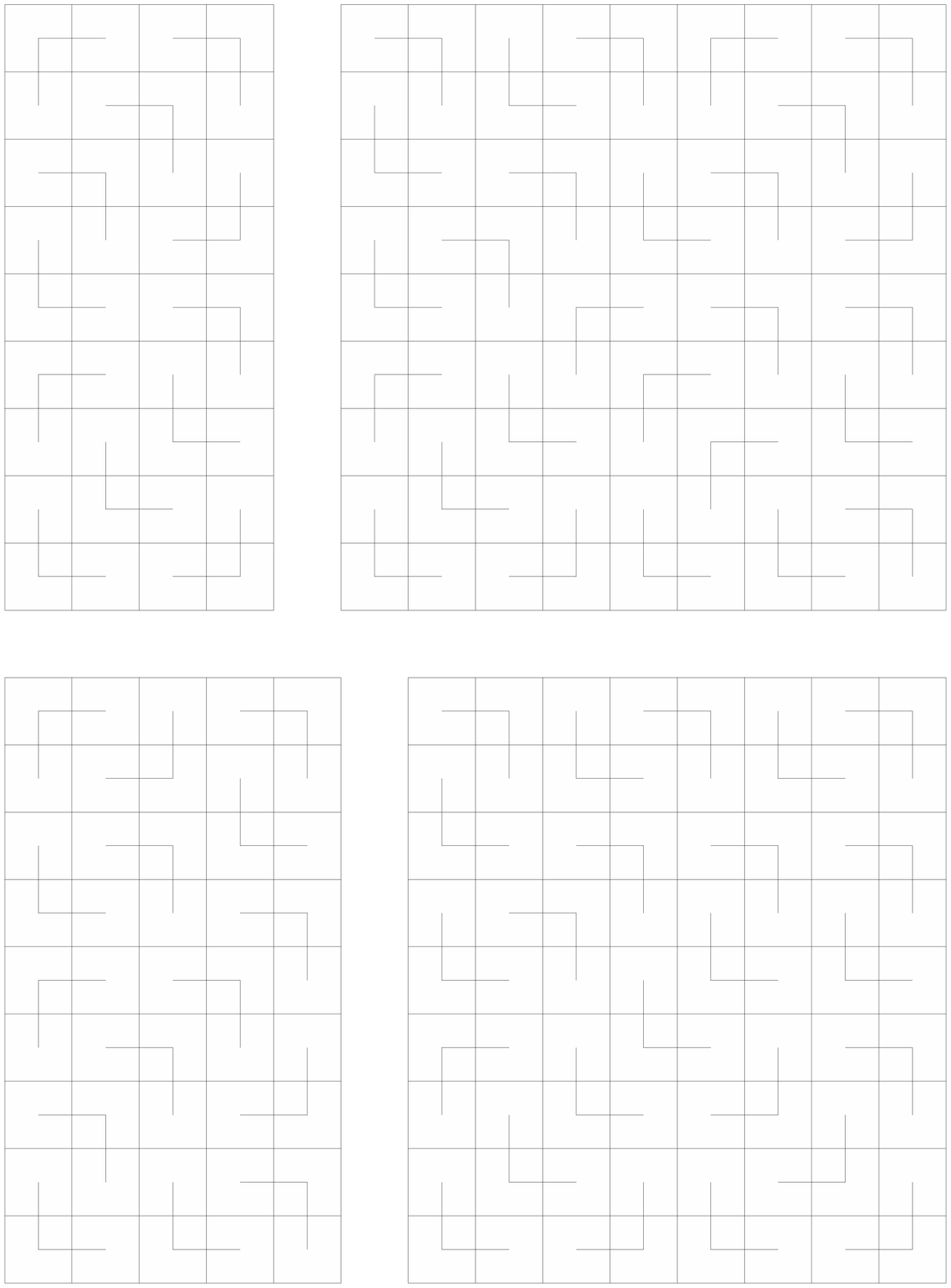}}}
{\scalebox{.15}{\includegraphics{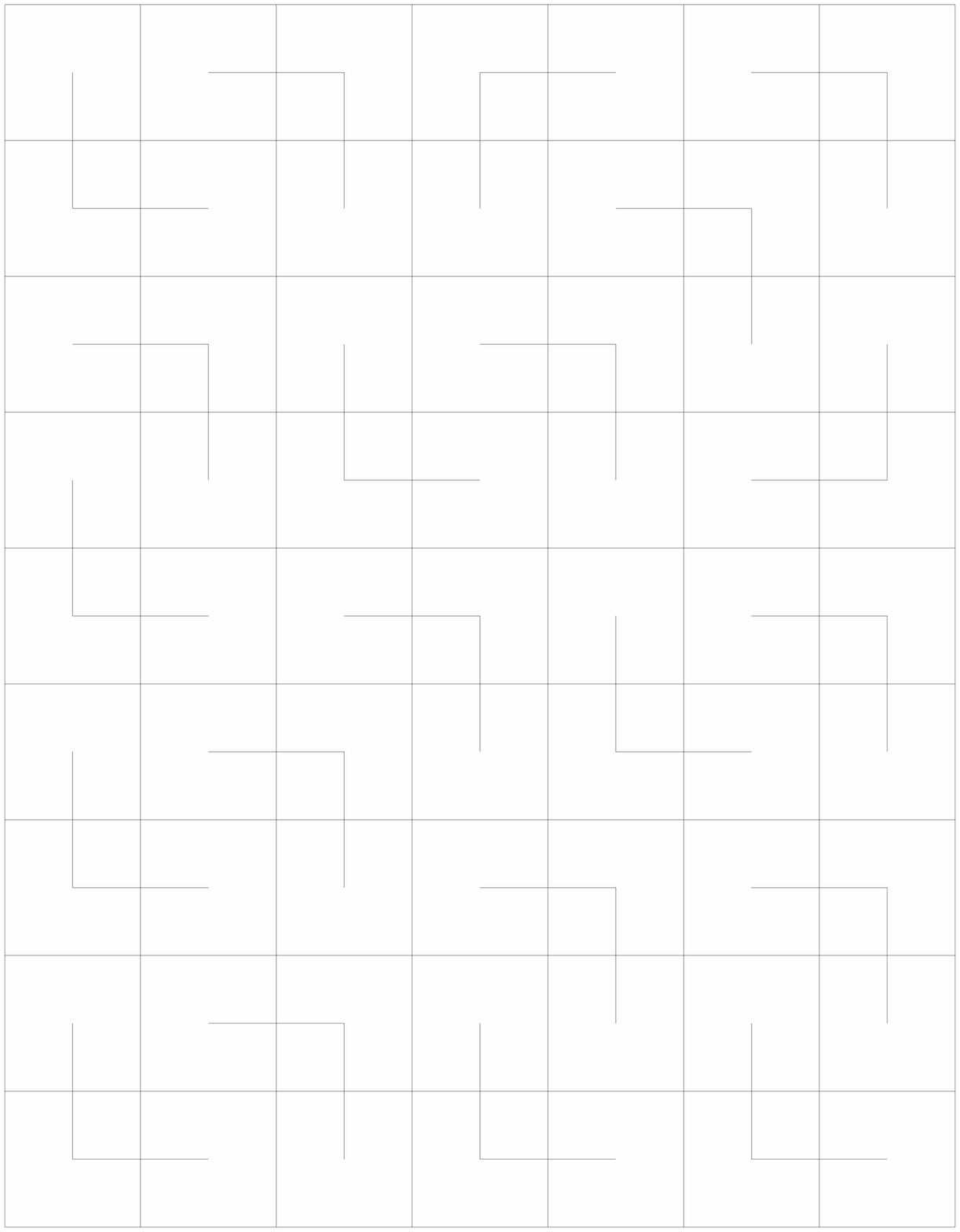}}}
}
\caption{Faultfree tromino tilings for the basis cases}
\end{figure}

%

\section{Crossing numbers of faultfree tilings}
\label{seccrossing}

The notion of {\it crossings} is inherent in the very definition of a faultfree tiling. 
In this section, we introduce the notion of {\it horizontal} and {\it vertical crossing 
numbers}, thereby developing a measure of faultfree-ness for a tiling. 
We define the {\it horizontal (vertical) crossing number} as the total 
number of crossings on the critical horizontal (vertical) grid line. 
We address the problem of finding an upper bound on 
the minimum crossing numbers across rows and columns. 
We show that both the horizontal and vertical 
crossing numbers can be as low as 2 for all rectangles $R(m,n)$, where $m,n\geq 10$ and $3|mn$, 
admitting faultfree tilings; such a tiling is achieved by our
method of extending faultfree tilings by six columns (or rows) as
explained in Section \ref{}. We have the following theorem.

\begin{theorem}{\bf [Minimum crossing number theorem]}
\label{}
All rectangles $R(m,n)$, where $m,n\geq 10$ and $3|mn$, admit a faultfree tromino tiling with horizonal and vertical crossing numbers
less than or equal to 2.
\end{theorem}
\begin{proof}
For constructing a faultfree tiling of $R(m,n)$ where $m,n\geq10$ and $3|mn$, 
we first construct a faultfree tromino tiling of $R(m-6,n-6)$ and then expand either 
side by six columns using the procedure discussed above. We note that in each of 
the generators, the previous edge now has two crossings. Apart from identifying one 
area on both sides of $R(m-6,n-6)$, we do not expand any other area by the generator 
approach, but tile the remaining area using $R(3,2)$s, or any other suitable approach.
Thus both the edges of $R(m-6,n-6)$, now have only two crossings. We do not consider 
the critical horizontal (vertical) grid line of $R(m-6,n-6)$ because that line may now 
have some more crossings due to expansion of both the edges. \hfill \qed
\end{proof}

An interesting question is whether faultfree tilings can have crossing number 1.
There exist tilings of rectangles with horizontal crossing number 1 (for example, 
consider $R(6,4)$, $R(6,5)$, $R(6,7)$, $R(9,4)$, $R(6,8)$ in Figure 4), although we do 
not know whether the rectangle $R(m,n)$, where $m,n>0$ and $3|mn$, can always admit tromino 
tilings with horizontal and/or vertical crossing numbers equal to 1.
\begin{openquestion}
Characterize $m$ and $n$, such that a faultfree tromino tiling of the rectangle $R(m,n)$, 
where $m,n>0$ and $3|mn$, has horizontal and/or vertical crossing numbers equal to 1.
\end{openquestion}

\subsection{Maximum Crossing Numbers}

Consider an arbitrary rectangle $R(m,n)$ and suppose that the horizontal and vertical crossing numbers are 
at least $k$. So the total number of crossings is at least $k(m+n-2)$ (counting the number 
of crossings on each grid line). Now, any tromino tiling of $R(m,n)$ has $\frac{mn}{3}$ trominoes and each tromino cuts 
a horizontal as well as a vertical grid line. So any tromino tiling of $R(m,n)$ has $2\times \frac{mn}{3}$ crossings. We have the following inequality:

\begin{eqnarray}
\frac{2mn}{3} \geq k(m+n-2)
\end{eqnarray}
 
This inequality has some interesting implications. 

For $m = 3$, $2n\geq k(n+1)$. Note that $k\leq 1$ since $\frac{2n}{n+1}<2$. 
Indeed we know from Lemma 1 that no rectangle $R(3,n)$ has a faultfree tiling.  

Further for $m = 6$, $4n\geq k(n+4)$. Here, each grid line can have 
at most 3 crossings, since $\frac{4n}{n+4}<4$. A natural question is whether there exists a 
tiling of $R(6,12)$ such that every grid line has {\it exactly} 3 crossings. 
We now present the following theorem:

\begin{theorem}{\bf [Maximum crossing number theorem]}
For any given rectangle $R(m,n)$, where $m = 3t$ and $t,n\geq1$, the horizontal (vertical) crossing number $k$ 
can be at most the 
minimum of $(2t-1)$ and $\frac{2n}{3}$,  i.e.,

\begin{eqnarray}
k \leq min\{2t-1, \frac{2n}{3}\} 
\end{eqnarray}

\end{theorem}
\begin{proof}
From the Chu-Johnbonsonbaugh Theorem \cite{chu}, we know that for a rectangle $R(m,n)$ to admit 
a tromino tiling, atleast one of its sides must be divisible by 3. Assuming $m = 3t$, where $t\geq1$, 
and using inequality (1) we get, 
\begin{eqnarray*}
   & => & 2nt \geq k(3t+n-2) \\
   & => & 2nt \geq kn + k(3t-2)
\end{eqnarray*}

Since $t\geq1$, the second term in the above inequality, $k(3t-2)$ is positive. If $k=2t$, then 
the right side becomes greater than the left side in the above inequality. So we conclude that $k \leq 2t-1$. 
Also, note that no tromino can fit in one column, so $n\geq2$. Using this fact and inequality (1), we get,
\begin{eqnarray*} 
   & => & 2nt \geq k(3t+n-2) \geq 3kt \\
   & => & k \leq \frac{2n}{3} 
\end{eqnarray*}
   
Surprisingly, the first (second) upper bound is not dependant on $n$ ($m$). Both $2t-1$ and $\frac{2n}{3}$ are upper 
bounds, whence the result. \hfill \qed
\end{proof}

\section{Counting the number of faultfree tilings of $R(4,3t)$ and $R(5,3t)$}
\label{r43t}

We now demonstrate how our incremental generative scheme 
of Section \ref{secffttr} can be used to count the total number of 
faultfree tilings for $R(4,3t)$, $t\geq 2$. 
For such rectangles, the leftmost three columns cannot 
match the pattern in Figure 3(b) since the $4th$ vertical grid
line would then be a faultline. So the leftmost three columns 
must match the pattern in Figure 3(a).
We count the total number of ways of increasing $t$ by $1$, and then 
multiply with the total number of faultfree tilings 
of $R(4,3t)$ to get the total number of 
faultfree tilings of $R(4,3t+3)$. 

Consider a faultfree tiling of $R(4,6)$. Since the 
tiling of the leftmost three columns in such a tiling must 
match the pattern in Figure 2(i)(a) (see Lemma 2), the
only possible way of tiling $R(4,6)$ faultfreely is
shown in Figure 4. Using Lemma 2 symmetrically
to the rightmost three columns, we note that 
their tiling should also match the pattern in Figure 2(i)(a), whence
the tiling in Figure 4. 

\begin{figure}[htbp]
\centerline{\scalebox{.3}{\includegraphics{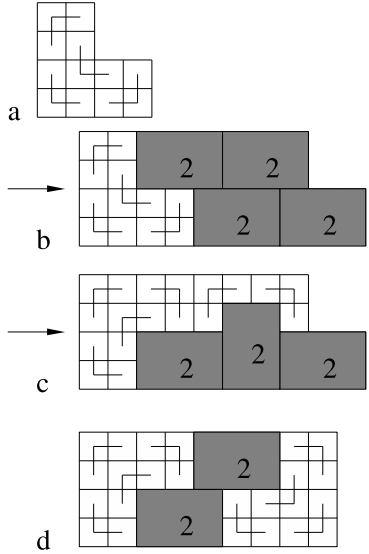}}}
\caption{Number of ways of extending $R(4,6)$ to make $R(4,12)$.}
\end{figure}

We now present the following lemma.
\begin{lemma}
In any faultfree tromino tiling of $R(4,3t)$, where $t>2$, the tiling of the leftmost four columns 
must match the pattern in Figure 6(a) or its flipped counterpart.
\end{lemma}
\begin{proof}
For $t>2$, we show that the pattern shown in Figure 5(a) will 
have $R(2,3)$ by its side as shown in 
Figure 6(a). 
We set our coordinate system such that 
the square to the top left becomes $(1,1)$. Suppose that the 
tiling of the leftmost four columns of $R(4,3t)$ did not match the pattern 
shown in Figure 6(a). Then the only other possible way of faultfreely tiling $R(4,3t)$ would be 
that a tromino covers the squares $(1,3)$, $(1,4)$, $(2,3)$ 
and another tromino covers $(2,4)$, $(2,5)$ and $(3,5)$. The reader can now easily see that the 
two trominoes covering $(1,5)$ and $(4,5)$ will make the $7th$ 
vertical grid line a fault line. So, we conclude that the 
tiling of the leftmost four columns in any faultfree tromino tiling of $R(4,3t)$, 
where $t>2$, must match the pattern shown in Figure 6(a) or its symmetric 
counterpart. \hfill \qed
\end{proof}

Now consider extending a tiling of $R(4,6)$ to a tiling of 
$R(4,12)$. 
We generate $R(4,12)$ from $R(4,6)$ using the pattern shown in 
Figure 5(a) as explained in the following lemma.

\begin{lemma}
The only possible ways of extending a faultfree tiling of $R(4,6)$
to that of $R(4,12)$ using the pattern in Figure 5(a) are by the generators shown in Figure 5.
\end{lemma}
\begin{proof}
From Lemma 3, we know that any tiling of the leftmost 
four columns of $R(4,12)$ must match the pattern shown in Figure 6(a) or 
its flipped counterpart. 
Due to symmetry, the rightmost four columns of $R(4,12)$ 
must also match the pattern
in Figure 6(a) or its flipped counterpart. So, we must have a $R(2,3)$ 
at the bottom right of the generators shown in Figure 5.
Consider the case when a tiling of the leftmost four columns of $R(4,12)$ 
matches the pattern shown in Figure 6(a). The only possible ways of tiling 
the remaining area are shown in Figure 5(b). This is because there are only 
two orientations for a tromino covering the squares $(3,5)$ and $(4,5)$. If it 
covers $(3,6)$ then the bottom $R(2,3)$ is completed by the tromino covering 
$(4,6)$ and if it covers $(4,6)$, then the bottom $R(2,3)$ is completed by the 
tromino covering $(4,7)$; whence the tiling shown in Figure 5(b). 
Now consider the case when a tiling of the leftmost 
four columns of $R(4,12)$ matches the flipped counterpart of the pattern shown 
in Figure 6(a). In this case the only possible ways of tiling the remaining area 
are shown in Figure 5(c). There are two orientations for a tromino covering the 
squares $(4,6)$ and $(4,7)$. If it covers $(3,6)$ then the tromino covering $(3,7)$ 
must cover $(2,6)$ and $(2,7)$, otherwise, in the case it covers $(2,7)$ and $(2,8)$, 
the reader can easily see that the square $(1,8)$ becomes inaccessible.  
If the tromino covering $(4,6)$ and $(4,7)$ covers $(3,7)$ then, following a similar 
reasoning as above, the reader can see that the tromino covering $(3,6)$ must complete a $R(3,2)$ with it. 
Since the above two cases exhaust all possible 
arrangements of trominoes, we conclude that the generators shown in Figure 5 
are the only possible ways of extending a faultfree tromino tiling of $R(4,6)$ to that 
of $R(4,12)$, using the pattern shown in Figure 5(a). \hfill \qed
\end{proof}

\begin{figure}[htbp]
\centerline{
{\scalebox{.3}{\includegraphics{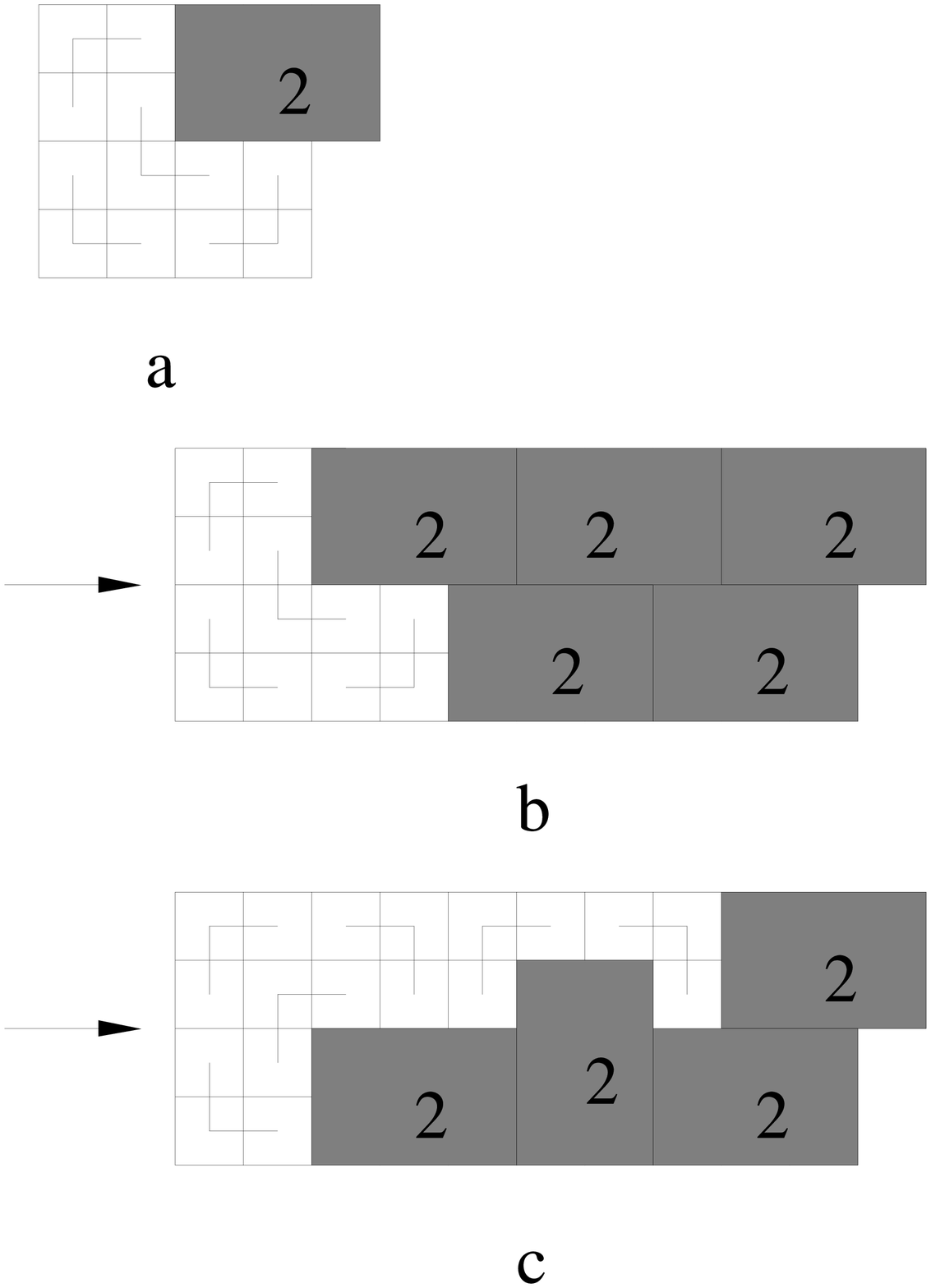}}}
{\scalebox{.3}{\includegraphics{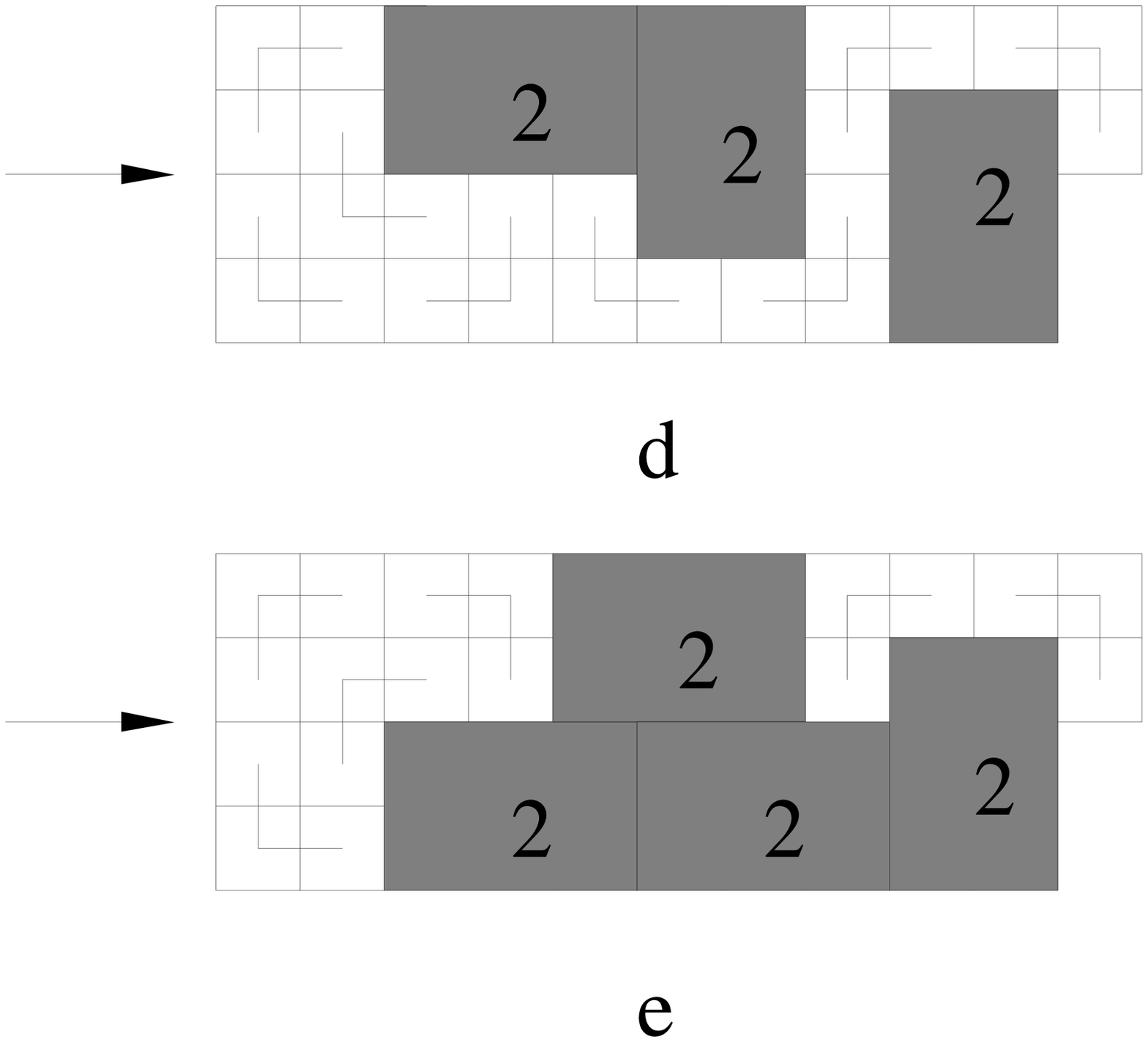}}}
}
\caption{Number of ways of extending $R(4,3t)$ when $t>2$.}
\end{figure}  

Following Lemma 4, we observe that 
the total number of ways of tiling $R(4,12)$ = $2\times(16+8)$ = $48$ (we multiply 
our result by $2$ since there are two ways of 
tiling $R(4,6)$ faultfreely). The total number of ways of 
tiling $R(4,9)$ faultfreely = $2\times4$ = 8 
(as shown in Figure 5, we multiply by $2$ to take into account the symmetric 
counterpart of Figure 5(d)). Since $t$ is 
greater than 2, in any faultfree tromino tiling of $R(4,3t)$ the tiling of the 
leftmost four columns must match the pattern shown in Figure 6(a). 
As in the case of $R(4,6)$ (see Figure 4), the tiling of the
rightmost four columns of $R(4,9)$ must also match Figure 6(a), 
yielding the only possible faultfree tiling of $R(4,9)$ (see Figure 5(d)).

Now consider extending a faultfree tiling of $R(4,3t)$ to that of
$R(4,3t+6)$, where $t>2$. 
We know from Lemma 3 that the tiling of the leftmost four columns of $R(4,3t)$ 
must match the pattern in Figure 6(a) or its flipped counterpart.
So, we use the four generators in Figure 6 for extending a faultfree
tiling of $R(4,3t)$ to a faultfree tiling of $R(4,3t+3)$, 
where $t>2$, using the pattern shown in Figure 6(a). We state the following lemma.
\begin{lemma}
The only possible ways of extending a faultfree tiling of $R(4,3t)$
to that of $R(4,3t+6)$ using the pattern in Figure 6(a) are by the generators shown in Figure 6. 
\end{lemma}
\begin{proof}
From Lemma 3, we know that in any faultfree tromino tiling of $R(4,3t+6)$ the tiling of 
the leftmost four columns must match the pattern shown in Figure 6(a) (say Case 1) or its 
flipped counterpart (say Case 2). Consider 
Case 1 when this pattern is the same as that shown in Figure 6(a). 
We claim that Figure 6(b) and (d) show all possible ways of tiling the remaining area. 
Two orientations are possible for the tromino covering $(3,5)$ and $(4,5)$, 
it may cover $(3,6)$ or $(4,6)$. If it covers $(3,6)$ then a $R(2,3)$ is 
completed by the tromino covering $(4,6)$. If it covers $(4,6)$ then two more 
cases arise (say Cases 1(a) and 1(b)),
depending on how $(4,7)$ is covered. In Case 1(a), the tromino covering 
$(4,7)$ also covers $(3,6)$ and $(3,7)$; so we again get a $R(2,3)$ at the 
bottom. Now the remaining area can only be tiled as shown in Figure 6(b). 
This is because the tromino covering $(4,10)$ must also cover $(3,10)$. If 
this is not the case, then the reader can easily see that it has to cover the 
squares $(4,9)$ and $(3,9)$. But in this case $(4,8)$ becomes inaccessible. 
In case of the other two orientations (when the tromino covers $(3,10)$, $(4,10)$, 
$(3,9)$ and $(3,10)$, $(4,9)$, $(4,10)$) the reader can verify 
that we will get a $R(2,3)$ at the bottom right. It is easy to see that the remaining area can only be 
tiled by two $R(2,3)$ rectangles. In Case 1(b), the tromino covering $(4,7)$ 
also covers $(3,8)$ 
and $(4,8)$; then we claim that the only ways of tiling the remaining area 
are shown in Figure 6(d). As can be easily seen, the tromino 
covering $(3,6)$ must also cover $(3,7)$. If this tromino 
covers $(2,6)$ then a $R(3,2)$ is completed by the tromino covering $(1,6)$ and 
if it covers $(2,7)$ then $R(3,2)$ is completed by the tromino covering $(2,6)$. 
As proved in Lemma 3, the remaining area can only be tiled as shown in Figure 
6(d). The reader can verify that all other orientations of the tromino covering 
$(4,7)$ do not permit a faultfree tromino tiling. 

Now consider Case 2 where a tiling of the leftmost four columns of $R(4,3t+6)$ 
matches the flipped counterpart of the pattern shown in Figure 6(a). In this case 
we show that Figures 6(c) and (e) depict all possible ways of 
tiling the remaining area as follows. Two 
orientations are possible for the tromino covering the squares $(1,5)$ and $(2,5)$, 
it may cover $(2,6)$ or $(1,6)$. If it covers $(2,6)$ then a $R(2,3)$ is completed 
by the tromino covering $(1,6)$. If it covers $(1,6)$ then two cases (say Cases 2(a) and
2(b)) arise depending on how $(1,7)$ is covered. In Case 2(a), the tromino covering $(1,7)$ 
also covers $(2,6)$ and $(2,7)$; so we again get a $R(2,3)$. 
In this case the only possible ways of tiling the remaining area are shown in Figure 
6(e). Two orientations are possible for the tromino covering $(3,6)$ and $(4,6)$. 
If it covers $(3,7)$ then a $R(2,3)$ is completed by the tromino covering $(4,7)$. If 
it covers $(4,7)$ then the tromino covering the squares $(3,7)$ and $(3,8)$ may 
cover $(2,8)$ or $(4,8)$. If it covers $(4,8)$ then we get a $R(2,3)$ again. If it 
covers $(2,8)$ then the trominoes covering $(1,8)$ and $(4,8)$ make the $10th$ 
vertical grid line a fault line. So this orientation is not permissible. As has 
been proved in Lemma 3 the remaining area can only be tiled as shown in Figure 6(e).
In Case 2(b), the tromino covering $(1,7)$ also covers $(1,8)$ and $(2,8)$; 
then all possible 
ways of tiling the remaining area are shown in Figure 6(c). If the tromino covering 
the squares $(2,6)$ and $(2,7)$ covers $(3,6)$ then a $R(3,2)$ is completed by the 
tromino covering $(4,6)$ and if it covers $(3,7)$ then a $R(3,2)$ is completed by the 
tromino covering $(3,6)$. It has been proved in Lemma 3 that the remaining area can 
only be tiled as shown in Figure 6(c). The reader can verify that all other orientations 
of the tromino covering $(1,7)$ do not permit a faultfree tromino tiling. Thus, Figure 
6 shows all possible ways of extending by six columns a faultfree tromino tiling of $R(4,3t)$
to that of $R(4,3t+6)$ using the pattern shown in Figure 6(a).       \hfill \qed
\end{proof} 

Now we show how to generate and count all faultfree tilings of $R(4,3t)$, where $t\geq3$, 
by our incremental generative scheme. 

\begin{theorem}
The number of distinct faultfree tilings of $R(4,3t)$, where $t\geq 3$, is 
$8.6^{t-3}$.
Starting with $R(4,9)$ and $R(4,12)$, and using the generators in Figure 6,
it is possible to generate all the distinct faultfree 
tilings of $R(4,3t)$, where $t\geq3$, exactly once.
\end{theorem}
\begin{proof}
We perform induction on t. From Lemmas 3 and 4 
we know how to generate all the 8 and 48 
possible (distinct) faultfree tilings of $R(4,9)$ and $R(4,12)$, respectively. 
Note that $8.6^{3-3}=8$ and $8.6^{4-3}=48$, satisfying the count
$8.6^{t-3}$ for $t=3$ and $t=4$, respectively.
We consider the tilings of $R(4,9)$ and $R(4,12)$ as
the basis cases. For the 
inductive part, suppose we have generated all the 
$8.6^{t-3}$ 
possible distinct faultfree tilings of 
$R(4,3t)$ and we have to 
generate all possible faultfree tilings of $R(4,3t+6)$. 
From Lemma 3, we know that in any faultfree tromino tiling of $R(4,3t)$ the tiling of the 
leftmost four columns must match the pattern shown in Figure 6(a). 
Suppose we remove this pattern from $R(4,3t)$. Then, by our inductive hypothesis, 
we have all the $4.6^{t-3}$ distinct possible ways of 
tiling the remaining portion of $R(4,3t)$ because the 
pattern removed (as in Figure 6(a)) has a $R(2,3)$, which has two 
ways of tiling. 
The generators shown in Figure 6(d) and (e) reconfigure the 
trominoes in this $R(2,3)$ 
attached to the initial pattern, while those in Figure
6(b) and 6(c) do not. So, the number of ways of tiling 
$R(4,3t+6)$ distinctly, is $4.6^{t-3}\times(32+16+8+16)=8.6^{t-1}=8.6^{(t+2)-3}$,
establishing the induction. Observe that the generative scheme for 
counting the possible tilings is
also constructive, thereby providing a method to generate all the possible 
tilings.
\hfill \qed


\end{proof}

In \cite{moore}, Moore showed that there are $8.6^{t-3}$ ways of tiling 
$R(4,3t)$, $t\geq 3$ faultfreely. Theorem 5 summarises an 
alternate derivation of this formula using our incremental generative approach thereby
also providing a method for the enumeration of the tilings. 

\subsection{Faultfree tilings of $5\times 3t$ rectangles} 

We now count all faultfree tilings of $R(5,3t)$, where $t\geq2$. We start with the 
following lemma. 
\begin{lemma}
In any faultfree tiling of $R(5,3t)$, the tiling of the leftmost three 
columns must match the patterns in Figure 3(c) or 3(e), or their flipped counterparts.
\end{lemma}
\begin{proof}
Consider the arrangement of trominoes around the critical tromino crossing the $3rd$ 
vertical grid line. Since there are only five rows, by Lemma 2, the tiling of the leftmost 
three columns must match Figure 3(c), 3(d), or 3(e). The pattern 3(d) will never occur as it 
will make $(5,1)$ and $(5,2)$ inaccessible. So any tiling of the leftmost three columns must 
match the pattern shown in Figure 3(c) or 3(e), or their flipped counterparts.    \hfill \qed
\end{proof} 

\begin{figure}[htbp]
\centerline{
{\scalebox{.3}{\includegraphics{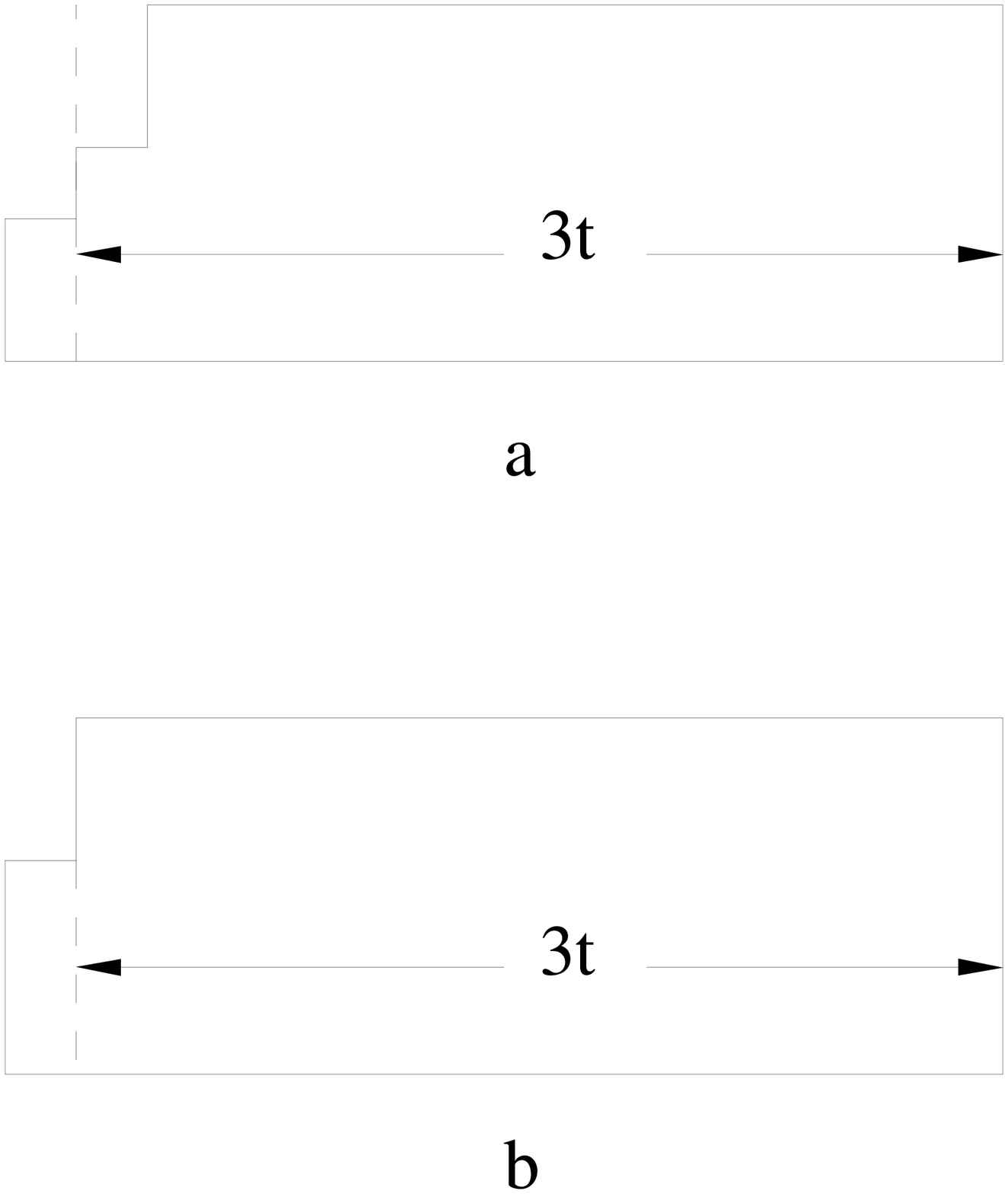}}}
{\scalebox{.3}{\includegraphics{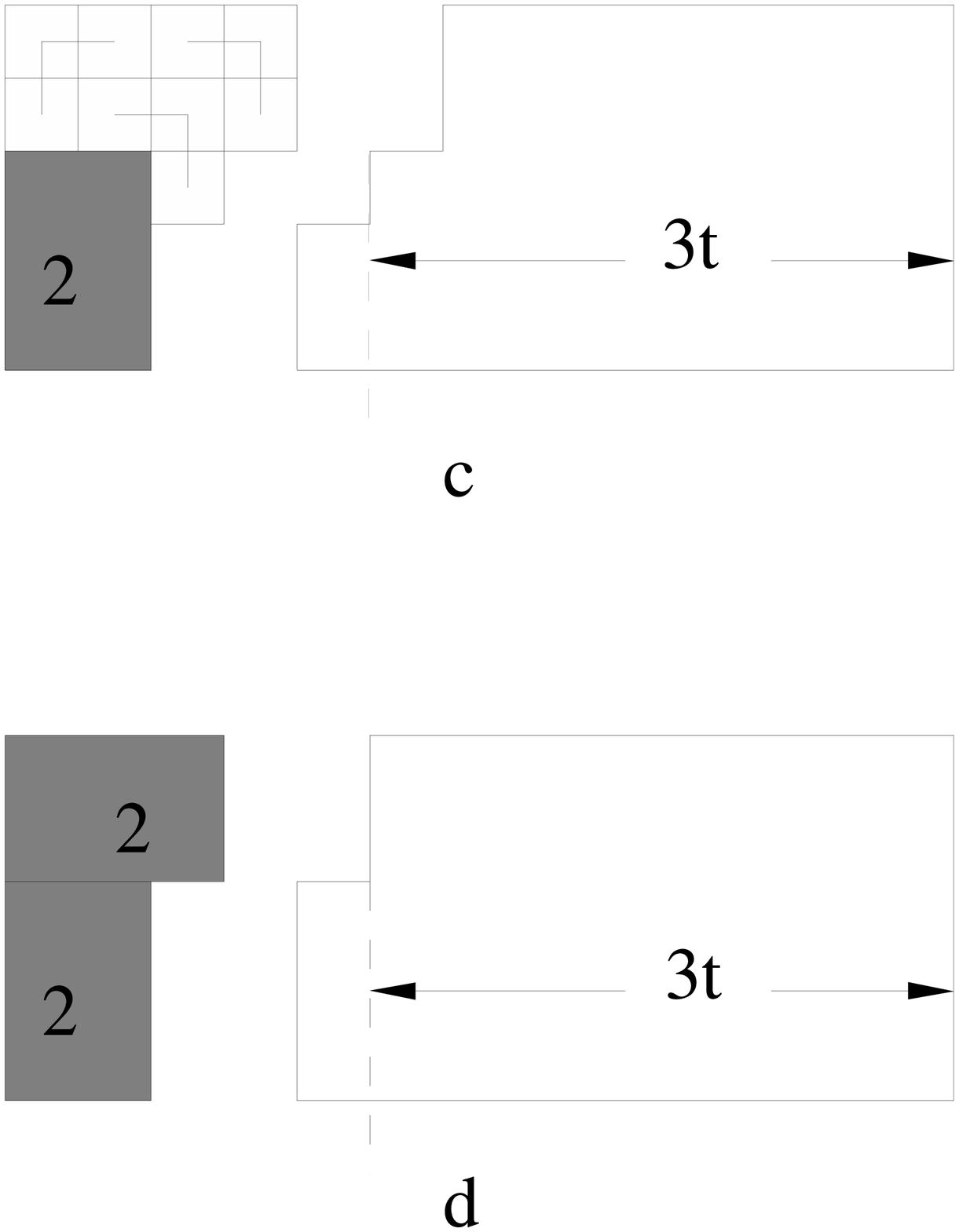}}}
}
\caption{Possible interfaces in $R(5,3t)$.}
\end{figure}  

Consider the case when the tiling of the leftmost three columns matches the pattern shown in 
Figure 3(e). The reader can easily see that the tromino covering $(1,3)$ will also cover $(1,4)$ 
and $(2,4)$. Now consider the interfaces shown in Figure 7(a) and (b). As shown in Figure 7(c) and 
(d), when the tiling of the leftmost three columns of $R(5,3t)$ matches Figure 3(e), the rest of 
the area matches the interface shown in Figure 7(a), and when the tiling of the leftmost three 
columns matches Figure 3(c), the rest of the area matches the interface shown in Figure 7(b). 
The interface shown in Figure 7(a) is called {\it slope} and that shown in Figure 7(b) is called 
{\it little jog}. Suppose $N(t)$, $N_1(t)$ and $N_2(t)$ denote the number of faultfree tilings of 
$R(5,3t)$, little jog and slope when there are $n = 3t$ columns to the right of the dotted line. It 
is easy to see that $N_1(t)$ and $N_2(t)$ remain the same if we count tilings of their flips instead. 
Our goal is to find $N(t)$. In order to express these, we write them as 
generating functions. 

\begin{eqnarray*}
G(z) & = & \sum_{t}N(t)z^{t}
\end{eqnarray*}  

From Figure 7, we see that $R(5,3t)$ can be obtained from little jog in 8 ways (4 from its flipped 
counterpart), and it can be obtained from slope in 4 ways (2 from its flipped counterpart). So, 
we have the following equation:

\begin{eqnarray}
G(z) & = & 8zG_1(z) + 4zG_2(z)
\end{eqnarray} 

In \cite{moore} Moore had derived the generating functions $G_1(z)$ and $G_2(z)$ as

\begin{eqnarray}
G_1(z) & = & \frac{8z(1+2z+5z^{2})}{1-2z-31z^{2}-40z^{3}-20z^{4}}
\end{eqnarray}

\begin{eqnarray}
G_2(z) & = & \frac{2z(1+14z+5z^{2})}{1-2z-31z^{2}-40z^{3}-20z^{4}}
\end{eqnarray}

Substituting $G_1(z)$ and $G_2(z)$ in equation (1), we get the final generating
function for exactly counting the number of faultfree tilings of $R(5,3t)$, $t\geq 2$.

\begin{eqnarray}
G(z) & = & \frac{24z^{2}(3+10z+15z^{2})}{1-2z-31z^{2}-40z^{3}-20z^{4}}  
\end{eqnarray}

The coefficients $N(t)$, starting with $t=2$, of the first few terms of $G$'s Taylor expansion are 

\begin{eqnarray*}
72, 384, 3360, 21504, 163968, 1136640, 8283648, 58791936, 423121920, ...
\end{eqnarray*}

\section{A general scheme for counting faultfree tilings of $R(m,n)$}

We now present an approach to count the total number of faultfree tilings of the rectangle $R(m,n)$ 
where $m,n\geq4$ and $3|mn$. Our approach is similar to that used in Section 4; we 
count the total number of ways of extending a smaller faultfree tiling to a bigger one 
and multiply with the total number of faultfree tilings of the smaller one to get the 
total number of faultfree tilings of the bigger faultfree tiling. Each of the patterns 
depicted in Figure \ref{patterns} is expanded along its rows by 
one of its generators, 
independent of expansions in other rows of the rectangle. 
So, the portion of the extended tilings of the rectangle 
due to distinct instances of generators cannot share a tromino. Also
note that all tiling extensions done in this section are in the 
horizontal direction.      
Naturally, tilings obtained in such a scheme cannot exhaust the
full set of possible tilings of the final rectangle. Therefore, 
we get at best a lower bound estimate of the number of tilings of 
$R(m,n)$. 

\subsection{Faultfree tilings of $6\times 6t$ rectangles}

We start with the case $m=6, n=6t$, where $t\geq1$. We have the following lemma.   
\begin{figure}[htbp]
\centerline{
{\scalebox{.15}{\includegraphics{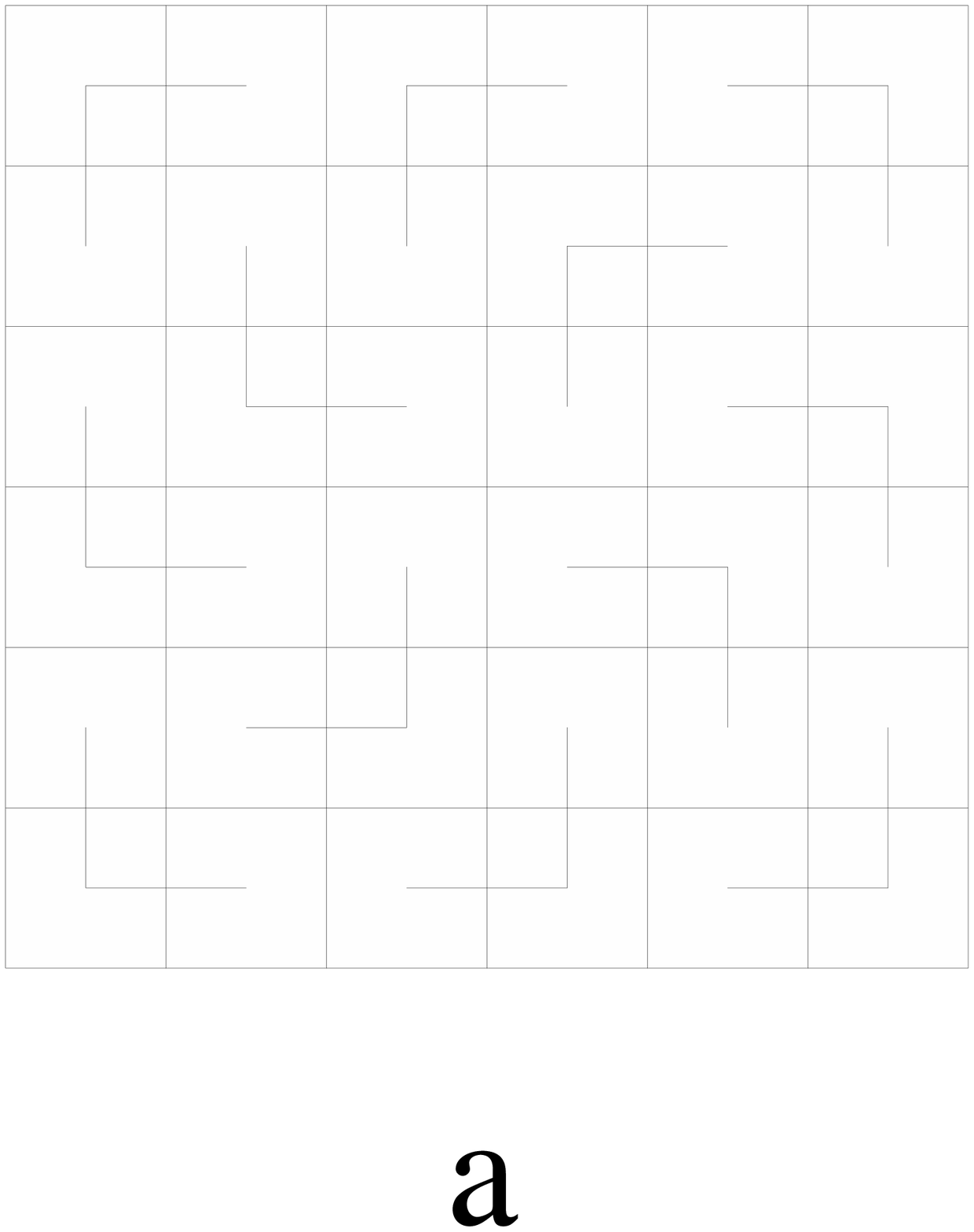}}}
{\scalebox{.15}{\includegraphics{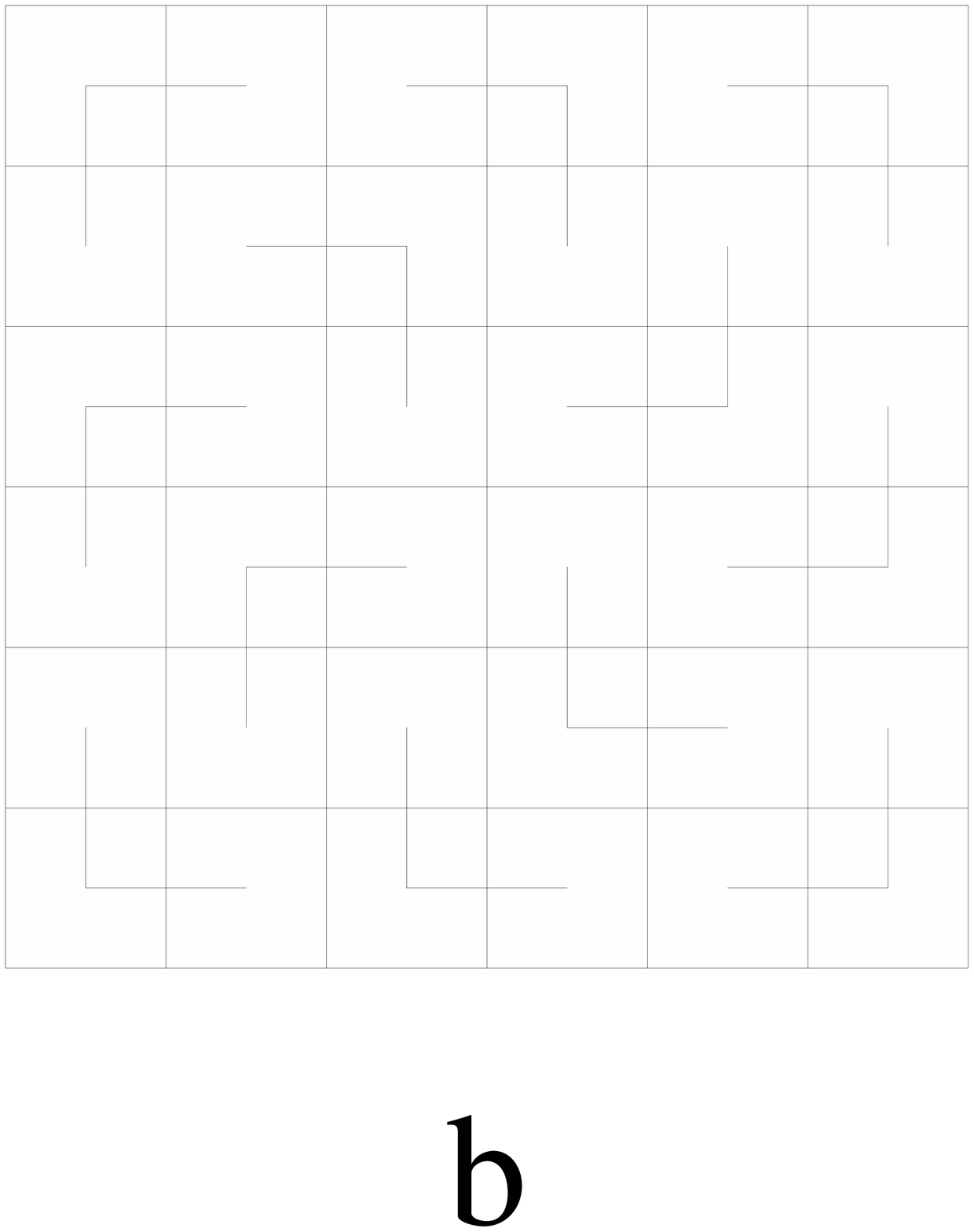}}}
}
\caption{Only possible ways of tiling $R(6,6)$ faultfreely.}
\end{figure}  
\begin{lemma}
The only possible ways of tiling $R(6,6)$ faultfreely are those shown 
in Figure 8.
\end{lemma} 
\begin{proof}
From Lemma 2, we know that at 
least one of the four patterns, shown in Figure 2(i), will appear in the leftmost three columns in any 
faultfree tromino tiling of $R(6,6)$. The patterns shown in Figure 3(e) 
and 3(c) do not arise as they make respectively, $(1,1)$ or $(6,1)$ inaccessible.  
Suppose the pattern shown in Figure 2(d) is present 
in the first five rows. Consider the tromino covering $(5,1)$ and $(5,2)$. 
If this tromino covers $(6,2)$ then $(6,1)$ becomes inaccessible, and if it covers 
$(6,1)$ then $(6,2)$ becomes inaccesible. So we conclude that this pattern too  
does not arise. Now consider Figure 2(b) and suppose that it is present in the 
first two rows. The tromino covering $(1,4)$ has two permissible orientations. If it 
(i) covers $(1,5)$ and $(2,4)$ then a $R(2,3)$ is completed by 
the tromino covering $(1,6)$, and if it (ii) covers $(2,4)$ and $(2,5)$, then a
$R(2,3)$ is completed by the tromino covering 
$(1,5)$. In either case, the $3rd$ horizontal grid line becomes 
a fault line. So this pattern 
cannot appear in the first two rows. The only other possibility is that this pattern appears 
in rows 3 and 4. In this case, if the tromino (i) covering $(1,1)$ and $(2,1)$ covers $(1,2)$, 
then a $R(2,3)$ is completed by the tromino covering 
$(2,2)$, and if it (ii) covers $(2,2)$ then $R(2,3)$ is completed by 
the tromino covering $(1,2)$. 
Applying a similar reasoning to the tromino covering $(5,1)$ and $(6,1)$, the reader can 
easily see that the $4th$ vertical grid line will become a fault line. So, we conclude that 
this pattern too does not arise. 

Now consider the pattern shown in Figure 2(a). This pattern 
arises either in the first four rows or the 
last four rows. Otherwise, either
$(1,1)$ or $(6,1)$ will become inaccessible. Without loss of generality we  
assume that this pattern arises in the first 
four rows. Since Figure 2(b) does not arise, the tromino covering $(5,1)$ and $(6,1)$ must 
cover $(6,2)$ and the tromino covering $(5,2)$ and $(5,3)$ must cover $(4,3)$. Instead of 
covering $(3,3)$, if the critical tromino covers $(2,3)$, then it is easy to see that the 
trominoes covering $(1,3)$, $(3,3)$ and $(6,3)$ will make the $5th$ vertical grid line a fault line. 
So we know that the critical tromino will cover $(3,3)$. Due to symmetry, the rightmost three 
columns will also have the same tiling as above, whence the tiling in Figure 8(a). If the 
pattern in Figure 2(a) appears in the last four rows, then following similar reasoning the 
reader can verify that we will get the tiling shown in Figure 8(b). So we conclude that Figure 
8 shows all possible ways of tiling $R(6,6)$ faultfreely.       \hfill \qed
\end{proof}

We now proceed to derive a lower bound on the number of faultfree tilings of $R(6,6t)$. 
From Lemma 7, we know all possible ways of tiling $R(6,6)$ faultfreely. Each of these need be 
extended in only two directions, left and right, by the generative scheme suggested below;
due to rotational symmetry or $R(6,6)$, we do not need to separately consider extensions of the top 
and the bottom. We first state the following lemma. In 
the scheme suggested below, we extend either only on the left, or only on the right, or on both sides, by six columns.

\begin{figure}[htbp]
\centerline{\scalebox{.3}{\includegraphics{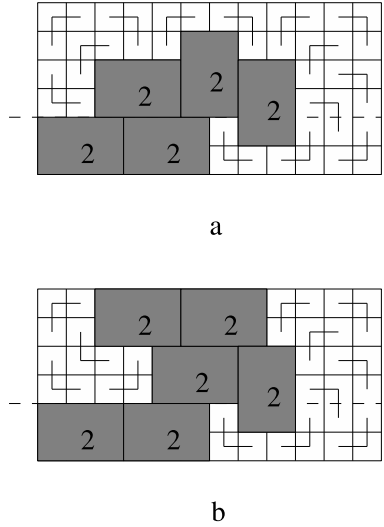}}}
\caption{Number of ways of expanding the pattern of Figure 3(a).}
\end{figure}
\begin{lemma}
The only possible ways of extending the pattern shown in Figure 3(a), by six columns,  
are those depicted in Figure 9.
\end{lemma}
\begin{proof}
As has been proved in Section 4, in any tiling of a rectangle $R(6,12)$, the tiling of 
the leftmost four columns and topmost four rows 
of each of the generators shown in Figure 9 must match the pattern shown 
in Figure 6(a) or its flipped counterpart. We also claim that in any 
faultfree extension of the pattern in Figure 3(a) the tromino covering 
the squares $(3,8)$ and $(3,9)$ must cover $(4,8)$. If this is not 
the case then it must cover $(2,8)$. But then the trominoes covering the 
squares $(1,8)$ and $(4,8)$ will make the $7th$ vertical grid line a fault 
line (since the remaining two rows are simply padded with $R(2,3)$s), 
whence the given orientation of the tromino. Consider the case when 
a tiling of the leftmost four columns of the generator matches the pattern 
shown in Figure 6(a). Following Lemma 4 the remaining area can only be tiled 
as shown in Figure 9(b). Now consider the case when a tiling of the 
leftmost four columns matches the flipped counterpart of the pattern shown 
in Figure 6(a). Again following Lemma 4 the remaining area can only be tiled 
as shown in Figure 9(a). Thus Figure 9 shows all possible ways of extending the 
pattern shown in Figure 3(a).   
  
Note that the tromino covering $(4,9)$, $(5,9)$ and $(5,8)$ completes 
$R(3,2)$ with the tromino covering $(3,9)$, $(3,8)$ and $(4,8)$. Thus, 
the total number of ways of extending $R(6,6)$ to $R(6,12)$ is double 
of that deduced above  since this $R(3,2)$ can be tiled in two ways.   \hfill \qed
\end{proof}  

First consider extending $R(6,6)$ only on one side (either entirely to 
the left, or entirely to the right). Without loss of 
generality, we assume that extension is done only to the left side. 
We use Lemma 8 for the very first extension of the top (or bottom) 
four rows by six columns. There are in all (i) 2 ways of tiling $R(6,6)$, 
(ii) 24 ways of extending the top or bottom four rows by six columns, and 
(iii) 4 ways of padding the remaining two rows with two $R(2,3)$ rectangles.  
So we get a total of $2\times24\times4=192$ ways of extending $R(6,6)$ to $R(6,12)$.  
Observe that we generate the pattern shown in Figure 6(a) (in the top or bottom 
four rows) after the first extension. For subsequent extensions by six columns,
we invoke any of the 36 possible generators in Figure 6 for four top or bottom rows. 
We always extend the remaining two rows by padding with $R(2,3)$ rectangles. So 
we have a total of $192.(36\times4)^{t-2}=192.144^{t-2}$, where 
$t\geq2$, ways of extending $R(6,6)$ to $R(6,6t)$ when extension is done only on the left side. 
Due to symmetry, there are $192.144^{t-2}$ ways of extending $R(6,6)$ to $R(6,6t)$ when 
extension is done only on the right side. So we conclude that there are in all 
$2\times192.144^{t-2} = 384.144^{t-2}$ ways of extending $R(6,6)$ to $R(6,6t)$ when extension 
is done only on one side (either entirely to the left, or entirely to the right). 

Now consider extension of $R(6,6)$ by $6r$ columns on one side and 
$6s$ columns on the other side,
where $r,s\geq1$, $6(r+s+1)=6t$, and $t\geq3$. Again, the first extension of 
$R(6,6)$ to $R(6,12)$ is by Lemma 8, and the remaining extensions use generators 
in Figure 6, along with padding with $R(2,3)$ rectangles. 
So, the total number of ways of such an extension on either side is
\begin{eqnarray}
M(t) & = & \sum_{r+s+1=t}(24\times4).(36\times4)^{r-1}.2.(24\times4).(36\times4)^{s-1} \nonumber \\
     & = & 18432\sum_{r+s+1=t}144^{r+s-2} \nonumber   \\
     & = & 18432\sum_{r+s+1=t}144^{t-3} \nonumber  \\
\end{eqnarray} 

The number of times we sum $(144)^{t-3}$ is the number of 
solutions of the equation $r+s+1=t$, 
where $r,s\geq1$. The number of solutions of this 
equation is the coefficient of $x^{t-1}$ in 
$\frac{x^2}{(1-x)^2}$. This comes out to be $^{2+(t-3)-1}C_{t-3}=t-2$. 
So, we conclude that the number of ways of tiling $R(6,6t)$ faultfreely is atleast
$N(t)=M(t)+384.144^{t-2}=
18432.(t-2).144^{t-3} + 384.144^{t-2}=128.(t+1).(144)^{t-2}$, where $t\geq3$. This formula 
also holds for the case $t=2$. We summarize our result in the following theorem. 
 
\begin{theorem} 
The number of ways of tiling $R(6,6t)$, where $t\geq2$, faultfreely is atleast 
\begin{eqnarray}
N(t) & = & 128.(t+1).144^{t-2}               
\end{eqnarray}      
\end{theorem} 
  




Suppose $Q(z) = \sum_{t}q(t)z^t$ is the generating function where the coefficient 
of $z^t$ ($q(t) = 384.144^{t-2}$) gives the number of ways of tiling $R(6,6t)$ faultfreely, 
where $t\geq2$, when $R(6,6)$ is extended only on one side (either entirely to the left or 
entirely to the right). We can use $Q(z)$ to find the generating function 
$F(z) = \sum_{t}f(t)z^t$, where the coefficient of $z^t$ ($f(t) = 128.(t+1).144^{t-2}$) gives 
the total number of ways of tiling $R(6,6t)$ faultfreely when $R(6,6)$ is extended only on one 
side (either entirely to the left or entirely to the right) or on both sides. We state the 
following fundamental result as a theorem. 

\begin{theorem}
If $c$ is the total number of ways of tiling $R(6,6)$ faultfreely, then
\begin{eqnarray}
F(z) & = & \frac{Q^2(z)}{4cz} + Q(z)
\end{eqnarray} 
\end{theorem}  
\begin{proof}
Counting for the case of permissible extensions on both sides may be
viewed as convolution of extensions on only one side (either entirely 
to the left or entirely to the right). The first term is due to convolution; 
since $F(z)$ considers extensions on both left and 
right sides, we divide the square of the one-sided generating function $Q(z)^2$ 
by $2^2=4$. Also, the base rectangle 
$R(6,6)$ is considered twice, once by each $Q(z)$ term. So we merge the two 
rectangles into one, decreasing the number of columns by six and the number of 
tilings of the initial $R(6,6)$ by $c$. Therefore we divide $Q^2(z)$
by $4cz$. Note that the first term in $Q(z)$ is a multiple of $z^2$. So when $Q(z)$ is convoluted, we 
get a generating function which considers extension of the base rectangle $R(6,6)$ by 
$6r$ columns on one side and $6s$ columns on the other side, such that $6(r+s+1)=6t$, 
where $r,s\geq1$. For incorporating the cases where either $r=0$ or $s=0$, we add $Q(z)$ to this 
convolution.     \hfill \qed
\end{proof}

Once we have $Q(z)$, and we can use Theorem 7 to determine $F(z)$. We demonstrate 
the usefulness of the above theorem for the case of tiling $R(6,6t)$, where $t\geq2$, faultfreely. 
Here we have $Q(z)$ as the generating function where the coefficient of $z^t$ ($q(t)$) is 
$384.144^{t-2}$. Then we have

\begin{eqnarray}
Q(z) & = & \frac{384z^2}{(1-144z)} 
\end{eqnarray}

Now, by Theorem 7, we have 

\begin{eqnarray}
F(z) & = & \frac{384^2z^4}{(1-144z)^2\times4\times2\times z} + \frac{384z^2}{(1-144z)} \nonumber \\
     & = & \frac{384z^2}{(1-144z)}\{\frac{384z}{(1-144z)\times8} + 1\}  \nonumber \\
     & = & \frac{384z^2}{(1-144z)}\{\frac{48z+1-144z}{(1-144z)}\} \nonumber \\
     & = & \frac{384z^2(1-96z)}{(1-144z)^2}
\end{eqnarray}     

The coefficient of $z^t$ in $F(z)$ is 
\begin{eqnarray}
f(t) & = & 384\times144^{t-2}\times^{2+(t-2)-1}C_{t-2} - 96\times384\times144^{t-3}\times^{2+(t-3)-1}C_{t-3} \nonumber \\
     & = & 384\times144^{t-2}\times(t-1) - 96\times384\times144^{t-3}\times(t-2) \nonumber \\
     & = & 384\times144^{t-3}\times{48t+48} \nonumber \\
     & = & 128.(t+1).144^{t-2}
\end{eqnarray}

This is exactly the same formula as derived in Theorem 6 (equation (8)). Thus, 
Theorem 7 establishes an elegant method for deriving the generating function $F(z)$ 
for tiling extensions on only one side (either entirely to the left or entirely to 
the right) or on both sides, from $Q(z)$ which is for tiling extensions on only one side 
(either entirely to the left or entirely to the right). 

\subsection{Faultfree tilings of $R(7,6t)$}

Now we consider the case $m=7, n=6t$, where $t\geq1$. We have the following lemma.   
\begin{figure}[htbp]
\centerline{
{\scalebox{.15}{\includegraphics{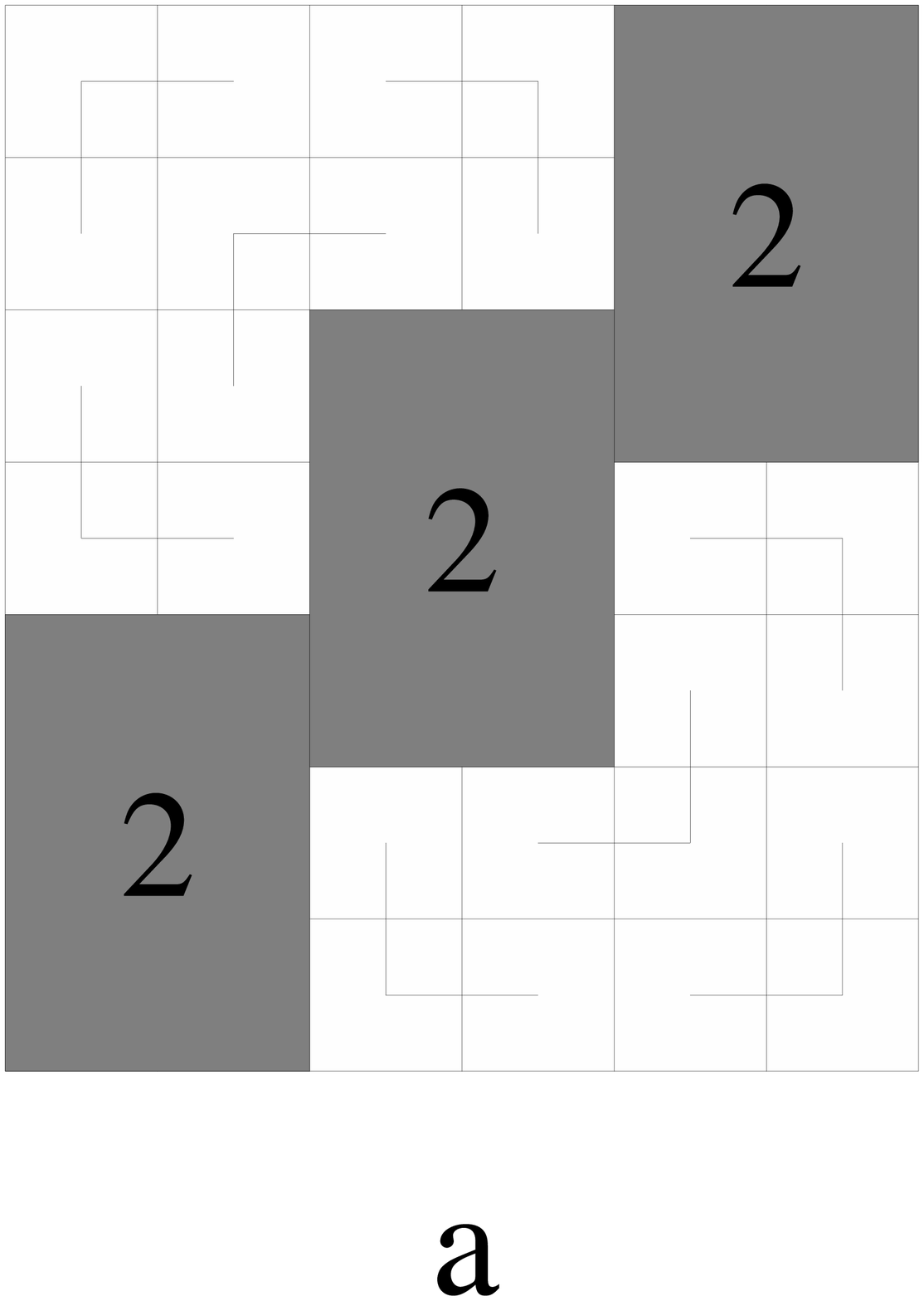}}}
{\scalebox{.15}{\includegraphics{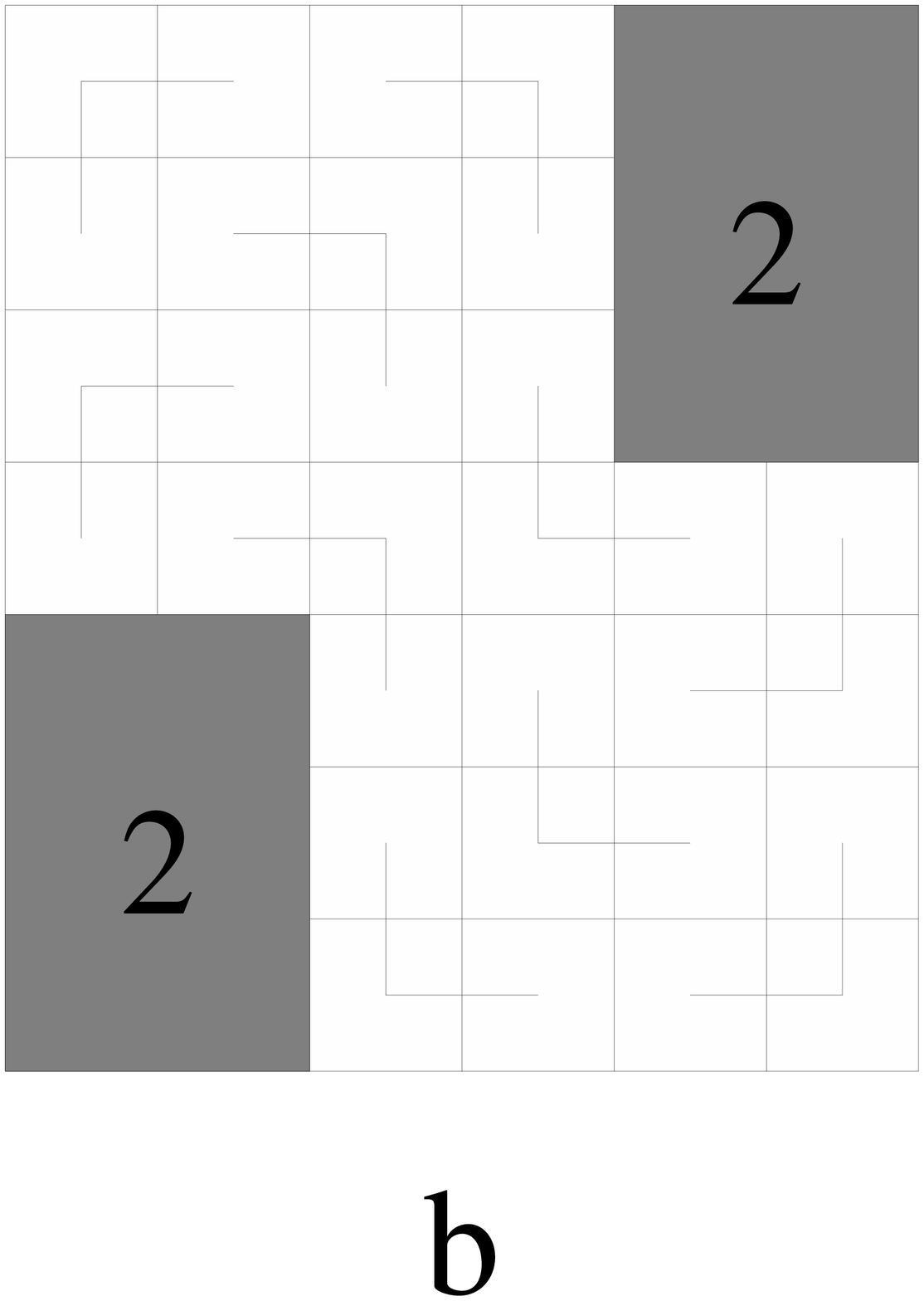}}}
}
\caption{Only possible ways of tiling $R(7,6)$ faultfreely.}
\end{figure}  
\begin{lemma}
The only possible ways of tiling $R(7,6)$ faultfreely are those shown in Figure 10 and their flipped 
counterparts. 
\end{lemma}
\begin{proof}
From Lemma 2, at least one of the patterns in Figure 2 will arise in the leftmost three columns 
of $R(7,6)$. If the patterns in Figure 3(c) or 3(e) are present in the leftmost three columns, 
then they will be either in the topmost five rows or the bottommost five rows. They can't occupy 
rows 2-6 as then the squares $(1,1)$ and $(7,1)$ will become inaccessible.
  
First consider that the pattern in Figure 3(c) arises in the bottommost five rows. 
If the tromino covering $(1,1)$ and $(2,1)$ covers $(2,2)$ then $R(2,3)$ is completed by the 
tromino covering $(1,2)$ and if it covers $(1,2)$ then $R(2,3)$ is completed by the tromino 
covering $(2,2)$. So, we have a $R(2,3)$ at the top of the pattern in Figure 3(c). 
Consider the tromino covering $(1,4)$ and $(2,4)$ which has two permissible orientations.  
If this tromino covers $(2,5)$, then another $R(2,3)$ is completed by the tromino covering $(1,5)$ 
and if it covers $(1,5)$, then $(2,3)$ is completed by the tromino covering $(1,6)$. so, the 
$3rd$ horizontal grid line becomes a fault line. So, we conslude that the pattern in Figure 3(c) 
cannot be present in the bottommost five rows. Suppose it is present in the topmost five rows. 
The reader can easily see that we cannot have $R(2,3)$ at the bottom as then the $6th$ horizontal 
grid line will become a fault line. Then the only possible case is that the tromino covering $(6,1)$ 
and $(7,1)$ covers $(7,2)$ and the tromino covering $(6,2)$ and $(6,3)$ covers $(5,3)$. It is easy 
to see that the tromino covering $(7,3)$ covers $(7,4)$ and $(6,4)$. Consider the case when the 
tromino covering $(7,5)$ and $(7,6)$ completes a $R(3,2)$ with the tromino covering $(5,5)$. If the 
tromino covering $(5,4)$ and $(4,4)$ covers $(4,5)$ then the trominoes covering $(4,3)$ and $(4,6)$ 
make the $3rd$ horizontal grid line a fault line. If it covers $(4,3)$ then the trominoes covering 
$(3,3)$ and $(1,4)$ make $(1,6)$ inaccessible. The only other possibility is that the tromino  
covering $(7,5)$ and $(7,6)$ covers $(6,6)$ and the tromino covering $(6,5)$ and $(5,5)$ covers 
$(5,4)$. In this case the trominoes covering $(5,6)$, $(4,6)$, $(4,5)$ and $(3,3)$, $(4,3)$, $(4,4)$ 
leave an untileable $R(3,3)$. We conclude that the pattern in Figure 3(c) cannot be present in 
the leftmost three columns of $R(7,6)$. Due to symmetry, it cannot be present in the rightmost three 
columns also. For this same reason, the pattern in Figure 3(e) cannot be present in the topmost five 
rows. 
   
Now consider the case where the pattern in Figure 3(e) is present in the bottommost five rows. If the 
tromino covering $(1,1)$ and $(2,1)$ covers $(2,2)$ then $R(2,3)$ is completed by the tromino covering 
$(1,2)$, and if it covers $(1,2)$ then $R(2,3)$ is completed by the tromino covering $(2,2)$. For 
reasons similar to above, a $R(2,3)$ will also be present in the first two rows from columns 4 to 6. 
But then the $3rd$ vertical grid line will become a fault line. So, we conclude that the pattern in 
Figure 3(e) will not be present also. 
         
Consider the case when the pattern in Figure 2(d) is present in the first four rows. If the tromino 
covering $(5,1)$ and $(5,2)$ covers $(6,1)$ then $R(3,2)$ is completed by the tromino covering 
$(7,1)$ and if it covers $(6,2)$ then $R(3,2)$ is completed by the tromino covering $(6,1)$. Due, to 
symmetry, the rightmost three columns will have this same pattern. Thus, we get the tiling in Figure 
10(b). 
  
Consider the case when the pattern in Figure 2(a) is present in the tiling of the leftmost three 
columns. This pattern can be present only in the topmost or bottommost four rows as in any other 
case either $(1,1)$ or $(7,1)$ will become inaccessible. This same pattern will be present in the 
rightmost three columns, since the pattern in Figure 2(d) only occurs with itself, as shown above. 
If the pattern in Figure 2(a) occurs in the topmost right three columns, then the tromino covering 
$(2,5)$ and $(3,5)$ will cover $(2,4)$ and the reader can easily see that the $5th$ horizontal grid 
line will become a fault line. If the pattern in Figure 2(a) occurs in the bottommost three columns 
with the tromino that covers $(5,5)$ and $(6,5)$ also covering $(6,4)$, then it can be verified that 
either the squares $(1,3)$, $(1,4)$ or the squares $(5,3)$, $(5,4)$ will become untileable. So the 
pattern in Figure 2(a) can only occur in the first four rows in which the tromino covering $(2,2)$ and 
$(3,2)$ also covers $(2,3)$. Due to symmetry, this same pattern will also be present in the rightmost 
three columns, giving rise to the tiling shown in Figure 10(a). So Figure 10 shows all possible ways 
of tiling $R(7,6)$. \hfill \qed
\end{proof}

\begin{figure}[htbp]
\centerline{
{\scalebox{.25}{\includegraphics{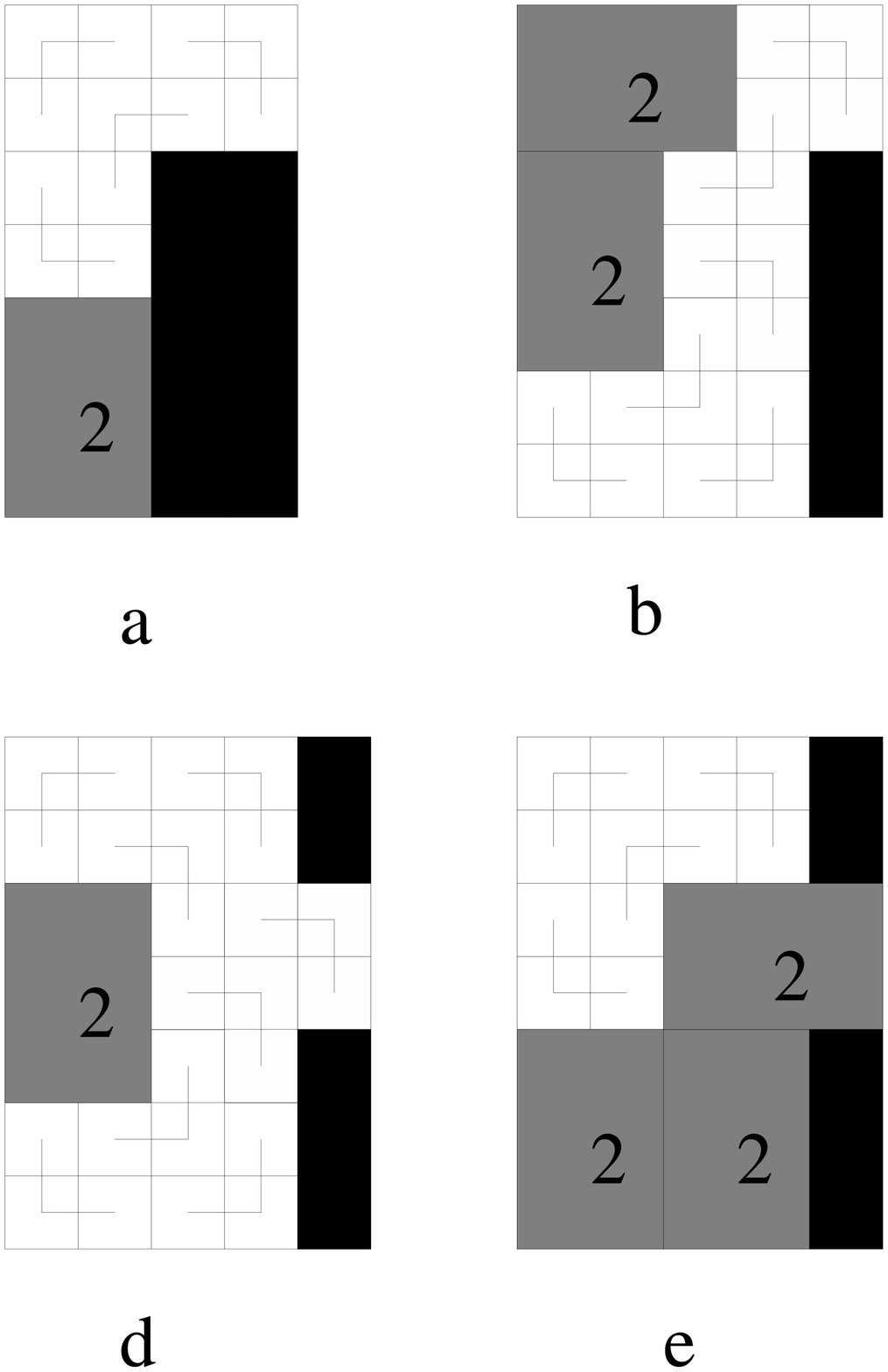}}}
{\scalebox{.25}{\includegraphics{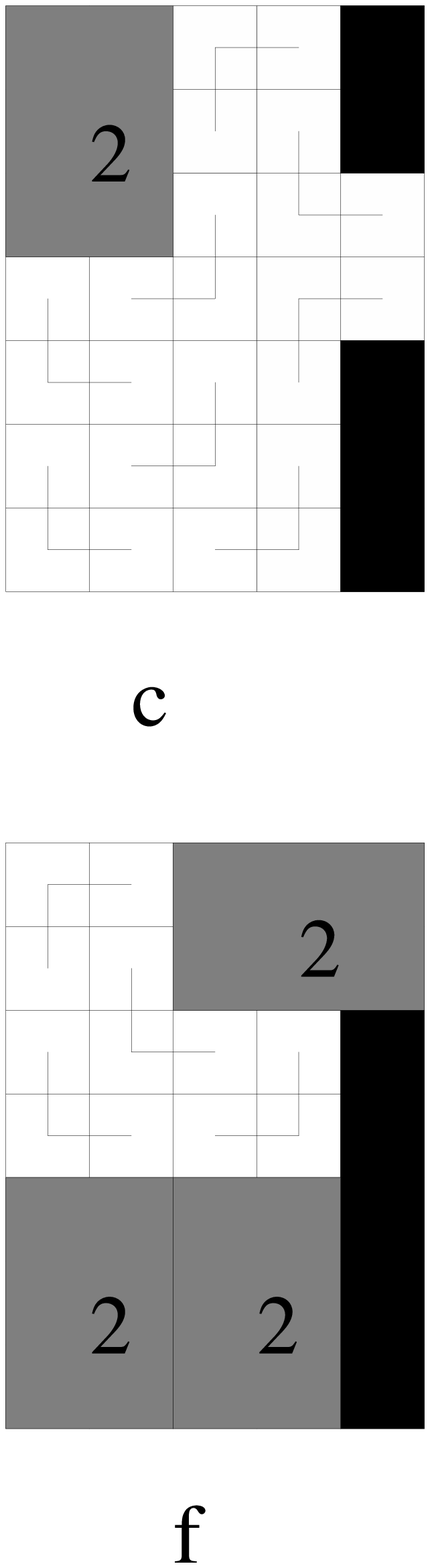}}}
}
\caption{Possible tilings of the leftmost five columns of $R(7,6t)$.}
\end{figure}  

Let us consider the pattern in Figure 10(a). In this case, Lemma 9 shows that there are only 
16 ways of tiling $R(7,6)$ faultfreely. Note that the pattern 
shown in Figure 5(a) is present either in the first four or the last four rows. So for the very 
first extension, we use the generators of Figure 5. We extend the remaining three rows by padding 
with $R(3,2)$ rectangles. In this process, the tiling of the leftmost (rightmost) five columns 
in the resulting tilings of $R(7,12)$ matches the 
patterns in Figures 11(e) or (f). 
We can choose to further extend these patterns by the generators in Figure 6. We also have
the option of retiling the leftmost (rightmost) five columns for tilings with patterns matching 
Figure 11(e) and (f), into patterns in Figure 12. So, for $R(7,12)$, we get distinct patterns in
(i) Figure 12 (c)-(e) from Figure 11(e), and (ii) Figure 12(a) and (b) from Figure 11(f). 
We generate distinct faultfree tilings of $R(7,12)$ using either the generators in
Figures 5 and 6, or the retilings leading to patterns in Figure 12. 
The reader should note that we consider extension of tilings of $R(7,6t)$ only on one side 
(either entirely to the left or entirely to the right), since we intend to derive a similar analogue 
of the generating function $Q(z)$ (see Theorem 7) for tilings of $R(7,6t)$, and further use it and 
Theorem 7 to derive the analogue of the generating function $F(z)$ for faultfree tilings of $R(7,6t)$.    
Now consider tilings of $R(7,6t)$, $t\geq 3$. From the retiled 
patterns in Figure 12, we can get the other four patterns in Figure 11 (a)-(d); 
now we may choose to extend tilings matching the patterns in Figure 11(b) and (d) using Moore's
approach, those in Figure 11(e) and (f) by our usual generators of Figure 6, and that in Figure 
11(a) by our generators of Figure 5. We do not extend tilings ending
in the pattern of Figure 11(c). 
Using this scheme we generate faultfree tilings of $R(7,6t)$, $t\geq 3$, iteratively.

\begin{figure}[htbp]
\centerline{
{\scalebox{.25}{\includegraphics{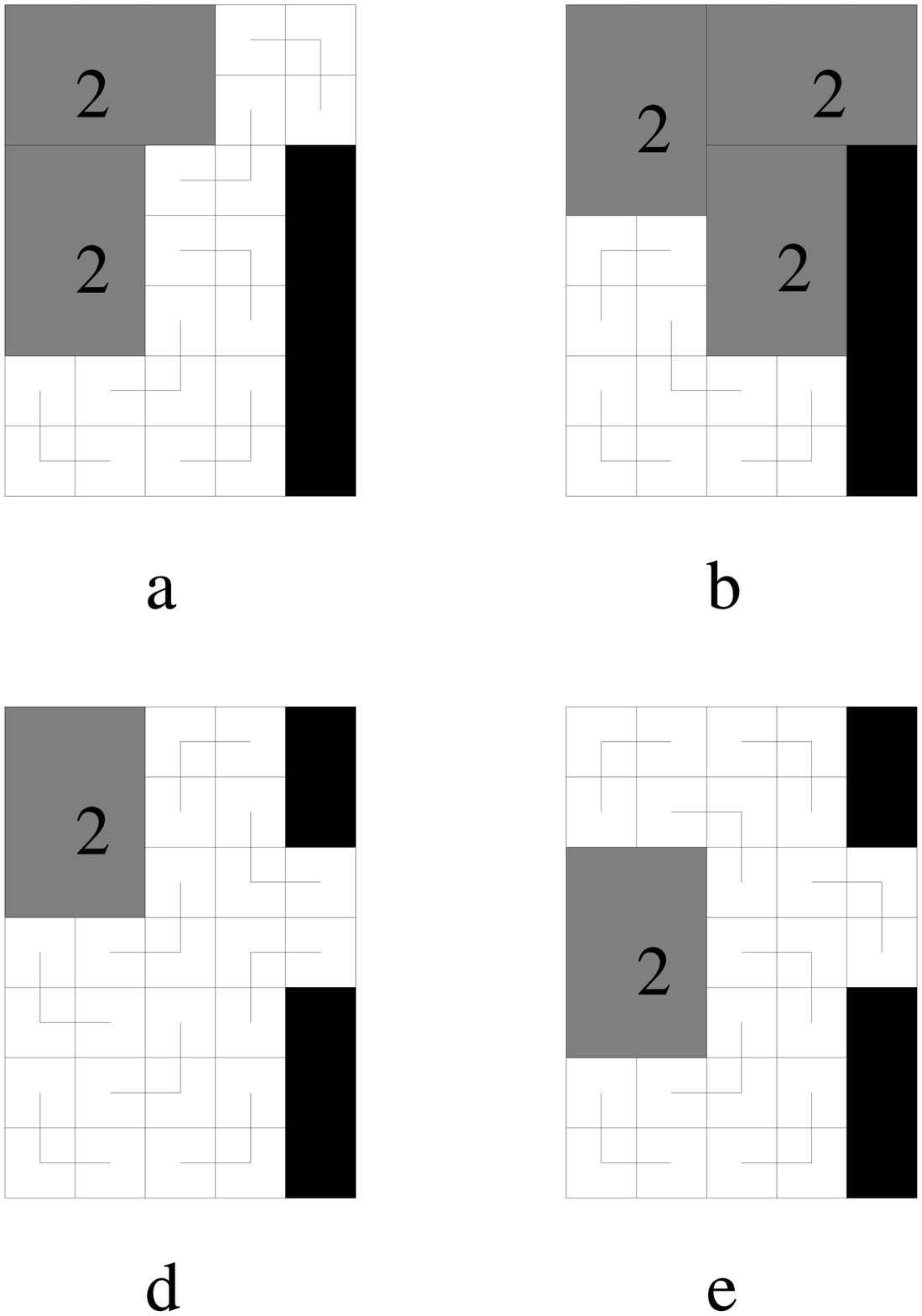}}}
{\scalebox{.25}{\includegraphics{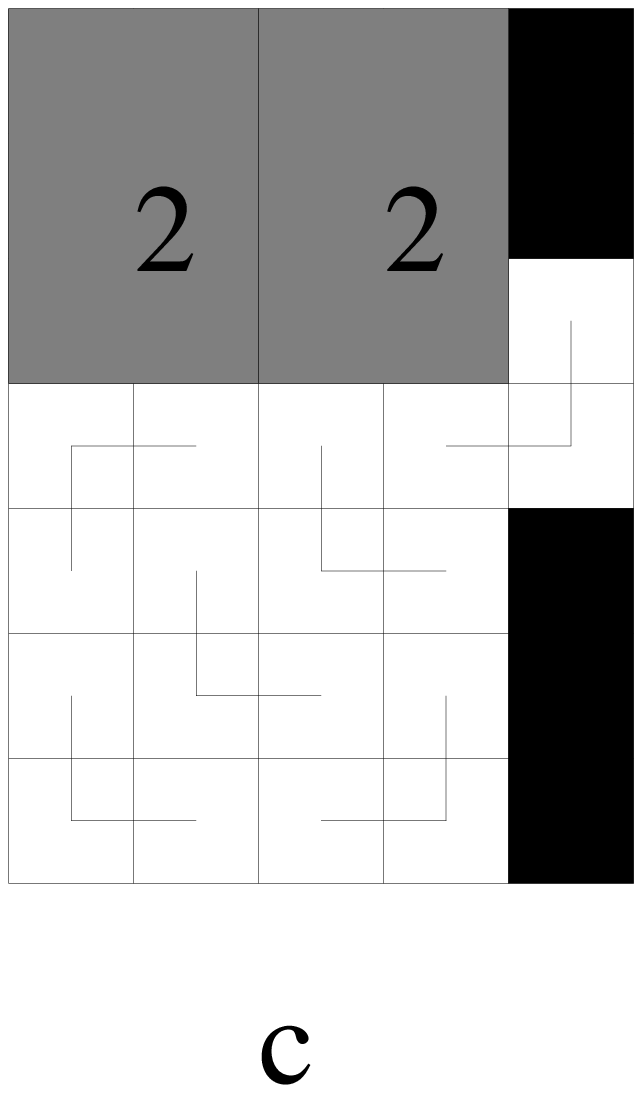}}}
}
\caption{Only possible ways of retiling Figures $11(e)$ and $11(f)$.}
\end{figure} 

\begin{lemma}
The only possible retilings of the patterns in Figure 11(e) and 11(f) 
preserving faultfree-ness are those shown in Figure 12.
\end{lemma} 
\begin{proof}
First consider retiling the pattern in Figure 11(f). Consider the two possible orientations of the 
tromino covering $(1,5)$ and $(2,5)$. If it covers $(2,4)$ then a $R(2,3)$ is completed by the tromino 
covering $(1,4)$. If it covers $(1,4)$, then the four possible orientations. Suppose that it covers 
$(1,3)$ and $(2,3)$ thereby completing a $R(2,3)$. Consider the tromino covering $(1,1)$ and $(1,2)$. 
If it covers $(2,2)$, then a $R(3,2)$ is completed by the tromino covering $(2,1)$. If it covers $(2,1)$, 
then the tromino covering $(2,2)$ and $(3,2)$ has two possible orientations. If it covers $(3,3)$, then it 
is easy to see that we get back the tiling in Figure 11(f). So it must cover $(3,1)$. Now from Lemma 2, the 
tiling of the bottommost four rows must match the pattern in Figure 2(a) or its symmetric counterpart. 
Note that the tromino covering $(5,2)$ and $(6,2)$ cannot cover $(5,3)$ as then the tromino covering 
$(5,4)$ will make either $(6,3)$ and $(6,4)$, or $(3,3)$ and $(3,4)$ untileable. The remaining area can 
now only be tiled in two ways as shown in Figure 12(b). In two of the other possible orientations of the 
tromino covering $(2,4)$ (1. it covers $(2,3)$ and $(3,4)$, and 2. it covers $(2,3)$ and $(3,3)$), the 
tromino covering $(1,3)$ makes $(1,1)$ inaccessible. So these two orientations are not possible. Consider 
the case when the tromino covering $(2,4)$ covers $(3,4)$ and $(3,3)$. If the tromino covering $(1,3)$ and 
$(2,3)$ covers $(1,2)$ then a $R(2,3)$ is completed by the tromino covering $(1,1)$ and if it covers $(2,2)$, 
then a $R(2,3)$ is completed by the tromino covering $(1,2)$. So in both these cases, a $R(2,3)$ is present 
in the first two rows. Consider the tromino covering $(3,1)$ and $(3,2)$. If it covers $(4,1)$, then the 
remaining area is $R(4,4)^{-}$ which has a unique tiling, so we get the tiling in Figure 12(a). If it covers 
$(4,2)$, then a $R(3,2)$ is completed by the tromino covering $(4,1)$. The remaining area can now only be 
covered as shown in Figure 12(a).  
 
Now consider retiling the pattern in Figure 11(e). Suppose that a tromino covers $(3,5)$, $(3,4)$, $(2,4)$ 
and another tromino covers $(4,5)$, $(4,4)$, $(5,4)$. Note that the tromino covering $(1,4)$ must cover 
$(1,3)$ and $(2,3)$. Consider the tromino covering $(6,4)$ and $(7,4)$. If it covers $(6,3)$ then a 
$R(2,3)$ is completed by the tromino covering $(7,3)$ in which case $(7,1)$ becomes inaccessible. So it 
must cover $(7,3)$. For similar reasons as above, the tromino covering $(6,3)$ will not cover $(6,2)$ and 
$(7,2)$. If the tromino covering $(7,1)$ and $(7,2)$ covers $(6,2)$ then a $R(3,2)$ is completed by the 
tromino covering $(6,1)$, in which case $(4,1)$ becomes inaccessible. So it must cover $(6,1)$. If the 
tromino covering $(6,2)$ and $(6,3)$ covers $(5,2)$, then it can be easily seen that the tromino covering 
$(5,1)$ will make $(5,3)$ inaccessible. So it must cover $(5,3)$. From Lemma 2, the remaining area can only 
be tiled by the pattern in Figure 2(c), whence the tiling in Figure 12(d). Consider the case when a 
single tromino covers $(3,4)$ and $(4,4)$. Suppose that it covers $(3,4)$. The tromino covering $(4,4)$ has 
four possible orientations. If it covers $(4,3)$ and $(5,3)$, then a $R(3,2)$ is completed by the tromino 
covering $(5,4)$, making $(7,4)$ inaccessible. If it covers $(5,4)$ and $(5,3)$, then either $(1,1)$ or 
$(7,1)$ becomes inaccessible. If it covers $(4,3)$ and $(3,3)$, then we get back the tiling in Figure 11(e). 
So we conclude that it must cover $(4,3)$ and $(5,4)$. Due to symmetry, the remaining area can now only be 
covered as shown in Figure 12(e). Consider the final case when the tromino covering $(3,5)$ and $(4,5)$ covers 
$(4,4)$. Consider the five possible orientations of the tromino covering $(4,3)$. If it covers $(3,3)$ and 
$(3,4)$, then it can be easily seen that we get back the tiling in Figure 11(e). If it covers $(4,2)$ and 
$(5,3)$ then the tromino covering $(5,4)$ makes $(7,4)$ inaccessible. If it covers $(5,3)$ and $(5,2)$, then 
$(7,4)$ becomes inaccessible for similar reasons as above. If it covers $(4,2)$ and $(5,2)$, the the tromino 
covering $(5,1)$ makes $(7,1)$ inaccessible. So we conclude that it must cover $(5,3)$ and $(5,4)$. Consider 
the tromino covering $(1,4)$. It cannot cover $(2,3)$ and $(2,4)$, otherwise $(3,4)$ will become inaccessible. 
If it covers $(1,3)$ and $(2,4)$, then a $R(3,2)$ is completed by the tromino 
covering $(3,4)$, and if it covers $(1,3)$ and $(2,3)$, then a $R(3,2)$ is completed by the tromino covering 
$(2,4)$. So in both cases we get a $R(3,2)$. Now consider the tromino covering $(1,1)$ and $(1,2)$. If it covers 
$(2,1)$, then a $R(3,2)$ is completed by the tromino covering $(2,2)$, and if it covers $(2,2)$, then a $R(3,2)$ 
is completed by the tromino covering $(2,1)$. The remaining area can now only be covered as shown inn Figure 
12(c). Thus, Figure 12 shows all possible ways of retiling the patterns in Figures 11(e) and 11(f). \hfill \qed
\end{proof}

We claim that
all tilings generated thus are distinct, thereby yielding a lower bound on the number of faultfree 
tilings of $R(7,6t)$. We show this in the following lemma. 

\begin{lemma}
All tilings of $R(7,6t)$, where $t\geq3$, generated thus using the scheme of extension and retiling are distinct
\end{lemma}
\begin{proof}
Consider distinct tilings of $R(7,6t)$, $t\geq 3$ whose leftmost five columns
match the patterns in Figure 11(a), (e) and (f). We claim that 
these tilings are extended by 
six columns into unique tilings of $R(7,6t+6)$ using our generators in Figures 5 and 6 for the
top or bottom four rows, with the remaining three rows simply being padded with $R(3,2)$ rectangles.
To establish the distinctness of these extended tilings, observe that the 
tilings of the leftmost five columns again match the patterns in Figures 11(e) and (f); 
observing the first four columns, we note that the pattern in Figure 11 (a) is generated only 
by retiling the pattern in Figure 11(e) (resp., Figure 11(f))
into Figure 12(c) (resp., Figure 12(b)). Note that the pattern in Figure 
11(a) subsumes only patterns in Figure 12 (b)-(c). So, tilings generated for $R(7,6t+6)$
whose leftmost five columns match Figure 11(e) are distinct from those that match Figure
11(a).  \hfill \qed
\end{proof}

We can further retile the patterns of Figure 
11(e) and (f) in the leftmost (rightmost) five columns to get 
the patterns in Figure 12 (c)-(e) and Figure 12 (a)-(b), respectively.
We now derive the generating functions for the number of tilings of $R(7,6t)$ where the 
tiling of the leftmost five columns match the patterns shown in Figures 11 (a)-(f), and 
which are obtained by the approach discussed above, namely, extension by six columns 
using the generators of Figures 5 and 6, and then retiling the leftmost five columns as 
shown in Figure 12. The reader should again note that all extensions are done only on one 
side (either entirely to the left or entirely to the right), since we intend to derive the 
analogue of the generating function $Q(z)$ (ref. Theorem 7) for faultfree tilings of $R(7,6t)$. 
 
Let $H_1(t)$ ($J_1(t)$, $K_1(t)$, $P_1(t)$, $T_1(t)$, $S_1(t)$) 
denote the number of tilings of $R(7,6t)$, using the above approach, where the tilings of 
the leftmost five columns match the patterns shown in Figure 11(a) (11(b), (c), 
(d), (e), (f)), or their flipped counterparts respectively. 
We define the generating function $H(z)=\sum_{t=1}H_1(t)z^t$. Similarly, we 
define the generating functions $J(z)$, $K(z)$, $P(z)$, $T(z)$ and $S(z)$. Note that all 
the 16 faultfree tilings of $R(7,6)$ have their leftmost three columns matching the pattern 
in Figure 11(a). So the coefficient of $z$ in $H(z)$ is 16. The pattern in Figure 11(a) is 
obtained from that in 11(f) by retiling as shown in Figure 12(b). Here we disturb the 
three $R(3,2)$ rectangles in 11(f) and again construct three $R(3,2)$ rectangles as in Figure 12(b). So, 
corresponding to each tiling of $R(7,6t)$, where the leftmost five columns match the pattern 
in Figure 11(f), we get a tiling where the leftmost four columns match the pattern in Figure 11(a). Now consider 
retiling the pattern in Figure 11(e) to that in Figure 11(a). Here we disturb the three 
$R(3,2)$ rectangles in Figure 11(e) and construct back only two $R(3,2)$ rectangles as shown in Figure 12(c). 
So corresponding to every two tilings of $R(7,6t)$ where the leftmost five columns match the 
pattern in Figure 11(e), we get one tiling where the leftmost four columns match the pattern 
in Figure 11(a). So we get the following equation:

\begin{eqnarray}
H(z) & = & 16z + S(z) + \frac{1}{2}T(z)
\end{eqnarray}
        
Now consider all tilings of $R(7,6t)$ where the leftmost five columns match the pattern in 
Figure 11(f). We can get such tilings from tilings of $R(7,6t-6)$ where the leftmost five 
columns again match the pattern in Figure 11(f) in 160 ways (we extend the pattern in Figure 
6(a) which is present in the leftmost top or bottom four rows of $R(7,6t-6)$ in 20 ways using 
the generators in Figure 6(b) and (d); the remaining three rows being extended by simply padding 
with three $R(3,2)$ rectangles). Now consider tilings of $R(7,6t-6)$ where the leftmost five columns match 
the pattern in Figure 11(e). We can get a tiling of $R(7,6t)$ where the leftmost five columns 
match the pattern in Figure 11(f) in 128 ways from such a tiling (we extend the pattern in 
Figure 6(a) which is present in the leftmost top or bottom four rows of $R(7,6t-6)$ in 16 ways 
using the generators in Figure 6(c) and (e); the remaining three rows being extended by 
padding with three $R(3,2)$ rectangles). Finally consider tilings of $R(7,6t-6)$ where the leftmost four columns 
match the pattern in Figure 11(a). We can extend such tilings to tilings of $R(7,6t)$ where the 
leftmost five columns match the pattern in Figure 11(f) in 64 ways (we extend the pattern in 
Figure 5(a) which is present in the leftmost top or bottom four rows of $R(7,6t-6)$ in 8 ways using 
the generator in Figure 5(c); the remaining three rows being extended by padding with three $R(3,2)$ rectangles). 
So, we get the following equation:

\begin{eqnarray}
S(z) & = & 160zS(z) + 128zT(z) + 64zH(z)
\end{eqnarray}

Consider tilings of $R(7,6t)$ where the leftmost five columns match the pattern in Figure 
11(e). Following similar reasoning as above and observing the symmetry of Figures 11(e) and 11(f) and, 
the reader can easily verify that we get such tilings from tilings of $R(7,6t-6)$ where the 
leftmost five columns match the pattern in Figure 11(f) in 128 ways, 
and from those that match the pattern in Figure 11(e) itself in 160 ways. Consider tilings of $R(7,6t-6)$ where 
the leftmost four columns match the pattern in Figure 11(a). From such tilings we can get tilings 
of $R(7,6t)$ where the leftmost five columns match the pattern in Figure 11(e) in 128 ways (we extend the pattern in 
Figure 5(a) which is present in the leftmost top or bottom four rows in 16 ways using the generator 
in Figure 5(b); the remaining three rows being extended by padding with three $R(3,2)$ rectangles). So we have the 
following equation:

\begin{eqnarray}
T(z) & = & 128zS(z) + 160zT(z) + 128zH(z)
\end{eqnarray}

Solving equations (13), (14) and (15), we get 
\begin{eqnarray}
H(z) & = & \frac{16z(1 - 32z)(1 - 288z)}{1 - 448z + 9216z^2}  \\
S(z) & = & \frac{1024z^2(1 + 96z)}{1 - 448z + 9216z^2} \\
T(z) & = & \frac{2048z^2(1 - 96z)}{1 - 448z + 9216z^2}
\end{eqnarray}

Let us consider tilings of $R(7,6t)$ where the leftmost five columns match the pattern in Figure 11(b). 
Such tilings can only be obtained by retiling tilings of $R(7,6t)$ where the leftmost five columns match 
the pattern in Figure 11(f) (as shown in Figure 12(a)). 
Note that in such retilings, we disturb the three $R(3,2)$ rectangles in Figure 11(e) 
and construct back only two. So we have 

\begin{eqnarray}
J(z) & = & \frac{1}{2}S(z)  \nonumber \\
     & = & \frac{512z^2(1 + 96z)}{1 - 448z + 9216z^2}
\end{eqnarray}        

Finally, consider tilings of $R(7,6t)$ where the leftmost five columns match the pattern in Figure 11(c) 
or (d). Both these kinds of tilings of tilings can be obtained only by retiling tilings of $R(7,6t)$ where 
the leftmost five columns match the pattern in Figure 11(e) (as shown in Figures 12(d) and (e)). 
Note that in both these retilings we disturb 
the three $R(3,2)$ rectangles in Figure 11(e) and construct back only one. So we have 

\begin{eqnarray}
K(z)  =  P(z)  =  \frac{1}{4}T(z)  =  \frac{512z^2(1 - 96z)}{1 - 448z + 9216z^2}
\end{eqnarray}		    

Once we retile the leftmost five columns of $R(7,6t)$ to get patterns matching 
Figure 11(b) or (d), we apply Moore's rules for any further extension of tilings.
Note that in such an extension, the topmost or bottommost five rows are extended using
Moore's rules and the remaining two rows are extended simply by padding using $R(2,3)$ rectangles. 
Observe that once we do such extensions, we never generate patterns
in Figure 11. 
So, these
tilings obtained with Moore's rules are distinct from those obtained using 
patterns in Figure 11 (a), (e) and (f) using our generators in Figure 5 and 6, and
accounted as in the generating functions in equations (16)-(18). The reader must note that 
the tilings of $R(7,6t)$ where the leftmost five columns match the pattern in  
Figure 11(b) will be distinct from those generated by Moore's approach and having the 
tiling of the leftmost five columns similar to the pattern in Figure 11(b). This is 
because in the extension by Moore's scheme, we extend only the topmost (bottommost) five 
rows, padding the remaining two rows by $R(2,3)$ rectangles. Thus, those tilings generated by Moore's 
scheme will have a $R(2,3)$ in the bottom (top) two rows from columns 4 to 6, while those tilings 
generated by our generators in Figure 5 and 6 do not have this $R(2,3)$ after retiling patterns 
in Figure 11(e) and (f) to get the pattern in Figure 11(b). 
In Moore's extension, we extend $R(7,6)$ first by $6r$ columns (say), where $1\leq r\leq t$, 
using our generators in Figures 5 and 6, we then retile the tiling of the leftmost five columns 
matching patterns in Figures 11(f) or (e) to get the patterns in Figures 11(b) and 
(d), and then extend the patterns in Figure 7(c) and (d), which are present (in the 
leftmost four columns) in the topmost (bottommost) five rows, by $6(t-r)$ columns by Moore's 
scheme (as shown in Figure 13(b) and (a)), padding the remaining tow rows 
by $R(2,3)$ rectangles. Note that any two tilings of $R(7,6t)$, one generated by extension of $R(7,6r)$ 
where the leftmost five columns match the pattern in Figure 11(b), and the other generated by 
extension of $R(7,6r)$ where the leftmost five columns match the pattern in Figure 11(d), will 
be distinct. This is because in both these $R(7,6r)$ rectangles, the squares $(3,4)$ and $(3,5)$ 
are tiled in different ways, and the reader can easily verify that Moore's extension 
does not disturb this particular arrangement. 
Moore had derived generating functions $G_1'(z)$ and $G_3'(z)$, where the coefficient of $z^t$ 
gives all faultfree tilings of the 5-row patterns shown in Figure 13, $G_1'(z)$ being 
for the pattern in Figure 13(a) and $G_3'(z)$ being for that in Figure 13(b). 
These generating functions are :

\begin{eqnarray}
G_1'(z) & = & \frac{4z(1+z)(1-z-10z^2)}{1-2z-31z^2-40z^3-20z^4} \\
G_3'(z) & = & \frac{2z(1+14z+5z^2)}{1-2z-31z^2-40z^3-20z^4}
\end{eqnarray} 

\begin{figure}[htbp]
\centerline{
{\scalebox{.25}{\includegraphics{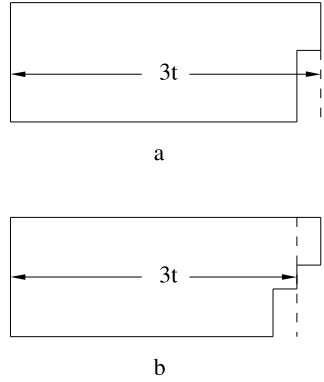}}}
}
\caption{Patterns corresponding to Moore's generating functions.}
\end{figure}

So, if we count tilings of $R(7,6t)$ using our scheme described above, we will essentially be 
convoluting our generating functions $J(z)$ and $P(z)$ 
with Moore's functions $G_1'(z)$ and $G_3'(z)$. Since $J(z)$ and $P(z)$ consider all tilings of $R(7,6t)$ 
whereas Moore's functions consider tilings of $R(5,3t)$, we change $z$ to $z^2$ in both 
$J(z)$ and $P(z)$, to make them compatible for convolution with Moore's generating functions. 
Note that the coefficients of all odd powers of z in $J(z)$ and $P(z)$ thus obtained will be zero. 
We first consider the generating function $J(z)$. For the sake of clarity and insight, we write the 
generating functions $J(z)$, $G_1'(z)$ and $G_3'(z)$ in their expanded form, as shown below : 

\begin{eqnarray}
J(z)    & = & a_2z^2 + a_4z^4 + a_6z^6 + ....... \\
G_1'(z) & = & b_1z + b_2z^2 + b_3z^3 + ..........\\ 
G_3'(z) & = & c_1z + c_2z^2 + c_3z^3 + ..........
\end{eqnarray}

Consider the first case in our scheme by Moore's extension mentioned above, in which we first obtain 
the pattern in Figure 11(b) and then extend using Moore's approach. If we choose to extend the retiled 
$R(7,6)$ by Moore's scheme, then we can pad the remaining two rows in $2^{2t-2}$ ways (since every $R(2,3)$ 
can be tiled in two ways). If we choose to extend the topmost five rows, we can do so in $\frac{b_{2t-1}}{4}$ ways, and 
if we choose to extend the bottommost five rows, we can do so in $\frac{c_{2t-1}}{2}$ ways. So, the total number of ways 
of extending the retiled $R(7,6)$ by Moore's scheme is $2^{2t-2}a_2(\frac{b_{2t-1}}{4}+\frac{c_{2t-1}}{2})$. 
Suppose we choose to 
extend the retiled $R(7,12)$ by Moore's scheme, then we can pad the remaining two rows in $2^{2t-4}$ ways. The 
topmost five rows can be extended in $\frac{b_{2t-3}}{4}$ ways, and the bottommost five rows can be extended 
in $\frac{c_{2t-3}}{2}$ ways. So the total number of ways of extending the retiled $R(7,12)$ is $2^{2t-4}a_4(\frac{b_{2t-3}}{4}+\frac{c_{2t-3}}{2})$. Going 
on like this, the total number of ways of tiling $R(7,6t)$ by Moore's scheme, when the tiling of the leftmost five 
columns of $R(7,6r)$, where $1\leq r \leq t$, matches the pattern in Figure 11(b) is :

\begin{eqnarray}
\alpha _{2t} & = & 2^{2t-2}a_2(\frac{b_{2t-1}}{4}+\frac{c_{2t-1}}{2}) + 2^{2t-4}a_4(\frac{b_{2t-3}}{4}+\frac{c_{2t-3}}{2}) + ..... \nonumber \\
             & = & \sum_{i=1}^{2t-1}(\frac{2^{2t-i}a_ib_{2t-i+1}}{4} + \frac{2^{2t-i}a_ic_{2t-i+1}}{2}) \nonumber \\
	     & = &  2^{2t}\{\frac{1}{4}\sum_{i=1}^{2t-1}\frac{a_i}{2^i}b_{2t-i+1} + \frac{1}{2}\sum_{i=1}^{2t-1}\frac{a_i}{2^i}c_{2t-i+1}\}   
\end{eqnarray}

Let L(z) be the generating function in which the coefficient of $z^t$ is $\frac{\alpha_t}{2^t}$. Since our extension 
scheme always extends by six columns, so $\alpha_t$ will be zero when $t$ is odd. From equation 23, we get,

\begin{eqnarray}
L(z) & = & J(\frac{z^2}{4})\{\frac{G_1'(z)+2G_3'(z)}{4z}\}
\end{eqnarray}

The generating function where the coefficient of $z^t$ is $\alpha_t$ is $L(2z)$. 
Let $\beta_t$ be the total number of ways of tiling $R(7,6t)$ by Moore's scheme when the tiling of 
the leftmost five columns of $R(7,6r)$, where $1\leq r \leq t$, matches the pattern in Figure 11(d). 
Let $M(z)$ be the generating function where coeffiecient of $z^t$ is $\frac{\beta_t}{2^t}$. Proceeding exactly 
in the same manner as above, we get,

\begin{eqnarray}
M(z) & = & P(\frac{z^2}{4})\frac{G_3'(z)}{z}
\end{eqnarray}

So, the generating function, where coefficient of $z^t$ is $\beta_t$ is $M(2z)$. The total number of 
ways of faultfreely tiling $R(7,6t)$ by the above approach, when extension of $R(7,6)$ is done only 
on one side (either entirely to the left or entirely to the right) will be twice the sum of all 
the generating functions derived above, since we have two choices for the direction of 
extension (either left or right). So, the generating function $Q_1(z)$ where the coefficient 
of $z^{2t}$, where $t\geq2$, is the total number of ways of tiling 
$R(7,6t)$ faultfreely, when $R(7,6)$ (as in Figure 10(a)) is extended only on one side (either 
entirely to the left or entirely to the right) is,

\begin{eqnarray}
Q_1(z) & = & 2\{L(2z)+M(2z)+H(z^2)+K(z^2)+S(z^2)+T(z^2)\}
\end{eqnarray}    

Now let us consider the pattern in Figure 10(b). In this case, Lemma 9 shows that there are 8 ways of 
tiling $R(7,6)$ faultfreely. The reader should observe that in the leftmost (rightmost) four columns 
of $R(7,6)$ (as shown in Figure 10(b)), the pattern in Figure 2(c) is present in the bottommost (topmost) five 
rows while the pattern in Figure 2(d) is present in the topmost (bottommost) five rows. 
Consider the pattern in Figure 2(c), we have the following lemma. 

\begin{figure}[htbp]
\centerline{
{\scalebox{.23}{\includegraphics{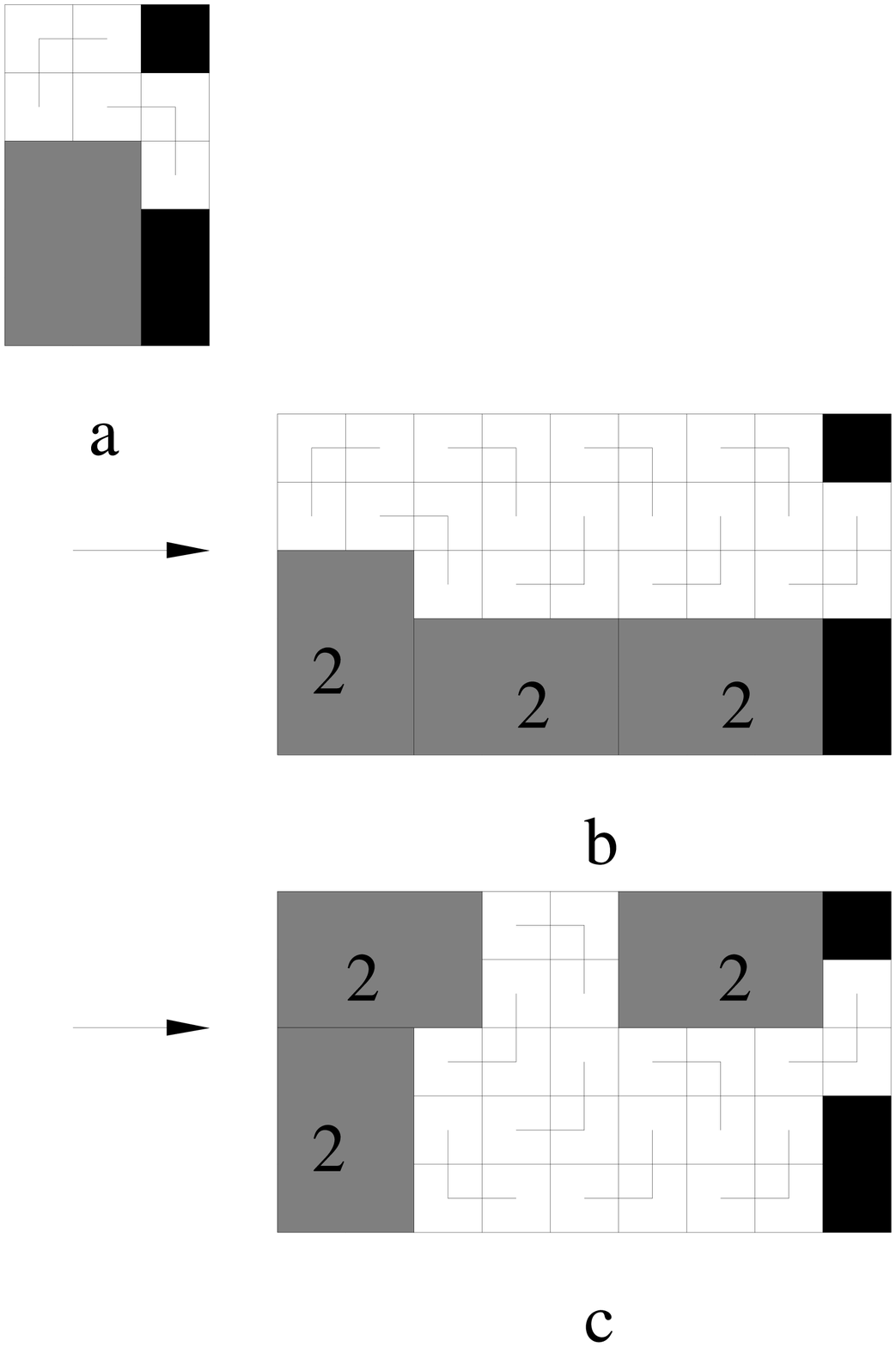}}}
{\scalebox{.23}{\includegraphics{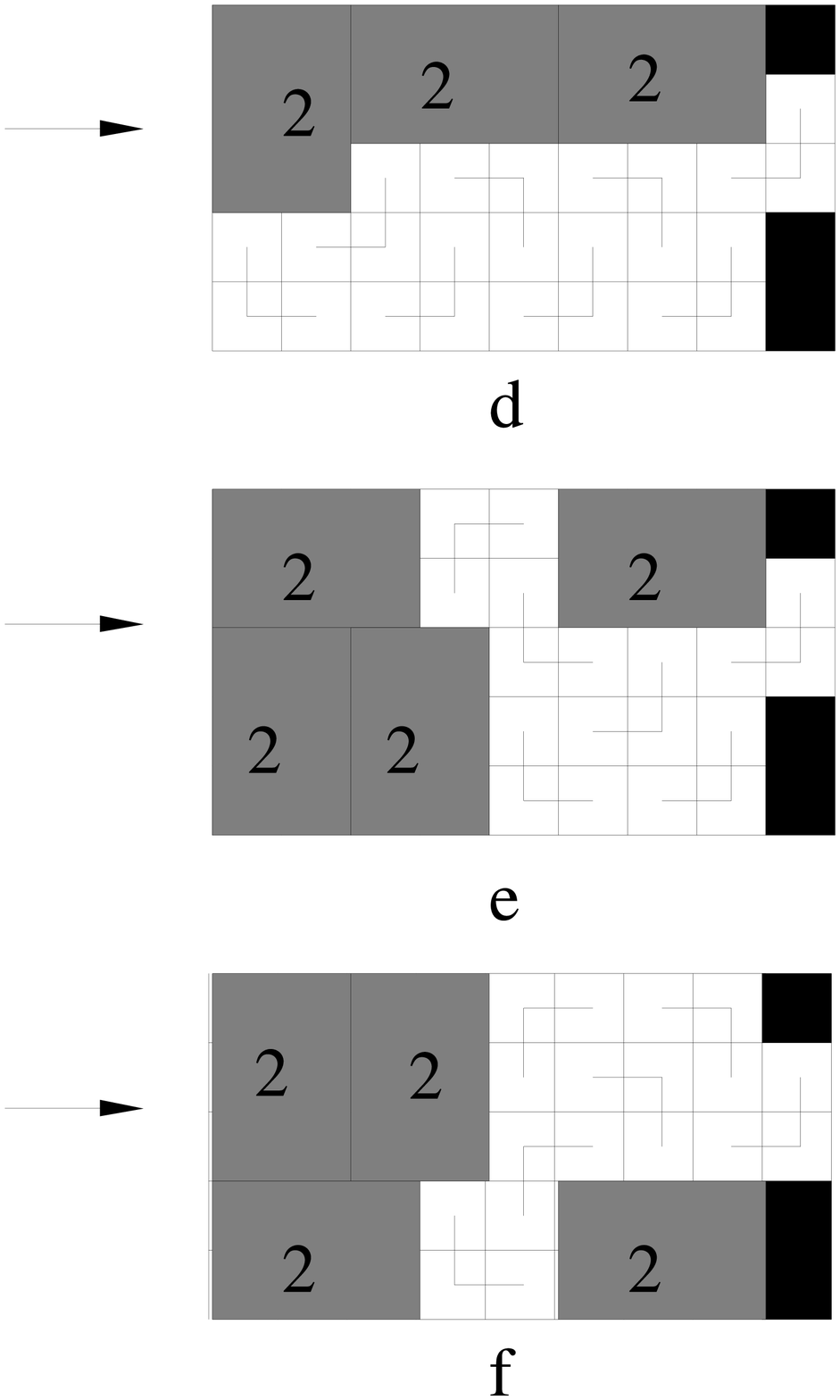}}}
{\scalebox{.23}{\includegraphics{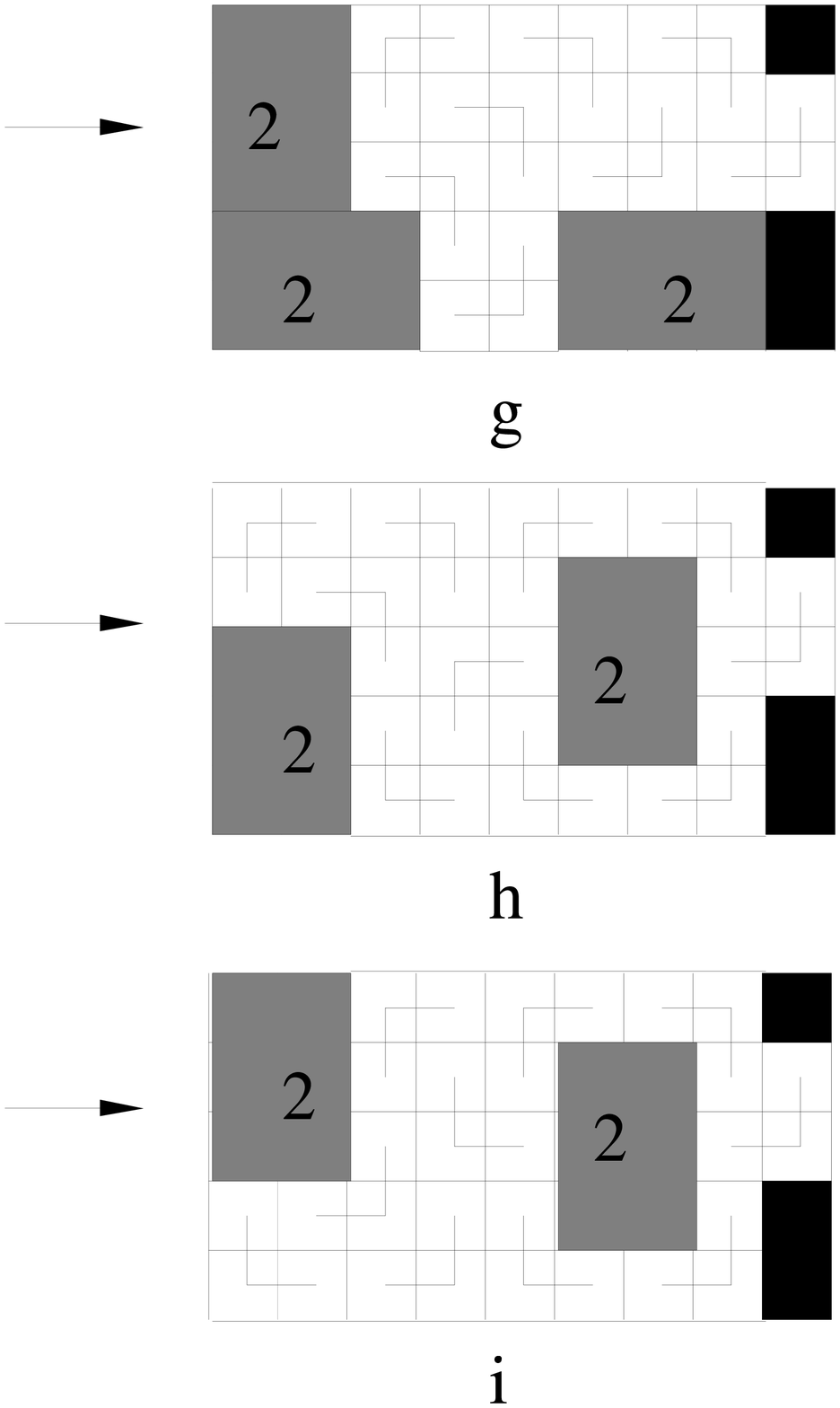}}}
}
\caption{Number of ways of expanding the pattern in Figure 3(e).}
\end{figure}

\begin{figure}[htbp]
\centerline{
{\scalebox{.23}{\includegraphics{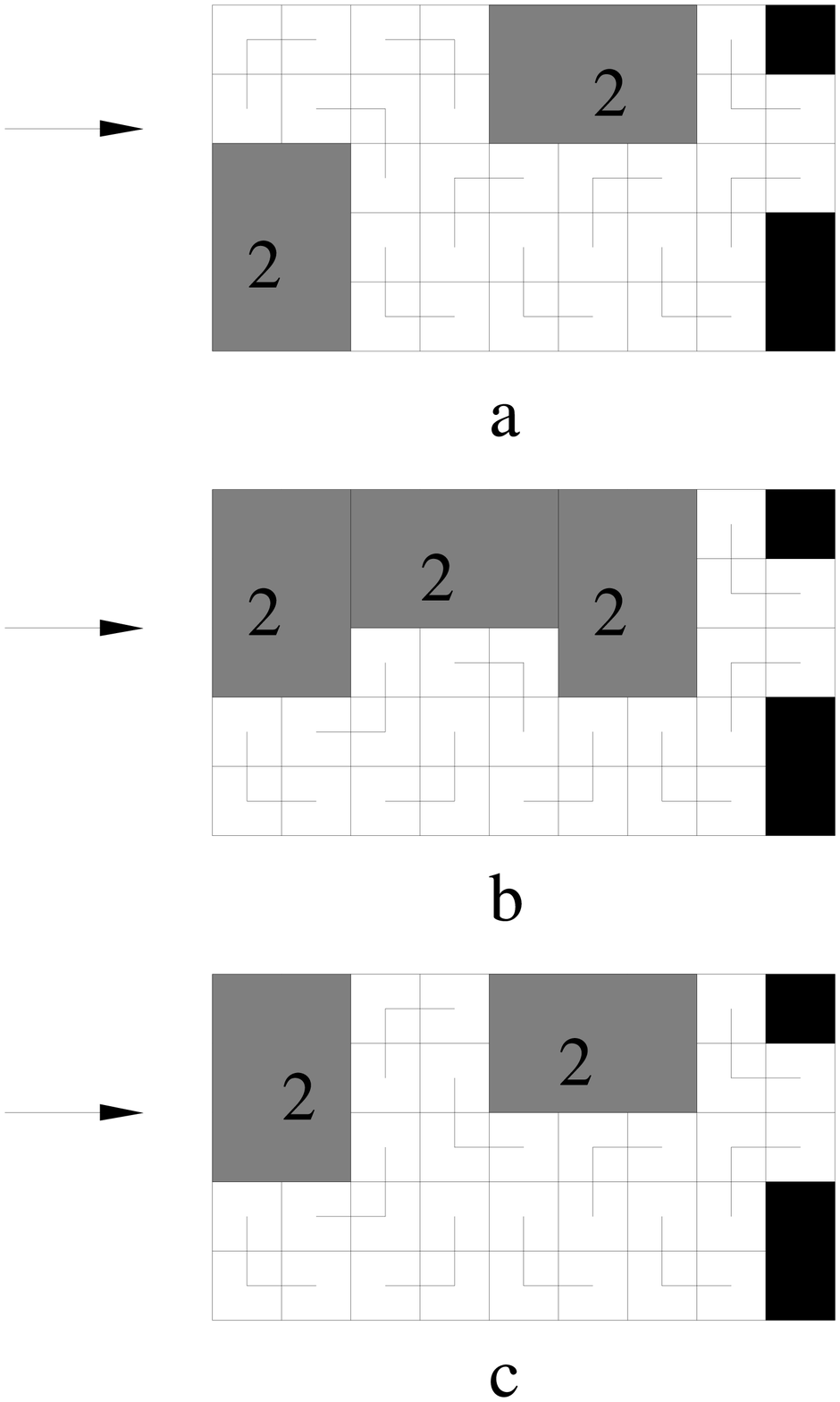}}}
{\scalebox{.23}{\includegraphics{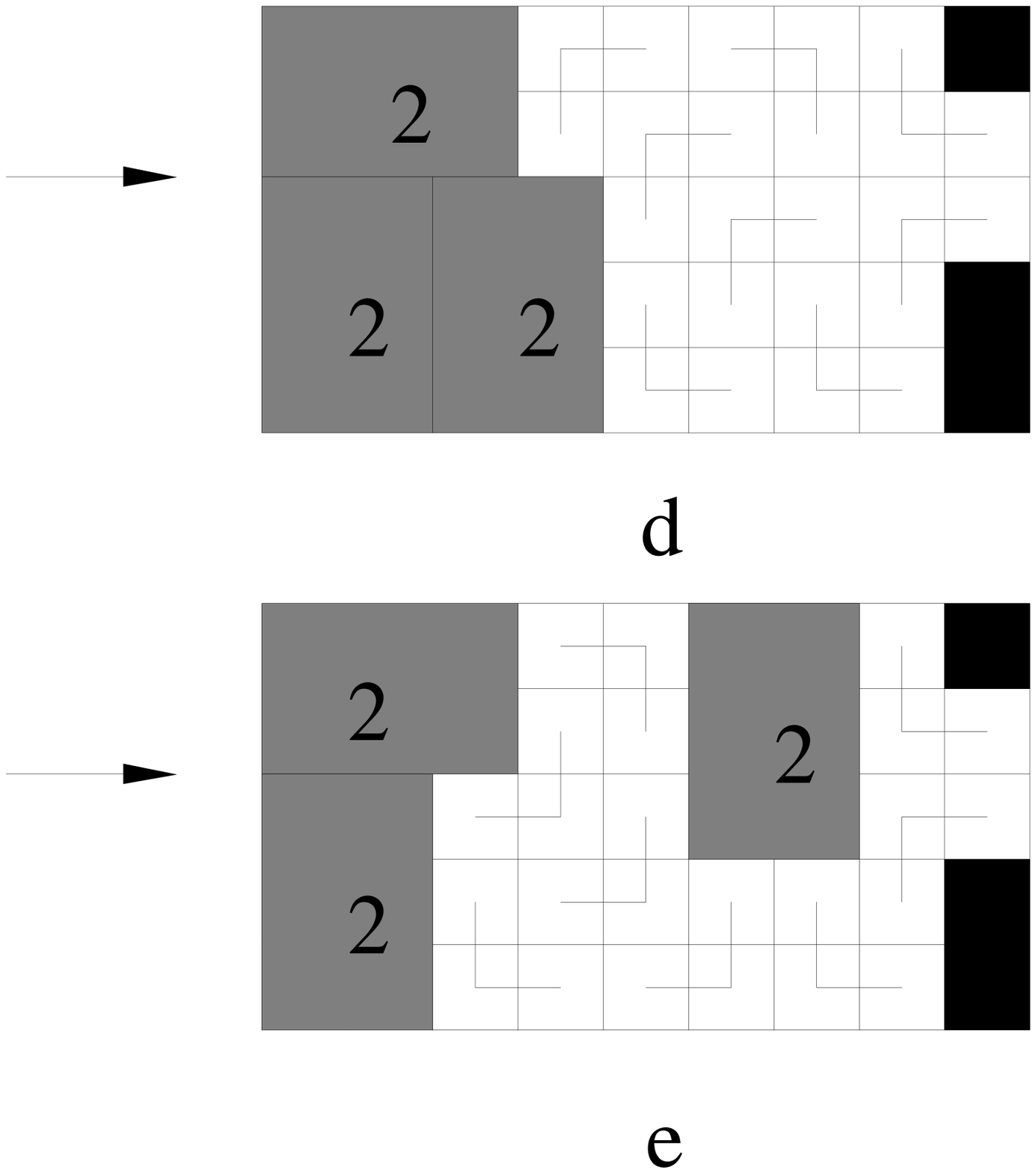}}}
}
\caption{Number of ways of expanding the pattern in Figure 3(e).}
\end{figure}

\begin{lemma}
The only possible ways of extending the pattern shown in Figure 2(c), 
so that the resulting tiling is also faultfree, are those depicted in Figures 14 and 15.
\end{lemma}
\begin{proof}
As argued in Lemma 6, any tiling of the leftmost three columns of $R(5,3t)$ 
must match this same pattern or its flipped counterpart, or the pattern 
shown in Figure 3(c) or its flipped counterpart. Note that the same tromino cannot cover 
$(2,8)$, $(2,9)$ and $(3,9)$ as then a tromino will have to cover $(1,8)$, $(1,7)$ and 
$(2,7)$. Now the reader can easily see that the $7th$ vertical grid line will become a 
fault line. 
   
Consider the case when a tiling of the leftmost three columns of the generator extending 
Figure 3(e) matches this same pattern. The tromino covering $(4,3)$ and $(5,3)$ has two 
permissible orientations. If it covers $(4,4)$ then $R(2,3)$ is completed by the tromino 
covering $(5,4)$. If it covers $(5,4)$ and the tromino covering $(4,4)$ also covers $(4,5)$ 
and $(5,5)$ then we again get $R(2,3)$. The remaining can be tiled now only in the ways 
shown in Figure 14(b). Suppose the tromino covering $(4,3)$ and $(5,3)$ covers $(5,4)$. 
If the tromino covering $(4,4)$ covers $(3,4)$ and $(4,5)$ then it is not hard to see that 
the $7th$ vertical grid line will become a fault line. So, the only other orientation 
permissible for the tromino covering $(4,4)$ is that it covers $(3,4)$ and $(3,5)$. If a 
tromino covers $(2,9)$, $(3,9)$ and $(3,8)$, then from symmetry, we find that the only 
possible ways to tile the generator are those shown in Figure 14(h). The only case 
remaining is that a tromino covers $(1,8)$, $(2,8)$ and $(2,9)$ and another tromino covers 
$(4,8)$, $(3,8)$ and $(3,9)$. Again, it is easy to verify that the only possible ways of 
tiling the remaining area are those shown in Figure 15(a). 

Now consider the case when a tiling of the leftmost three columns of a generator matches 
the flipped counterpart of the pattern shown in Figure 3(e). The tromino covering $(1,3)$ 
and $(2,3)$ has two permissible orientations. If it covers $(2,4)$ then $R(2,3)$ is 
completed by the tromino covering $(1,4)$. If it covers $(1,4)$ and the tromino covering 
$(2,4)$ covers $(1,5)$ and $(2,5)$ then we again get $R(2,3)$. Now if a tromino covers 
$(2,9)$, $(3,9)$ and $(3,8)$ then the only possible ways of tiling the remaining area are 
those shown in Figure 14(d). If instead, a tromino covers $(1,8)$, $(2,8)$ and $(2,9)$ and 
another tromino covers $(4,8)$, $(3,8)$ and $(3,9)$ then the only ways of tiling the 
remaining area are those shown in Figure 15(b). Suppose the tromino covering $(1,3)$ 
and $(2,3)$ covers $(1,4)$. If the tromino covering $(2,4)$ covers $(3,4)$ and $(2,5)$ then 
it can be easily verified that the $7th$ vertical grid line will become a fault line. So, 
apart from the cases discussed above, the only other permissible orientation for the tromino 
covering $(2,4)$ is that it covers $(3,4)$ and $(3,5)$. Now if a tromino covers $(2,9)$, 
$(3,9)$ and $(3,8)$, then from symmetry, we can see that the remaining can be tiled only in 
the ways shown in Figure 14(i). On the other hand, if a tromino covers $(1,8)$, $(2,8)$ and 
$(2,9)$ and another tromino covers $(4,8)$, $(3,8)$ and $(3,9)$ then the reader can verify 
that the only possible ways of tiling the remaining area are those shown in Figure 15(c). 

Considering the third case when a tiling of the leftmost three columns matches the pattern 
shown in Figure 3(c). The tromino covering $(3,3)$ has three permissible orientations. If 
it covers $(3,4)$ and $(2,4)$ then $R(3,2)$ is completed by the tromino covering $(2,3)$. 
If it covers $(2,3)$ and $(3,4)$ then $R(3,2)$ is completed by the tromino covering $(1,3)$. 
Suppose we have $R(3,2)$ as per the above two cases. The tromino covering $(4,4)$ and $(5,4)$ 
has two possible orientations, it may cover $(4,5)$ or $(5,5)$. If a tromino covers $(1,8)$, $(2,8)$ 
and $(2,9)$ and another tromino covers $(3,9)$, $(3,8)$ and $(4,8)$, then in both the above 
orientations of the tromino covering $(4,4)$ and $(5,4)$ we will get $R(2,3)$ at the bottom. 
The reader can now see that this case requires us to tile $R(3,3)$ which we know is impossible 
to tile \cite{chu}. So we conclude that a tromino must cover $(2,9)$, $(3,9)$ and $(3,8)$. 
If the tromino covering $(4,4)$ and $(5,4)$ covers $(4,5)$, then $R(2,3)$ is completed by the 
tromino covering $(5,5)$. But then a tromino will have to cover $(5,7)$, $(5,8)$ and $(4,8)$. The 
reader can easily see that the square $(1,5)$ will become inaccessible in this case. We conclude 
that the tromino covering $(4,4)$ and $(5,4)$ must cover $(5,5)$. There are two possible 
orientations of the tromino covering $(4,8)$ and $(5,8)$. If it covers $(4,7)$ then $R(2,3)$ is 
completed by the tromino covering $(5,7)$ and if it covers $(5,7)$, then from the above discussion, 
it follows that $R(2,3)$ is again completed by the tromino covering $(5,6)$. The remaining area 
is uniquely tileable as shown in Figure 14(f). Now consider the case when the tromino covering 
$(3,3)$ covers $(3,4)$ and $(4,4)$. If a tromino covers $(1,8)$, $(2,8)$ and $(2,9)$ and another 
tromino covers $(4,8)$, $(3,8)$ and $(3,9)$, then the reader can verify that $(5,6)$ will become 
inaccessible. So a tromino must cover $(2,9)$, $(3,9)$ and $(3,8)$. There are two orientations of 
the tromino covering $(4,8)$ and $(5,8)$. If it covers $(4,7)$ then $R(2,3)$ is completed by the 
tromino covering $(5,7)$ and if it covers $(5,7)$ then $R(2,3)$ is completed by the tromino 
covering $(5,6)$. As can be verified, the remaining area has a unique tiling as shown in 
Figure 14(g). 
 
Finally, we consider the case when a tiling of the leftmost three columns of a generator matches 
the flipped counterpart of the pattern shown in Figure 3(c). The tromino covering $(3,3)$ has 
three permissible orientations. If it covers $(3,4)$ and $(4,4)$ then $R(3,2)$ is completed by 
the tromino covering $(4,3)$. If it covers $(3,4)$ and $(4,3)$ then $R(3,2)$ is again completed 
by the tromino covering $(5,3)$. Suppose we have $R(3,2)$ as per the above cases. Further suppose 
that a tromino covers $(1,8)$, $(2,8)$ and $(2,9)$ and another tromino covers $(4,8)$, $(3,8)$ 
and $(3,9)$. If the tromino covering $(5,5)$ and $(5,6)$ covers $(4,6)$ then $R(3,2)$ is completed 
by the tromino covering $(4,5)$. It is not hard to see that the tromino covering $(3,7)$ will make 
$(1,7)$ inaccessible. So the tromino covering $(5,5)$ and $(5,6)$ must cover $(4,5)$. For the same  
reason, the tromino covering $(4,6)$ and $(3,6)$ must cover $(3,7)$. The remaining area has a 
unique tiling as shown in Figure 15(d). Now suppose that a tromino covers $(2,9)$, $(3,9)$ and 
$(3,8)$. If the tromino covering $(1,4)$ and $(2,4)$ covers $(2,5)$, then $R(2,3)$ is completed 
by the tromino covering $(1,5)$. The reader can verify that the chain of trominoes covering $(1,8)$, 
$(2,7)$ and $(3,5)$ will make $(5,5)$ inaccessible. So the tromino covering $(1,4)$ and $(2,4)$ 
must cover $(1,5)$, and for the same reason, the tromnino covering $(2,5)$ will not complete $R(2,3)$ 
with this tromino. If the tromino covering $(2,5)$ also covers $(3,5)$ and $(2,6)$, then it can be 
easily seen that the squares $(1,8)$ and $(2,8)$ will become inaccessible. So we conclude that the 
tromino covering $(2,5)$ must cover $(3,5)$ and $(3,6)$. If the tromino covering $(1,6)$ and $(2,6)$ 
covers $(1,7)$ then $R(2,3)$ is completed by the tromino covering $(1,8)$ and if it covers $(2,7)$, 
then $R(2,3)$ is completed by the tromino covering $(1,7)$. The remaining area is uniquely tileable 
as shown in Figure 14(e). Now suppose that the tromino covering $(3,3)$ covers $(3,4)$ and $(2,4)$. 
Also, suppose that a tromino covers $(1,8)$, $(2,8)$ and $(2,9)$ and another tromino covers $(4,8)$, 
$(3,8)$ and $(3,9)$. The tromino covering $(1,6)$ and $(1,7)$ has two permissible orientations. If 
it covers $(2,6)$ then $R(3,2)$ is completed by the tromino covering $(2,7)$. If it covers $(2,7)$, 
then $R(3,2)$ is again completed by the tromino covering $(2,6)$ and $(3,7)$. In case the tromino 
covering $(2,6)$ covers $(3,6)$ and $(3,5)$ then the reader can easily see that $(3,7)$ will become 
inaccessible. So this case is not possible. The remaining area has a unique tiling as shown in 
Figure 15(e). Lastly, suppose that a tromino covers $(2,9)$, $(3,9)$ and $(3,8)$. If the tromino 
covering $(1,8)$ and $(2,8)$ covers $(1,7)$, then $R(2,3)$ is completed by the tromino covering 
$(1,6)$, and if it covers $(2,7)$, then $R(2,3)$ is completed by the tromino covering $(1,7)$. 
Consider the tromino covering $(4,8)$ and $(5,8)$. If it covers $(4,7)$ then $R(2,3)$ is completed 
by the tromino covering $(5,7)$. Clearly, the square $(3,7)$ becomes inaccessible in this case. So 
the tromino covering $(4,8)$ and $(5,8)$ must cover $(5,7)$. For this same reason, the tromino 
covering $(4,7)$ will not complete $R(2,3)$ with this tromino. If the tromino covering $(4,7)$ 
covers $(3,7)$ and $(4,6)$, then the tromino covering $(5,6)$ will make $(3,6)$ inaccessible. So 
the tromino covering $(4,7)$ must cover $(3,7)$ and $(3,6)$. The reader can verify that the remaining 
area has a unique tiling as shown in Figure 14(c). Thus, Figures 14 and 15 show all possible ways of 
expanding the pattern shown in Figure 3(e). \hfill \qed  
\end{proof}

Now consider the pattern in Figure 2(d). We have the following lemma. 

\begin{figure}[htbp]
\centerline{
{\scalebox{.23}{\includegraphics{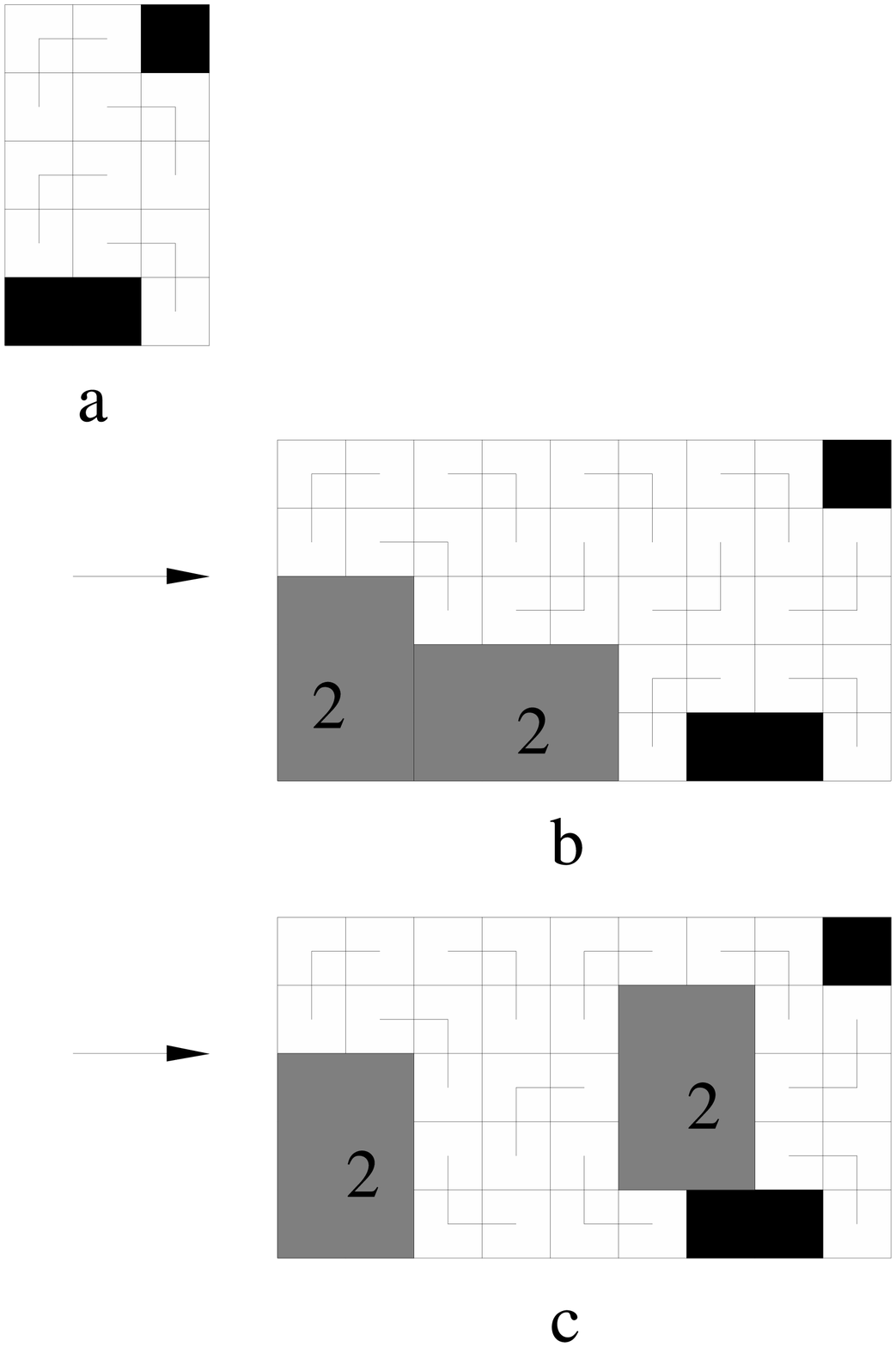}}}
{\scalebox{.23}{\includegraphics{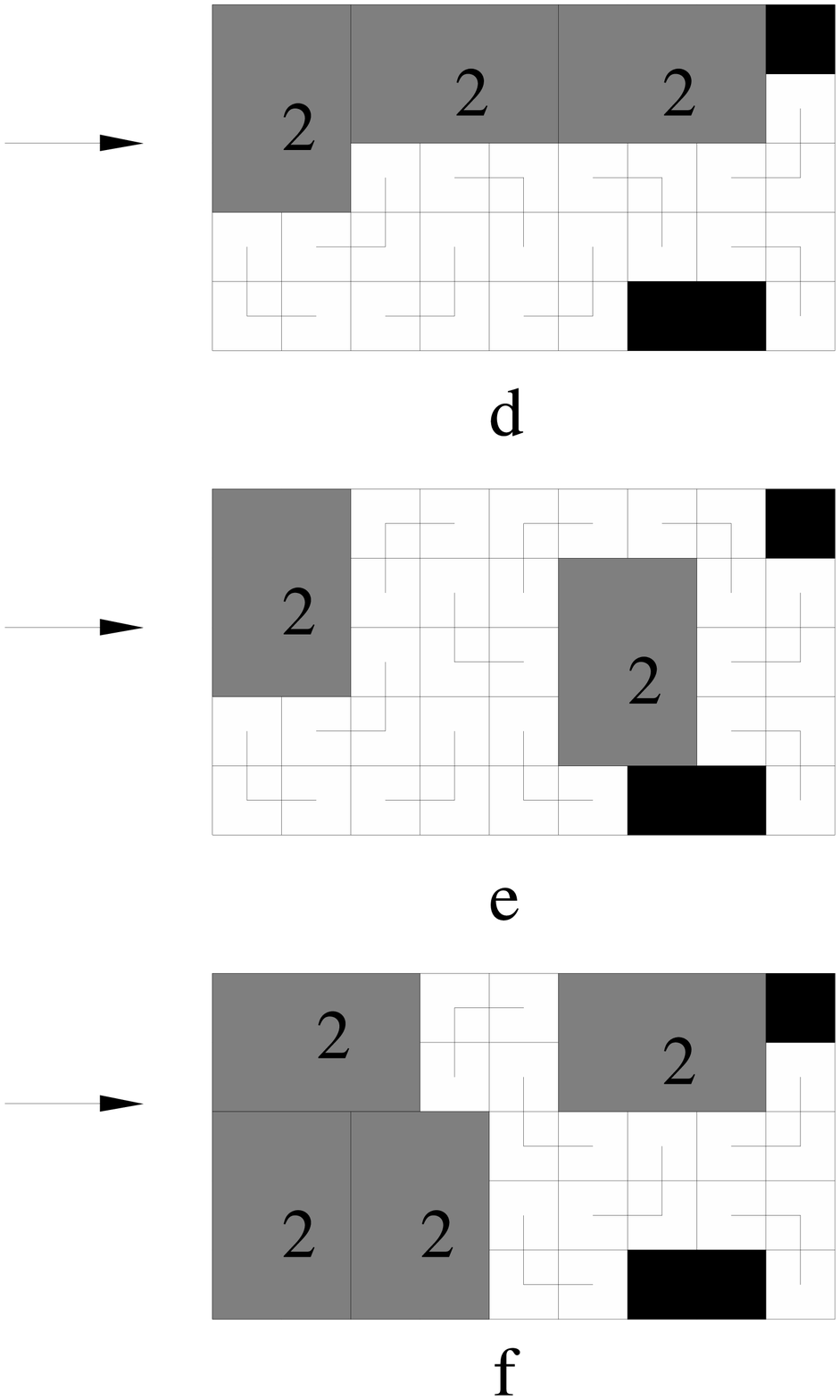}}}
{\scalebox{.23}{\includegraphics{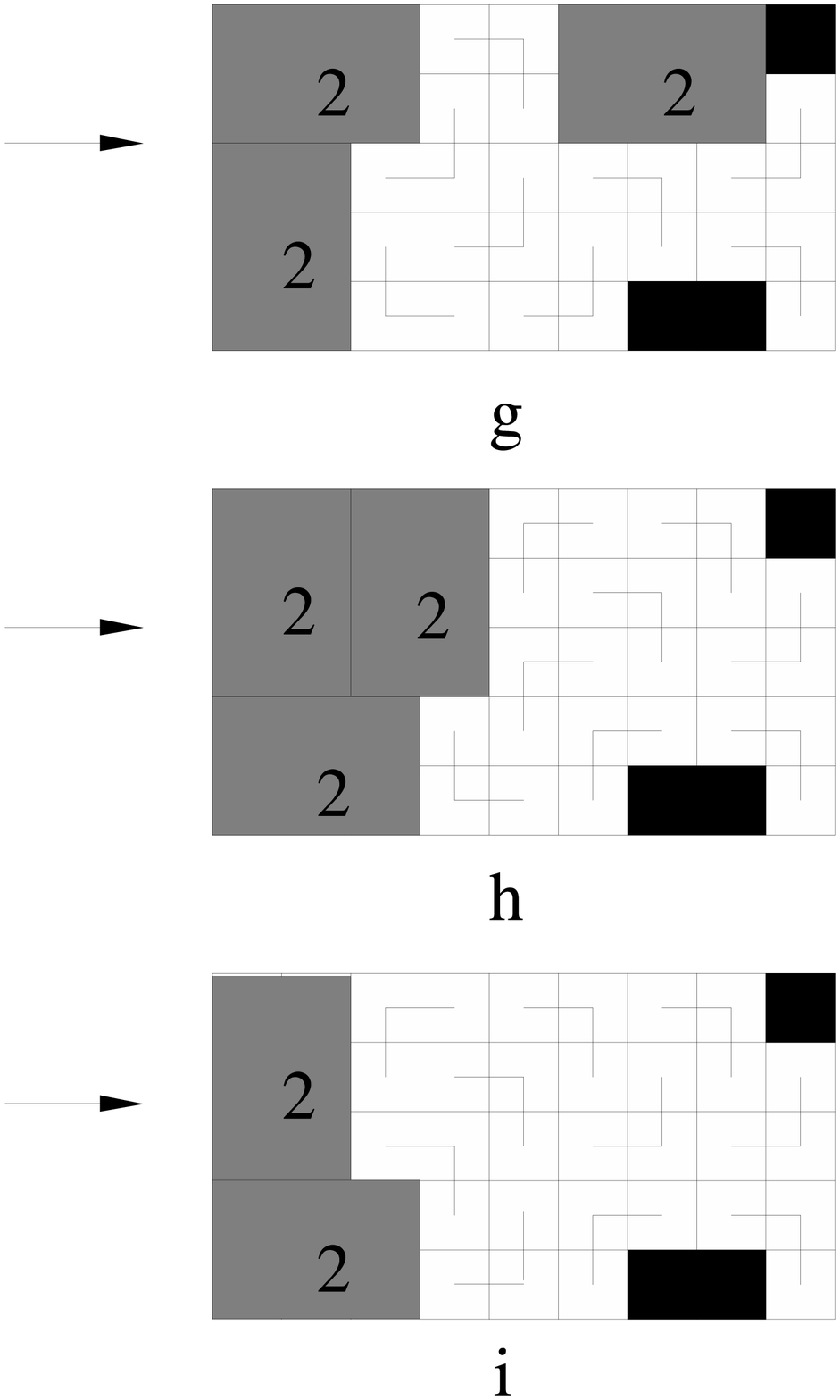}}}
}
\caption{Number of ways of expanding the pattern in Figure 2(d).}
\end{figure}

\begin{lemma}
The only possible ways of extending the pattern shown in Figure 2(d), 
so that the resulting tiling is also faultfree, are those depicted in Figure 16.
\end{lemma}
\begin{proof}
From Lemma 6, we know that any tiling of the leftmost three columns of $R(5,3t)$ must match 
the pattern in Figures 3(c) or 3(e), or their symmetric counterparts. In all these cases the 
reader can easily see that the tromino covering $(5,9)$ will also cover $(4,9)$ and $(4,8)$. 
   
First consider the case when any tiling of the leftmost three columns of a generator extending 
the pattern in Figure 2(d) matches the pattern in Figure 3(e). It can be easily seen that the 
tromino covering $(1,3)$ will also cover $(1,4)$ and $(2,4)$. The tromino covering $(4,3)$ and 
$(5,3)$ has two permissible orientations, it can either cover $(4,4)$ or $(5,4)$. If it covers 
$(4,4)$ then $R(2,3)$ is completed by the tromino covering $(5,4)$. If it covers $(5,4)$ then 
the tromino covering $(4,4)$ has three possible orientations. Consider the case when it covers 
$(4,5)$ and $(5,5)$ thereby completing a $R(2,3)$. Now it is evident that the tromino covering 
$(5,6)$ must cover $(4,6)$ and $(4,7)$. It can be verified easily that the rest of the area can 
only be tiled uniquely as shown in Figure 16(b). Now consider the second possible orientation of 
the tromino covering $(4,4)$, suppose that it covers $(3,4)$ and $(3,5)$. Note that the tromino 
covering $(4,5)$ and $(5,5)$ must cover $(5,6)$. If this is not the case, then it must cover $(4,6)$ 
and in this case $(5,6)$ will become inaccessible. If the tromino covering $(1,5)$ and $(2,5)$ 
covers $(2,6)$ then a $R(2,3)$ is completed by the tromino covering $(1,6)$. Also, the tromino 
covering $(1,8)$ must cover $(2,8)$ and $(2,9)$. The reader can now easily see that $(3,9)$ will 
become inaccessible in this case. So the tromino covering $(1,5)$ and $(2,5)$ must cover $(1,6)$. 
If the tromino covering $(4,6)$ and $(4,7)$ covers $(3,7)$ then a $R(3,2)$ is completed by the 
tromino covering $(3,6)$. If it covers $(3,6)$ then the tromino covering $(2,6)$ and $(2,7)$ must 
cover $(3,7)$. If this is not the case, then it must cover $(1,7)$. But in this case the tromino 
covering $(1,8)$ must cover $(2,8)$ and $(2,9)$, making the remaining area untileable. So a $(3,2)$ 
must be present in columns 6 and 7 from rows 2-4. The remaining area can now only be tiled as 
shown in Figure 16(c). Consider the third possible orientation of the tromino covering $(4,4)$, suppose 
that it covers $(3,4)$ and $(4,5)$. In this case, is can be easily seen that the tromino covering $(5,5)$ 
must cover $(5,6)$ and $(4,6)$. Now consider the two possible orientations of the tromino covering 
$(1,5)$. If it covers $(1,6)$ and $(2,5)$ then a $R(3,2)$ is completed by the tromino covering $(3,5)$, 
and if it covers $(1,6)$ and $(2,6)$ then a $(3,2)$ is completed by the tromino covering $(2,5)$. In 
either case the $7th$ vertical grid line becomes a fault line.  

Now consider the second case when a tiling of the leftmost three columns of a generator extending the 
pattern in Figure 2(d) matches the flipped counterpart of the pattern in Figure 3(e). Note that the 
tromino covering $(5,3)$ must cover $(5,4)$ and $(4,4)$. Consider the two possible orientations of the 
tromino covering $(1,3)$ and $(2,3)$. If it covers $(2,4)$ then a $R(2,3)$ is completed by the tromino 
covering $(1,4)$. If it covers $(1,4)$ then there are three possible orientations for the tromino covering 
$(2,4)$. Let us consider the first case and suppose that it covers $(2,5)$ and $(1,5)$ thereby completing 
a $R(2,3)$. Note that in this case a tromino must cover $(3,4)$, $(3,5)$, $(4,5)$ and another tromino 
must cover $(5,5)$, $(5,6)$ and $(4,6)$. Now the tromino covering $(4,7)$ and $(3,7)$ must cover $(3,6)$. 
If this is not the case then it must cover $(3,8)$, in which case the tromino covering $(3,6)$ will make 
$(1,6)$ inaccessible. It can be easily seen that the remaining area can only be tiled in two possible ways, 
as shown in Figure 16(d). Consider the second possible orientation of the tromino covering $(2,4)$, suppose 
that it covers $(3,4)$ and $(3,5)$. Note that the tromino covering $(4,5)$ and $(5,5)$ must cover $(5,6)$, 
in the other case when it covers $(4,6)$, $(5,6)$ will become inaccessible. If the tromino covering $(1,5)$ 
and $(2,5)$ covers $(2,6)$ then a $R(2,3)$ is completed by the tromino covering $(1,6)$. In this case the 
tromino covering $(1,8)$ must cover $(2,8)$ and $(2,9)$ making the square $(3,9)$ inaccessible. So the tromino 
covering $(1,5)$ and $(2,5)$ must cover $(1,6)$. Now if the tromino covering $(4,6)$ and $(4,7)$ covers $(3,7)$ 
then $R(3,2)$ is completed by the tromino covering $(3,6)$. If it covers $(3,6)$ then the tromino covering 
$(2,6)$ and $(2,7)$ must cover $(3,7)$. If this is not the case then it must cover $(1,7)$. The reader can now 
easily see that the tromino covering $(3,7)$, $(3,8)$ and $(2,8)$ will make $(1,8)$ inaccessible. The remaining 
area can now only be tiled as shown in Figure 16(e). Now consider the third possible orientation of the tromino 
covering $(2,4)$, suppose that it covers $(3,4)$ and $(2,5)$. In this case, the tromino covering $(1,5)$ must 
cover $(1,6)$ and $(2,6)$. Consider the two possible orientations of the tromino covering $(5,5)$ and $(5,6)$. 
If it covers $(4,5)$ then a $R(3,2)$ is completed by the tromino covering $(3,5)$ and if it covers $(4,6)$ then 
a $R(3,2)$ is completed by the tromino covering $(4,5)$. In both these cases the $7th$ vertical grid line becomes 
a fault line. 

Consider the third case when a tiling of the leftmost three columns of a generator extending the pattern in 
Figure 2(d) matches the pattern in Figure 3(c). The tromino covering $(3,3)$ has three permissible orientations. 
If it covers $(4,3)$ and $(3,4)$ then a $R(3,2)$ is completed by the tromino covering $(5,3)$. If it covers 
$(3,4)$ and $(4,4)$ then a $R(3,2)$ is completed by the tromino covering $(4,3)$. We claim that the tromino 
covering $(1,4)$ and $(2,4)$ must cover $(1,5)$. If this is not the case then it must cover $(2,5)$, in which 
case it completes a $R(2,3)$ with the tromino covering $(1,5)$. Now consider the two possible orientations of the 
tromino covering $(5,5)$ and $(5,6)$. If it covers $(4,6)$ then a $R(3,2)$ is completed by the tromino covering 
$(4,5)$. If it covers $(4,5)$ then a $R(3,2)$ is completed again by the trominon covering $(3,5)$. In both these 
cases the $7th$ vertical grid line becomes a fault line. So our claim holds. Note tha the tromino covering $(5,5)$ 
and $(5,6)$ must cover $(4,5)$ for similar reasons as above (otherwise the $7th$ vertical grid line will become 
a fault line). Consider the tromino covering $(2,5)$ and $(3,5)$. If it covers $(2,6)$, then the trominoes 
covering $(1,6)$ and $(1,8)$ will make $(3,9)$ inaccessible. So this tromino must cover $(3,6)$. Th remaining area 
can now only be covered in two ways, as shown in Figure 16(f). Now suppose that the tromino covering $(3,3)$ 
covers $(3,4)$ and $(2,4)$. It can be easily seen that the tromino covering $(1,4)$ must cover $(1,5)$ and 
$(2,5)$. Consider the tromino covering $(4,3)$ and $(5,3)$. If it covers $(4,4)$, then a $R(2,3)$ is completed 
by the tromino covering $(5,4)$. Now the tromino covering $(3,5)$ and $(3,6)$ must cover $(2,6)$, else if it 
covers $(4,6)$ then $(5,6)$ will become inaccessible. The reader can now easily see that the trominoes covering 
$(1,6)$ and $(1,8)$ will make $(3,9)$ inaccessible. So we conclude that the tromino covering $(4,3)$ and $(5,3)$ 
must cover $(5,4)$. The tromino covering $(4,4)$ must cover $(4,5)$ and $(3,5)$ for similar reasons as above. 
Now the tromino covering $(5,5)$ must cover $(5,6)$ and $(4,6)$. Also, the tromino covering $(4,7)$ and $(3,7)$ 
must cover $(3,6)$. If it covers $(3,8)$ instead, then the tromino covering $(3,6)$ will make $(1,6)$ inaccessible.
Now the remaining area can be covered only in two ways, as shown in Figure 16(g). 

Finally, consider the case when a tiling of the leftmost three columns of a generator extending the pattern in 
Figure 2(d) matches the flipped counterpart of the pattern in Figure 3(c). The tromino covering $(3,3)$ has three 
permissible orientations. If it covers $(3,4)$ and $(2,4)$ then a $R(3,2)$ is completed by the tromino covering 
$(2,3)$. If it covers $(3,4)$ and $(2,3)$, then a $R(3,2)$ is again completed by the tromino covering $(1,3)$. 
Suppose we have a $R(3,2)$ as per the above two cases. Consider the tromino covering $(4,4)$ and $(5,4)$. If it 
covers $(4,5)$, then a $R(2,3)$ is completed by the tromino covering $(5,5)$. In this case, the tromino covering 
$(1,5)$ has two possible orientations. If it covers $(1,6)$ and $(2,5)$ then a $R(3,2)$ is completed by the tromino
covering $(3,5)$, and if it covers $(1,6)$ and $(2,6)$, then a $R(3,2)$ is again completed by the tromino covering 
$(2,5)$. In both these cases, the $7th$ vertical grid line becomes a fault line. So the tromino covering $(4,4)$ 
and $(5,4)$ must cover $(5,5)$. The tromino covering $(4,5)$ will not cover $(4,6)$ and $(5,6)$ for similar reasons 
as above. If it covers $(4,6)$ and $(3,6)$, then the tromino covering $(3,5)$ will make $(1,5)$ inaccessible. So 
this tromino must cover $(3,5)$ and $(3,6)$. It is easy to see that the tromino covering $(5,6)$ must cover $(4,6)$ 
and $(4,7)$. The remaining area can now only be tiled as shown in Figure 16(h). Now consider the case when the 
tromino covering $(3,3)$ covers $(3,4)$ and $(4,4)$. Note that the tromino covering $(5,4)$ must cover $(5,5)$ and 
$(4,5)$, and the tromino covering $(5,6)$ must cover $(4,6)$ and $(4,7)$. Consider the two possible orientations 
of the tromino covering $(1,3)$ and $(2,3)$. If it covers $(2,4)$ then a $R(2,3)$ is completed by the tromino 
covering $(1,4)$. In this case, it is easy to see that the rest of the area becomes untileable. So the tromino 
covering $(1,3)$ and $(2,3)$ must cover $(1,4)$ and the tromino covering $(2,4)$ must cover $(2,5)$ and $(3,5)$ 
for similar reasons as above. The remaining area now can be uniquely tiled as shown in Figure 16(i). So we 
conclude that Figure 16 shows all possible ways of extending the pattern in Figure 2(d) so that the tiling is 
faultfree.  \hfill \qed
\end{proof}

For the very first extension of $R(7,6)$ to $R(7,12)$ we use the generators in Figure 14, 15 and 16. The 
remaining two rows are extended by simply padding with $R(2,3)$ rectangles. Note that after the 
very first extension, we 
must have the pattern in Figure 3(c) or 3(e), or their flipped counterparts, in the topmost or bottommost five 
rows (in the leftmost (rightmost) four columns). 
We can now use Moore's approach to extend these patterns by $(6t-12)$ columns, padding the remaining 
two rows with $R(2,3)$ rectangles. From Figures 14 and 15, it can be seen that in 20 ways we can extend the pattern 
in Figure 2(c) (observing the fact that we retile the original $R(3,2)$ rectangle in this pattern),
to get a tiling of $R(7,12)$ where the tiling of the leftmost three columns matches the pattern 
in Figure 3(e); and in 32 ways to get a tiling of $R(7,12)$ where the 
tiling of the leftmost three columns matches the pattern in Figure 3(c). Similarly, From Figure 16, we can see 
that the pattern in Figure 2(d) can be extended in 20 ways to get a tiling of $R(7,12)$, where the tiling of the 
leftmost three columns matches the pattern in Figure 3(e); and in 36 ways to get a tiling of $R(7,12)$, where the 
tiling of the leftmost three columns matches the pattern in Figure 3(c). Since $R(7,6)$ (as 
in Figure 10(b)) can be tiled in 8 ways (see Lemma 9), and the remaining two rows are being padded 
with $R(2,3)$ rectangles, we can see that $8\times40\times4 = 
1280$ tilings of $R(7,12)$ contain the pattern in Figure 3(e) in their leftmost three columns. Now these tilings 
can be extended in $\frac{c_{2t-3}}{2}\times2^{2t-4}$ ways, where $t\geq2$. So the generating function $A(z)$, 
where the coefficient of $z^{2t}$, $t\geq2$ gives all ways of tiling $R(7,6t)$ faultfreely, when extension of 
$R(7,6)$ (as in Figure 10(b)) is done only on one side (either entirely to the left, or entirely to the 
right) by the above approach is,

\begin{eqnarray}
A(z) & = & \frac{1280}{4}\sum_{t\geq2}c_{2t-3}.2^{2t-3}.z^{2t} \nonumber \\
     & = & 160z^3.\{G_3'(2z) - G_3'(-2z)\}
\end{eqnarray} 

Similarly, we can see that $8\times\ 68\times4 = 2176$ tilings of $R(7,12)$ contain the pattern in Figure 
3(c) in their leftmost three columns. These tilings can be extended in $\frac{b_{2t-3}}{4}\times2^{2t-4}$ 
ways, where $t\geq2$. So the generating function $B(z)$, where the coefficient of $z^{2t}$, $t\geq2$ gives 
all possible ways of tiling $R(7,6t)$ faultfreely, when extension of $R(7,6)$ (as in Figure 10(b)) is done 
only on one side (either entirely to the left, or entirely to the right) by the above approach is,

\begin{eqnarray}
B(z) & = & \frac{2176}{8}\sum_{t\geq2}b_{2t-3}.2^{2t-3}.z^{2t} \nonumber \\
     & = & 136z^3.\{G_1'(2z) - G_1'(-2z)\}
\end{eqnarray}   

Note that all tilings produced by extending $R(7,6)$ (as in Figure 10(b)) by the approach will be distict 
from those tilings produced by extending $R(7,6)$ (as in Figure 10(a)). This is because $R(7,6)$ as shown 
in Figure 10(a), contains a $R(3,2)$ in the middle (in columns 3 and 4, from rows 3-5) 
which is left undisturbed in the first extension to 
$R(7,12)$. On the other hand, $R(7,6)$, as shown in Figure 10(b), contains no $R(3,2)$ in the middle and it 
can be easily seen that the first extension to $R(7,12)$ does not disturb both the two trominoes, one covering 
$(3,4)$, $(4,4)$, $(4,5)$ and the other covering $(4,2)$, $(4,3)$, $(5,3)$, at the same time. Only one out 
of these two trominoes is disturbed in the very first extension. So we do not retile in the middle to get a 
$R(3,2)$ in these extensions. So, all tilings produced by extending Figure 10(b) are distinct from those 
produced by extending Figure 10(a). 

The generating function $Q_2(z)$, where the coefficient of $z^{2t}$, where $t\geq2$, 
is the total number of ways of tiling $R(7,6t)$ faultfreely, when $R(7,6)$ (as shown in Figure 10(b)) is extended only on one side (either entirely to the left, or entirely to the right) is,

\begin{eqnarray}
Q_2(z) & = & 2A(z) + 2B(z)
\end{eqnarray}

So, the generating function $F(z)$, where the coefficient of $z^{2t}$, where $t\geq2$, 
is the total number of ways of tiling $R(7,6t)$ faultfreely, when $R(7,6)$ is extended either only on one side (either entirely to the left, or entirely to the right) or on both sides is,

\begin{eqnarray}
F(z) & = & F_1(z) + F_2(z)
\end{eqnarray}

Here, $F_1(z)$ and $F_2(z)$ have been calculated in accordance with Theorem 7.

\section{Upper Bound on the Number of Tromino Tilings}

A natural question to ask is whether the number of domino tilings of a given rectangle $R(m,n)$ is more or less than the number of tromino tilings of $R(m,n)$. Intuitively,  one might say that we have more freedom in tiling a rectangle with dominoes than with trominoes, by this we mean that there are possibly more ways of retiling a given area with dominoes than with trominoes.  However, no such comparisons have been reported so far in the literature. We now present such a comparison, thereby producing an upper bound on the number of tromino tilings of a given rectangle $R(m,n)$, where $3|mn$ and $m,n>0$. 

In order to develop such a comparison, we first define a {\it monodic} tiling of $R(m,n)$ as a tiling with $\frac{mn}{3}$ dominoes and $\frac{mn}{3}$ monominoes. We now construct a one-to-one mapping function from the set of tromino tilings of a rectangle $R(m,n)$ to the set of monodic tilings of $R(m,n)$. Our mapping function is diagramatically explained in Figure 17. We call a domino with an arrow (as shown in Figure 17) a {\it directed} domino. We call a monodic tiling with directed dominoes as a {\it directed monodic tiling}. The conversion from a tromino tiling of $R(m,n)$ to a directed monodic tiling of $R(m,n)$ is clear from the mapping function (by converting every tromino to a combination of a directed domino and a monomino as shown in Figures 17(a)-(d)). We call a directed monodic tiling of $R(m,n)$ as {\it valid} if it was obtained from a tromino tiling of $R(m,n)$ via the mapping function depicted in Figure 17. For the reverse, converting a valid directed monodic tiling of $R(m,n)$ to a tromino tiling of $R(m,n)$, we simply attach every monomino to the right of the arrow of a directed domino, in a valid directed monodic tiling of $R(m,n)$, to form a tromino; thereby constructing a tromino tiling of $R(m,n)$. The proof that this mapping is one-to-one is easy and is left to the reader as an exercise.

\begin{figure}[htbp]
\centerline{
{\scalebox{.25}{\includegraphics{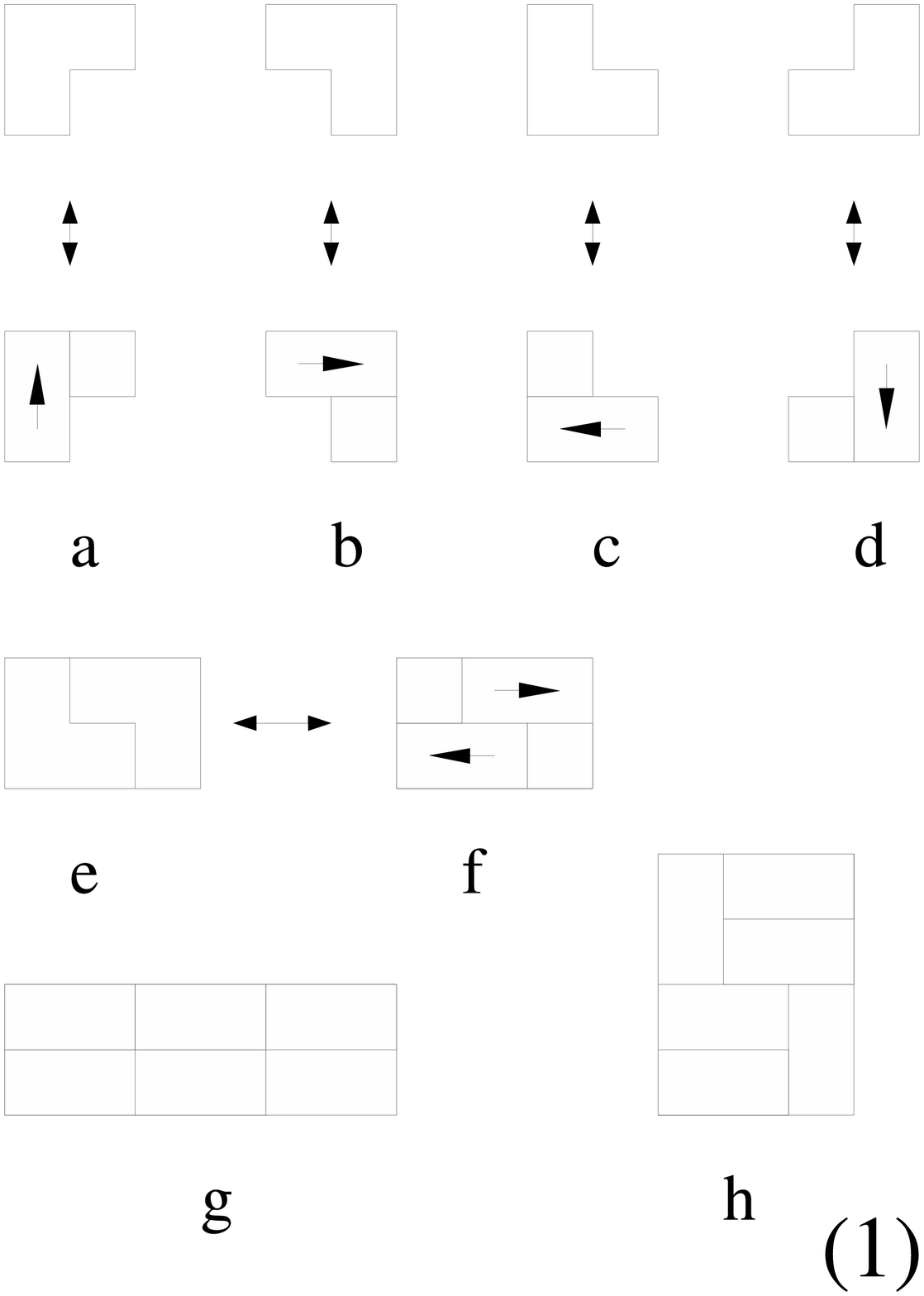}}}
{\scalebox{.23}{\includegraphics{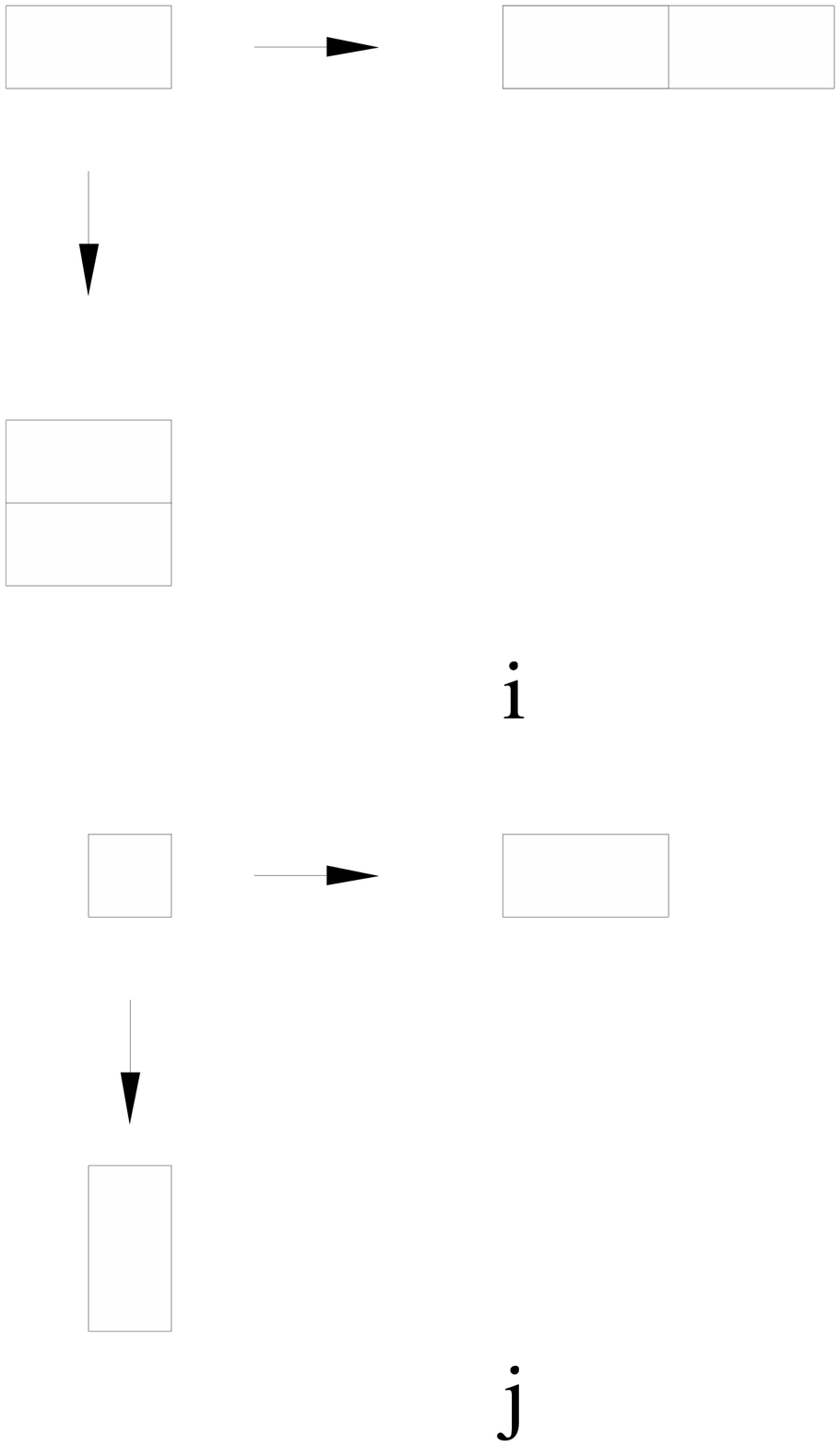}}}
{\scalebox{.15}{\includegraphics{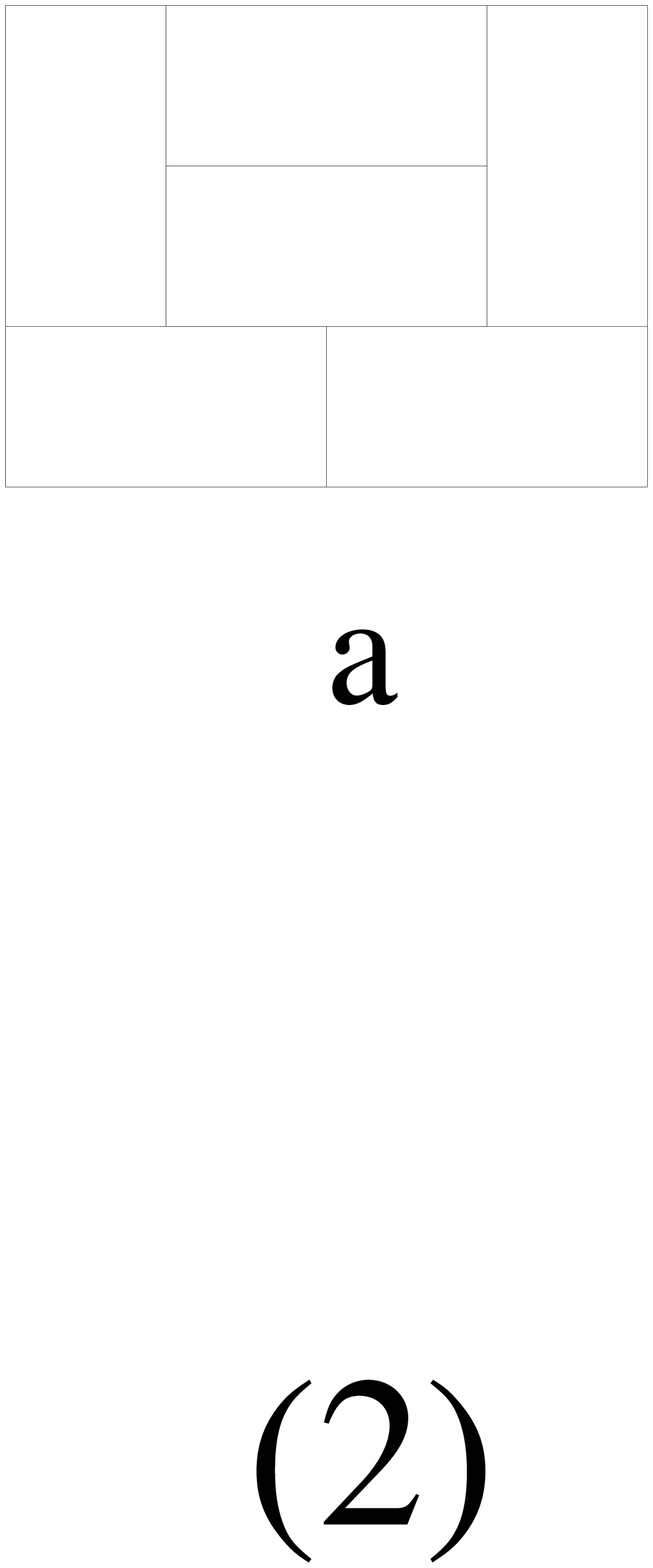}}}
}
\caption{(1) Mapping function from a tromino tiling to a domino tiling. 
(2) A domino tiling which is not obtained via the given mapping function.}
\end{figure}

As depicted in Figure 17(1), a tromino tiling of $R(m,n)$, after being converted to the corresponding directed monodic tiling, decided by the mapping function in Figure 17(1) (a)-(d), can be {\it stretched} either horizontally or vertically to a domino tiling. By {\it stretching} we double either the length $n$ or the breadth $m$ of the rectangle $R(m,n)$. In this process, every monomino gets converted to a domino, and a domino gets converted into either two horizontal dominoes lying side-by-side (in case of a horizontal stretching of $R(m,n)$), or it gets converted to two horizontal dominoes one lying on top of the other (in case of a vertical stretching of $R(m,n)$). This conversion is shown in Figures 17(1) (i)-(j). So a tromino tiling of $R(2,3)$ (as shown in Figure 17(1)(e)) gets converted to a domino tiling of either $R(2,6)$ or $R(4,3)$ (as shown in Figures 17(1)(g) and (h)). The reader must note that not every domino tiling of $R(m,2n)$ is obtained by stretching a monodic tiling of $R(m,n)$. Note, for example, the domino tiling of $R(3,4)$ shown in Figure 17(2). The reader can verify that this domino tiling cannot be obtained by stretching (either horizontally or vertically) any monodic tiling of $R(3,2)$. 

By our mapping function, every tromino tiling of $R(m,n)$ gets coverted to a unique directed monodic tiling. If we mutiply the total number of monodic tilings of $R(m,n)$ by $2^{\frac{mn}{3}}$, we get the total number of directed monodic tilings of $R(m,n)$. First, let us focus only on the case when the stretching done is horizontal. With slight modification, our streching process can also be converted to a one-to-one mapping function which converts every monodic tiling of $R(m,n)$ to a domino tiling of $R(m,2n)$. We call a domino tiling of $R(m,2n)$ {\it valid} if it was obtained by stretching a monodic tiling of $R(m,n)$. Any valid domino tiling can be converted back to a monodic tiling of $R(m,n)$ simply by {\it unstretching} (i.e., compressing lengthwise) the valid domino tiling of $R(m,2n)$. 


Suppose we have a monodic tiling of $R(m,n)$. We {\it colour} this tiling in the following manner, every monomino is coloured {\it blue} and every domino is coloured {\it red}. We now stretch $R(m,n)$ length-wise to produce $R(m,2n)$. The outcome will be a {\it coloured} domino tiling of $R(m,2n)$. Note that the dominoe(s) produced by the stretching of a coloured monomino (domino) will be of the same colour as that of the monomino (domino). And similarly, we define the unstretching of a coloured domino tiling. We have the following lemma.  

\begin{lemma}
Distinct coloured monodic tilings of $R(m,n)$ give rise to distinct coloured domino tilings of $R(m,2n)$.
\end{lemma}
\begin{proof}
We prove the above claim by the method of contradiction. Let us assume on the contrary that the above claim is false. First suppose that two distinct coloured monodic tilings of $R(m,n)$ give rise to the same coloured domino tiling of $R(m,2n)$ via stretching. Since the two monodic tilings of $R(m,n)$ being considered are distinct, we know that the position of at least one monomino is different in these two tilings. This means that the position of a monomino in the first tiling of $R(m,n)$ is covered by a domino in the second tiling of $R(m,n)$. Let this monomino be $(i,j)$. It can be easily seen that after stretching of $R(m,n)$, the squares $(i,2j)$ and $(i,2j+1)$ correspond to the square $(i,j)$ in $R(m,n)$. These two squares in $R(m,2n)$ will be covered by a blue domino via the first monodic tiling and a red domino via the second. But this is impossible because both the coloured monodic tilings of $R(m,n)$ being considered give rise to the same coloured domino tiling of $R(m,2n)$. Thus, we arrive at a contradiction. 
 
Now suppose that two distinct coloured domino tilings of $R(m,2n)$ give rise to the same coloured monodic tiling of $R(m,n)$ after unstretching. Again, since the two coloured domino tilings are distinct, the position of at least one square $(i,j)$ (say) in one domino tiling is covered by a red (blue) domino and which is covered by a blue (red) domino in the other. The reader can see that after unstretching the rectangle $R(m,2n)$, each square $(i,j)$ in $R(m,2n)$ maps to $(i,[j/2])$ in $R(m,n)$ (here [x] represents the integer part of x). When we unstretch the two coloured tilings, this means that the square $(i,[j/2])$ in $R(m,n)$, corresponding to the square $(i,j)$ in $R(m,2n)$, will be red (blue) via the first tiling and blue (red) via the second tiling. Thus, we again arrive at a contradiction. So, we conclude that distinct coloured monodic tilings of $R(m,n)$ map to distinct coloured domino tilings of $R(m,2n)$.   \hfill \qed
\end{proof}   

We have already shown that not all domino tilings of $R(m,2n)$ are valid. Thus, via the two mappings just described above, we have defined an injective function from the set of tromino tilings of $R(m,n)$ to a proper subset of coloured domino tilings of $R(m,2n)$. The number of coloured domino tilings of $R(m,2n)$ is simply $2^{mn}$ times the number of domino tilings of $R(m,2n)$ (since every domino can be coloured either red or blue). Since the direction of stretching is unimportant in the proofs above, the reader should convince himself that the above claims also hold for vertical stretching of rectangles. Let $N_T(m,n)$ $(N_D(m,n))$, denote the number of tromino (domino) tilings of $R(m,n)$. We summarize our result in the following theorem.

\begin{theorem}
For all rectangles $R(m,n)$, such that $3|mn$ and $m,n>0$, the following inequality holds:
\begin{eqnarray}
N_T(m,n)\leq 2^{\frac{4mn}{3}}\{min[N_D(m,2n),N_D(2m,n)]\}
\end{eqnarray}

where the number of domino tilings of $R(2m,2n)$ is given by the formula,

\begin{eqnarray}
N_D(2m,2n) & = & 4^{mn}\prod_{j=1}^{m}\prod_{k=1}^{n}\{\cos^2\frac{j\pi}{2m+1} + \cos^2\frac{k\pi}{2n+1}\}
\end{eqnarray}
\end{theorem}

\section{Final Remarks}

We are currently viewing the rules of our incremental generative 
scheme in the formal language theoretic grammar paradigm, and propose to formalize the relationship 
between the generating functions and the grammar rules. Combining these rules with inclusion-exclusion 
counting techniques and integer partitioning (related to {\it Sterling Numbers of the Second Kind}), we 
plan to derive lower (upper) bounds on the number of faultfree tromino tilings of $m\times n$ rectangles. 
We wish to study these bounds obtained asymptotically, and compare them with 
known and any new lower (upper) bounds. 





 

\end{document}